\documentclass[oneside]{report}
\usepackage{amsmath}

\usepackage[fontsize=11pt]{scrextend}
\usepackage{tgtermes} % times font
\usepackage{graphicx} % Required for inserting images
\graphicspath{ {./images/} }
\usepackage{blindtext}
\usepackage{subcaption}
\usepackage{mwe}
\usepackage{setspace} 
\usepackage{mathdots}
\usepackage{enumerate}
\usepackage{microtype}
\usepackage[shortlabels]{enumitem}
\usepackage{setspace}%Using \doublespacing in the preamble
\usepackage{geometry}
\usepackage{xcolor}
\usepackage{hyperref}

\hypersetup{
    colorlinks=true,
    linkcolor=blue,
    filecolor=magenta,      
    urlcolor=cyan,
    citecolor=red,
    pdftitle={Overleaf Example},
    pdfpagemode=FullScreen,
    pdftitle={},
pdfauthor={}
    }
\renewcommand{\arraystretch}{0.6}
\usepackage{mathrsfs}
\usepackage{amsthm}
\usepackage{sectsty}
\chapternumberfont{\Huge}
\newtheorem{theorem}{Theorem}[section]
\newtheorem{lemma}[theorem]{Lemma}
\newtheorem{proposition}[theorem]{Proposition}
\newtheorem{corollary}[theorem]{Corollary}
\theoremstyle{definition}
\newtheorem{definition}[theorem]{Definition}

\newtheorem{example}[theorem]{Example}

\newtheorem{remark}[theorem]{Remark}

\newtheorem{construction}[theorem]{Construction}
\newtheorem{question}[theorem]{Question}
\usepackage{amsfonts}
\usepackage{amssymb}
\usepackage{breqn}
\usepackage{mathrsfs}
\usepackage[mathscr]{}
\begin{document}
\begin{titlepage}
    \centering
    \vspace*{1cm}
    \LARGE
    \textbf{DYNAMICS OF TRANSVERSAL MATROIDS}

    \vspace{1.5cm}
    \LARGE
    \text{HIN LON LAO}

    \vfill
    \large
    A DISSERTATION SUBMITTED TO THE FACULTY OF GRADUATE STUDIES IN PARTIAL FULFILLMENT
    OF THE REQUIREMENTS FOR THE DEGREE OF DOCTOR OF PHILOSOPHY

    \vspace{0.8cm}
    GRADUATE PROGRAM IN MATHEMATICS AND STATISTICS
    \\YORK UNIVERSITY
    \\TORONTO, ONTARIO

    \vspace{2cm}
    May 2026

    \vspace{1cm}
    \copyright HIN LON LAO, 2026
\end{titlepage}

\chapter*{Abstract}
\addcontentsline{toc}{chapter}{Abstract}
\pagenumbering{roman}
\setcounter{page}{2}
We propose a way of generating matroids by randomly picking indicator vectors, and since we will see that it mostly gives uniform matroids, which are transversal, we then look at randomly generated transversal matroids. In both models we also investigate threshold phenomena of certain monotone increasing properties. Non-monotone properties of transversal matroids will also be investigated. Moreover, we are going to see how to spot circuits of a transversal matroid from its bipartite graph presentation. In the last Chapter, we have a survey of the characterization of transversal matroids that admit acyclic bipartite presentations. We will also investigate when we can simplify a presentation for a transversal matroid.

\chapter*{}
\addcontentsline{toc}{chapter}{Dedication}
\textit{In memory of my great father, Mr. Lao, Kwok Fai, who had no chance of seeing this. His sacrifices for me can hardly be described in words. His name will be forever remembered and enshrined with this dissertation.}

\chapter*{Acknowledgments}
\addcontentsline{toc}{chapter}{Acknowledgments}
I must express my deepest gratitude to my supervisor, Nantel Bergeron, who always guided me back to a positive mindset whenever I was overwhelmed by negative thoughts. Without his immense patience, my journey to complete this dissertation would have been a nightmare. I would also like to thank Ada Chan and Patrick Ingram for serving on my supervisory committee. Many thanks to Alex Fink for serving as the external examiner for my dissertation and defense, and for providing constructive feedback on this dissertation. In addition, I have to thank Nathan Smale, who has supported me immensely since my time as an undergraduate student at the University of Utah. I also thank Weng-Fong Ke and Ngau Lam for their invaluable guidance during my MSc studies at National Cheng Kung University.

Thanks to my mother, Kit Meng Cheong, who has managed many family matters since my father passed away, allowing me to concentrate on my aspirations. Furthermore, thanks to my friends, Shih-Jie Chang and Chi Hou Cheang, for their support when I was in low spirits.

Finally, I want to extend my thanks to all those whom I am unable to name individually, but who have helped me along this journey.

{
\hypersetup{hidelinks}
\clearpage
\phantomsection
\addcontentsline{toc}{chapter}{Table of Contents}
\tableofcontents

\clearpage
\phantomsection
\addcontentsline{toc}{chapter}{List of Figures}
\listoffigures
}

\chapter{The Use of Transversal Matroids in Linear Network Coding and Electric Network Theory}
\pagenumbering{arabic}
In this Chapter, we introduce the motivation for why we choose to investigate transversal matroids and randomly generated matroids in this dissertation. A matroid $M$ on a finite set $E$ is a pair $(E,\mathscr{I})$, where $\mathscr{I}$ is a collection of subsets of $E$, such that 
\begin{enumerate}
    \item $\emptyset \in \mathscr{I}$.
    \item If $I\in \mathscr{I}$, and $I'\subseteq I$, then $I'\in \mathscr{I}$.
    \item If $\mathcal{I}_1$ and $\mathcal{I}_2$ are in $\mathscr{I}$ with $|\mathcal{I}_1|<|\mathcal{I}_2|$, then there is an element $e$ of $\mathcal{I}_2-\mathcal{I}_1$ such that $\mathcal{I}_1\cup e\in \mathcal{I}$.
\end{enumerate}
Members of $\mathscr{I}$ are called independent sets of $M$. A subset of $E$ that is not independent is called dependent in $M$. A circuit of $M$ is a minimal dependent set of $M$.

Since the founding paper of Whitney \cite{Whitney} on matroid theory, the subject has undergone significant research and there is an abundance of related mathematical objects or techniques. For instance, rigid structures, Tutte polynomials, hyperplane arrangements etc. The series by White \cite{White,White2,White3} is recommended for readers who are interested in a deep dive into the rich applications of the subject. In this dissertation, we are mostly concerned with a special class of matroids called transversal matroids  discovered by Edmonds and Fulkerson \cite{Edmonds}. Let $I$ be a finite index set and $\mathscr{A}=(A_i:i\in I)$ be a collection of subsets (they need not be distinct) of a finite $E$. A \textit{transversal} or \textit{system of distinct representatives} of $\mathscr{A}$ is a subset $\{e_i:e_i\in A_i,\forall i\in I\}$ and $e_i\neq e_j$ if $i\neq j$. Hall \cite{Hall} wrote down a necessary and sufficient condition for $\mathscr{A}$ to have a transversal: $\mathscr{A}$ has a transversal if and only if $|\bigcup_{i\in J}A_i|\;\geq |J|$ for all $J\subseteq I$. A set is called a \textit{partial transversal} of $\mathscr{A}$ if it is a transversal of some subcollection of $\mathscr{A}$. Edmonds and Fulkerson showed in \cite{Edmonds} that the collection of partial transversals of $\mathscr{A}$ forms a collection of independent sets of a matroid on $E$. Such a matroid is called a \textit{transversal matroid} and $\mathscr{A}$ a \textit{presentation} of the transversal matroid. Many studies have been initiated on transversal matroids since then. In Section \ref{sec: Brualdi et al.}, we will see that there is a unique maximal presentation for each transversal matroid, and all presentations of a transversal matroid are characterized \cite{Brualdi}. A classical work of Ingleton and Piff \cite{Ingleton} showed that transversal matroids are precisely the duals of strict gammoids. We will also mention this fact in Section \ref{sec:T(G)^*}. It is worth pointing out that Ardila \cite{Ardila} reproved this by using representations of transversal matroids and strict gammoids (when there is a bijective map between the ground set of the matroid and the columns of a matrix over a field such that the map preserves linear independence, the matroid is said to be \textit{representable} over that field and the matrix is called a \textit{representation} of the matroid). Mason \cite{Mason} characterized transversal matroids using their cyclic sets (sets that are unions of circuits of the matroid), and Ingleton \cite{Ingleton1} showed that Mason's results can be streamlined to cyclic flats. Brylawski \cite{Brylawski} characterizes transversal matroids by using their geometric presentations. Although the author is not going to introduce geometric presentations of matroids here because they are outside the scope of this dissertation, interested readers may consult \cite[Section 1.5]{Oxley}, \cite[Section 2.5]{Bonin1} and \cite{Brylawski} for details. Moreover, Piff and Welsh \cite{Piff} showed that the transversal matroids are representable over a sufficiently large field. The lecture note by Bonin \cite{Bonin1} serves as a good comprehensive introduction to transversal matroids.

 In what follows, we look at some applications of transversal matroid theory. Transversal matroids play an important role in the theory of network coding. The motivation and background of linear network coding can be found in Li et al. \cite{Li}. To introduce the relations between transversal matroids and linear network code, we briefly look at the content of Sun et al. \cite{Sun} and Eslami et al. \cite{Eslami} here. 

A network is a directed multi-edge graph $\mathscr{N}=(V,E,s)$, with set of vertices $V$ and a set of edges $E$, and a unique source vertex $s$ (the vertex with no incoming edges). An edge in $E$ is an ordered pair $(a,b)$ of vertices of $V$; $a$ and $b$ are called the \textit{tail} and the \textit{head} of the edge, respectively. So $s$ is the unique vertex such that there is no edge whose head is $s$. Let $h$ denote the number of edges whose tails are $s$.  An ordered pair of edges $(e_i,e_j)$ is called an \textit{adjacent pair} if the head of $e_i$ equals the tail of $e_j$. A path from $s$ to a subset $E'\subseteq E$ is any sequence of edges $(e_1,...,e_k)$ such that $(e_i,e_{i+1}),1\leq i\leq k-1$ are adjacent pairs, the tail of $e_1$ is $s$, and $e_k\in E'.$ Let $\operatorname{maxflow}(E')$ denote the maximum number of edge-disjoint paths from $s$ to $E'$. Thus, the maximum value $\{\operatorname{maxflow}(E')\mid E'\subseteq E\}$ can take is $h$ (It can be achieved by taking $E'$ to be the set $\{(s,v)\in E\mid v\in V\}$). Thus, $h$ is often called the \textit{maxflow} of the network $\mathscr{N}$. Let $In(v)$ denote the set of edges whose heads are $v$, and $Out(v)$ denote the set of edges whose tails are $v$. For a given network, we assign a $k_{i,j}\in \mathbb{F}$ to every pair $(e_i,e_j)$ of edges, and $k_{i,j}=0$ when $(e_i,e_j)$ is not an adjacent pair. An element $k_{i,j}$ of $(k_{i,j})$ is called the \textit{encoding coefficient} for $(e_i,e_j)$. Furthermore, for each $e_i\in E$, there is a \textit{coding vector} $f_i\in \mathbb{F}^h$ assigned such that: 
\begin{enumerate}
    \item $\{f_i:e_i\in Out(s)\}$ forms the natural basis of $\mathbb{F}^h$;
    \item $f_j=\sum_{e_i\in In(v)}k_{i,j}f_i$ for every edge $e_j\in Out(v)$ and non-source $v\in V$.
\end{enumerate}
The assignment of coding vectors that satisfy the stated conditions above is called a \textit{linear network code} for the network. For $E'\subseteq E$, let $M_{E'}$ be a matrix whose columns are the coding vectors of the edges in $E'$ (There are many choices for $M_{E'}$, but we only need to arbitrarily choose one of them. Since permuting the columns of $M_{E'}$ does not affect the rank of $M_{E'}$, and that suffices for what we need below). Then a \textit{linear network code} on a given network $(V,E,s)$ is said to be \textit{generic} (abbreviated as \textit{GLN code}) if for every $E'\subseteq E$ with at most $h$ edges, rank$(M_{E'})$=maxflow($E'$).

\begin{figure}[h]
    \centering
    \includegraphics[width=0.56\linewidth]{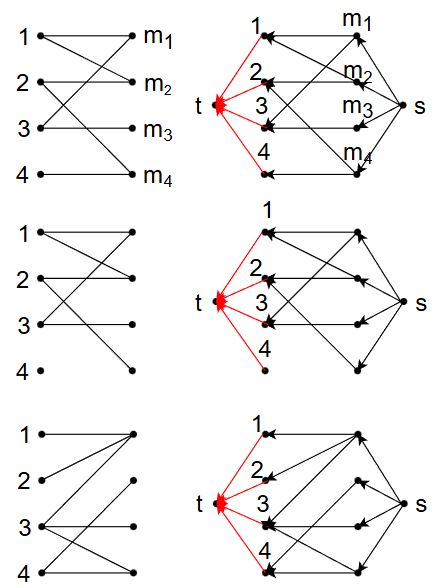}
    \caption{Transversal matroids as restriction of network matroids.}
    \label{fig:simple net}
\end{figure}

Network matroids are a class of representable matroids introduced by Sun et al. in \cite{Sun}. For a network $\mathscr{N}=(V,E,s)$, let $I$ be the collection of subsets $E'\subseteq E$ of cardinality at most $h$ with maxflow$(E')=|E'|$. Then $(E,I)$ forms a matroid which is called the \textit{network matroid} induced by $\mathscr{N}$. Moreover, there is a unique bipartite graph $G_{\mathscr{N}}$ such that its corresponding transversal matroid is dual to the network matroid of $\mathscr{N}$ \cite[Theorem 13]{Sun}. The way to construct $G_{\mathscr{N}}$ from $\mathscr{N}$ can be found in \cite[Section IV]{Sun} or \cite[Section III]{Eslami}. Moreover, in \cite[Section IV]{Eslami}, Eslami et al. discovered a way of constructing a GLN code on a network $\mathscr{N}$ by constructing a representation matrix for the transversal matroid induced by $G_{\mathscr{N}}$. It is worth pointing out that they have an algorithm to construct a representation for a transversal matroid \cite{Eslami1} such that it is faster than that in \cite{Piff}.  Indeed, if we can find a simpler presentation for this matroid (that is, a simpler bipartite graph inducing the same matroid), the process can be simplified. Due to Ingleton and Piff \cite{Ingleton}, a transversal matroid is isomorphic to a restriction of a network matroid (network matroids are in fact strict gammoids \cite[Theorem 3]{Sun}). In \cite[Section II]{Sun}, Sun et al. also restated that fact via direct graphical illustration.  We give an illustrative example here to tell their idea. In the top diagram of Figure \ref{fig:simple net}, the given bipartite graph represents a collection of subsets of $\{1,2,3,4\}$. Namely, it is $\mathscr{A}=\{\{1,3\},\{1,2\},\{3\},\{2,4\}\}$, the partial transversals of which form a collection of independent sets of a matroid on $\{1,2,3,4\}$. A set is a partial transversal of $\mathscr{A}$ if and only if it has a matching in the given bipartite graph. Then we turn the graph into a directed graph by directing all given edges towards elements of $\{1,2,3,4\}$, and we add vertices $s$ and $t$ to the graph, as well as directed edges $(s,m_i)$ and $(i,t)$ for all $1\leq i\leq 4$. The directed graph we get is a network. A subset $X\subseteq \{1,2,3,4\}$ is a partial transversal of $\mathscr{A}$ if and only if $\{(x,t): x\in X\}$ is independent in the network matroid. Thus, the transversal matroid induced by the bipartite graph is the same as the network matroid restricted to the edges colored in red.

The diagrams in the middle and at the bottom of Figure \ref{fig:simple net} illustrate why we are concerned with simple transversal matroids in Section \ref{sec:simple T(G)}. In the middle one, $\{4\}$ is an isolated vertex in the bipartite graph. Therefore, there is no path from the source $s$ to vertex 4 in the corresponding network. In the bottom one, $\{1,2\}$ is a pair of parallel edges (a two-element circuit) in the transversal matroid induced by the given bipartite graph. In the corresponding network, the directed edges oriented towards $\{1\}$ and $\{2\}$ come from exactly one vertex. Then the network is more efficient if vertex $\{1\}$ is identified with $\{2\}$.

A matroid on a set can be determined from its set of circuits; the set of circuits then becomes a significant object to investigate. In Proposition \ref{prop: cir. N=M+1}, Lemma \ref{lemma: cir. u. of cirs. iff no coloop} of this dissertation, we see how to spot the circuits of a transversal matroid from its presentation. Our main contributions include the following: In Theorem \ref{thm:type Z Thm}, we have a sufficient condition on a bipartite graph $G$ to spot a circuit in the transversal matroid induced by $G$. In Theorem \ref{thm:circuit, planting of pure basic type}, we characterize what a circuit of a transversal matroid looks like in an acyclic presentation of the transversal matroid. In Theorem \ref{thm:admits acyclic pre.}, we give a sufficient condition for two circuits of a transversal matroid $M$ such that we can find a presentation for $M$ in which the union of the two circuits is acyclic.

Now we look at another application of transversal matroids. The solvability of active linear electric networks (a different object from the network mentioned above) is an important problem of computer aided electric network design. Papers such as \cite{Iri,Iri1,Ozawa,Petersen,Recski,Recski1} were published to give necessary and sufficient conditions for the unique solvability, mainly using matroid theory, while they refer to different classes of electric networks. Recski and Iri gave a unified description of all models used in those papers in \cite{Recski2}, in which they also characterize the hierarchy classes of electric networks described in those papers using transversal matroids. Specifically, an electric network is an interconnection of various devices. In \cite{Recski2}, the only types of devices considered are called \emph{$n$-ports}, which are sets of admissible pairs $(\mathbf{u},\mathbf{i})$ subject to a system of equations of the form:
$$\mathbf{A_1u}+\mathbf{A_2i}=0$$
where $\mathbf{A_1}$ and $\mathbf{A_2}$ are $n\times n$ matrices with $r[\mathbf{A_1|A_2}]=n$, $\mathbf{u}=(u_1,u_2,...,u_n)^T,i=(i_1,i_2,...,i_n)^T$. For a matrix $[\mathbf{A_1|\mathbf{A_2}}]$ associated to an $n$-port $N$, there is a matroid $\mathbf{M}_N$ on the set of columns of $[\mathbf{A_1}|\mathbf{A_2}]$, whose collection of linearly independent columns is the collection of independent sets. Properties of $N$ can be deduced by investigating $\mathbf{M}_N$ \cite[Theorem 4]{Recski2}. If $\mathbf{M}_N$ is transversal, it is important to know whether $\mathbf{M}_N$ is fundamental transversal (A transversal matroid is called \textit{fundamental transversal} if it possesses a presentation $\mathscr{A}=\{A_1,...,A_r\}$ such that $A_i\setminus \bigcup_{j\neq i} A_j\neq \emptyset$ for all $1\leq i\leq r$).  Then certainly it is helpful to find a simpler presentation for the matroid. This is what our Chapter 5 aims to do. And in \cite[Remark 4]{Recski2}, the author provides a very good summary of how to put a collection of matrices of the above form together to form a big matrix $\mathbf{A}$ to describe an electric network with many $n$-ports. Then the author defines a transversal matroid out of that and gives a sufficient and necessary condition for the unique solvability of electric network, see Section 2, and Statement $i$ in \cite{Recski2}. This again lays out the importance of investigating possible ways of finding a simpler presentation for a transversal matroid. Furthermore, adding edges to a presentation of a transversal matroid may lead to changes in the electric network. Our Chapter 4 and a significant part of our Chapter 3 are devoted to investigating how adding an edge to a bipartite graph affects properties of the transversal matroid it induces. It is worth noting that Recski has written an excellent monograph \cite{Recski3} on the topic of electric network theory and its relation to matroid theory. It is recommended to readers who want to know more about the subject.

Significant parts of our Chapter 3 are devoted to ways of randomly generating matroids. The inspiration comes from random graph theory due to P. Erdős and A. Rényi \cite{Erdos1,Erdos2,Erdos3,Erdos4}. The theory originated from the probabilistic technique applied by Erdős to compute a lower bound of Ramsey numbers \cite{Erdos}. This also opens up another branch of mathematics called probabilistic method, whose essence is that, to show the existence of a structure with certain desired properties, one shows that the desired properties hold in an appropriate probability space with positive probability. Alon and Spencer wrote a classic text \cite{Alon} on probabilistic method. Our model in Chapter 3 is about picking indicator vectors randomly and generating a matroid from them. Since in most cases we get uniform matroids, which are transversal, from that, we naturally have a transition to Chapter 4, which is mostly about randomly generated transversal matroids and the threshold phenomena of monotone increasing properties (threshold phenomena in random graph theory were discovered by Erdős and Rényi in \cite{Erdos1,Erdos2,Erdos3}). It is worth pointing out that it is natural to want to compute how many uniform matroids there are in the class of transversal matroids. That is our Subsection \ref{subsec:no. of mini. pre.}.

\chapter{Preliminaries}
In this chapter, we are going to navigate through concepts that we will be frequently using throughout the dissertation.
\section{Preliminaries in Matroid}
The theory of matroids begins in the work of Whitney \cite{Whitney}, and it can be defined in many equivalent ways. The concepts and results in this section are foundational in matroid theory, and their proofs can easily be found in classical references such as \cite[Chapter 1,2]{Oxley} and \cite[Chapter 1]{Welsh}. Here, we begin with the bases approach:

\begin{definition}
A matroid $M$ on a finite set (ground set) $E$ is a collection $\mathscr{I}$ of subsets of $E$ satisfying the following conditions:
\begin{enumerate}[start=0,label={(\bfseries I\arabic*)}]
\item $\emptyset \in \mathcal{I}$.
\item If $\mathcal{I}\in \mathscr{I}$ and $\mathcal{I'}\subseteq \mathcal{I}$, then $\mathcal{I'}\in \mathscr{I}$.
\item If  $\mathcal{I}_1$ and $\mathcal{I}_2$ are in $\mathscr{I}$ with $|\mathcal{I}_1|<|\mathcal{I}_2|$, then there is an element $e$ of $\mathcal{I}_2-\mathcal{I}_1$ such that $\mathcal{I}_1\cup e\in \mathcal{I}$.
\end{enumerate}
\end{definition}
    
Throughout this dissertation, whenever we mention a matroid, we mean a matroid on a finite ground set. The members of $\mathscr{I}$ are called the \textit{independent sets} of $M$. A subset of $E$ which is not independent is called \textit{dependent}. A maximal independent set is called a \textit{basis}. A minimal dependent set is called a \textit{circuit} of $M$. The ground set, collection of independent sets, collection of bases and circuits of $M$ are denoted by $E(M),\mathscr{I}(M), \mathscr{B}(M),\mathscr{C}(M)$ respectively. A dependent set which is a singleton is called a \textit{loop}. Note that from (I2), one easily deduces that all bases of a matroid share the same cardinality. The following theorem is a common and convenient characterization of a collection of bases of a matroid (see \cite[Theorem 1.2.3]{Oxley} for a proof). 
  
\begin{theorem}\label{BE}
Let $\mathcal{B}$ be a class of subsets of a finite set $E$. Then $\mathcal{B}$ is the collection of bases of a matroid on $E$ if and only if it satisfies the following properties:
\begin{enumerate}[start=1,label={(\bfseries B\arabic*)}]
\item $\mathcal{B}$ is non-empty.
\item If $\mathcal{B}_1$ and $\mathcal{B}_2$ are in $\mathcal{B}$ and $x\in \mathcal{B}_1-\mathcal{B}_2$, then there exists an element $x'\in \mathcal{B}_2-\mathcal{B}_1$ such that $(\mathcal{B}_1-x)\cup x'\in \mathcal{B}$.
\end{enumerate}
\end{theorem}

A typical example is a uniform matroid: Let $E$ be an $n$-element set and $\mathcal{B}$ be the collection of all $m$-element subsets of $E$. It is easy to see that $\mathcal{B}$ forms a collection of bases of a matroid on the set $E$. This is called the \textit{uniform matroid} on an $n$-element set of rank $m$, denoted by $U_{m,n}$. (Note that technically, if we have another $n$-element set $E'$ such that $E\neq E'$, $m$-element uniform matroids on them are different because their ground sets are different. But they are clearly \textit{isomorphic}: there is a bijection from $E$ to $E'$ that preserves independence. Thus, people often neglect the difference and call it 'the' uniform matroid on an $n$-element set of rank $m$).

It is also common to characterize a matroid by its collection of circuits. The following theorem (see \cite[Corollary 1.1.5]{Oxley} for reference), gives a description of $\mathscr{C}(M)$ for a matroid $M$:
\begin{theorem}
Let $\mathscr{C}$ be a class of subsets of a finite set $E$. Then $\mathscr{C}$ is the collection of circuits of a matroid on $E$ if and only if it satisfies the following properties: 
\begin{enumerate}[start=1,label={(\bfseries C\arabic*)}]
\item $\emptyset \notin \mathscr{C}$.
\item If $\mathcal{C}_1$ and $\mathcal{C}_2$ are members of $\mathscr{C}$ and $\mathcal{C}_1\subseteq \mathcal{C}_2$, then $\mathcal{C}_1=\mathcal{C}_2$.
\item If $\mathcal{C}_1$ and $\mathcal{C}_2$ are distinct members of $\mathscr{C}$ and $e\in\mathcal{C}_1\cap\mathcal{C}_2$, then there is a member $\mathcal{C}_3$ of $\mathscr{C}$ such that $\mathcal{C}_3\subseteq (\mathcal{C}_1\cup\mathcal{C}_2)-e$.
\end{enumerate}
\end{theorem}

\begin{definition}
    The rank function $\rho:2^E\to \mathbb{Z^+}\cup\{0\}$ of a matroid $M$ on a finite set $E$ is defined by $\rho(X):=$the cardinality of a maximal independent set contained in $X$, $X\subseteq E$. $\rho (E)$ is called the rank of $M$, which is also denoted by $\rho(M)$.
\end{definition}

\begin{theorem}\label{thm:rank function}
The rank function $\rho$ of a matroid $M$ on a finite set $E$ satisfies the following properties:
\begin{enumerate}[start=1,label={(\bfseries R\arabic*)}]
\item $0\leq\rho(X)\leq|X|$ for all $X\subseteq E$.
\item $\rho(X)\leq\rho(Y)$ for all $X\subseteq  Y\subseteq E$.
\item if $X,Y\subseteq E$, then $$\rho(X\cup Y)+\rho(X\cap Y)\leq \rho(X)+\rho(Y).$$
\end{enumerate}
\end{theorem}

\begin{definition}\label{def:flat}
    For a matroid $M$, a set $F\subseteq E(M)$ is called a flat of $M$ if $\rho(F\cup x)=\rho(F)+1$ for all $x\in E(M)\setminus F$. A flat $F$ is called a hyperplane of $M$ if $\rho(F)=\rho(M)-1$.
\end{definition}

\section{Dual, graphic, transversal matroid}
The concept of dual matroid originated from Whitney \cite{Whitney}, and Minty used another approach called \textit{graphoids} in \cite{Minty} to study matroids and their duals. 
\begin{theorem}[{\cite{Whitney}}]
    If $\{B_i:i\in I\}$ is the set of bases of a matroid $M$ on a set $E$, then $\{E\setminus B_i:i\in I\}$ is the set of bases of a matroid $M^*$ on $E$.
\end{theorem}
This matroid is called the \textit{dual matroid} of $M$. A loop of $M^*$ is called a \textit{coloop} of $M$. A circuit of $M^*$ is called a \textit{cocircuit} of $M$.

\begin{theorem}[{\cite{Whitney}}]
    If $G$ is an undirected finite graph, not necessarily simple, then the cycles of $G$ are the circuits of a matroid $M(G)$ on the edge set $E(G)$ of $G$.
\end{theorem}
The matroid $M(G)$ is called the \textit{cycle matroid} of $G$. This matroid is called the \textit{polygon matroid} of $G$ in Tutte's paper \cite{Tutte}.

\begin{theorem}
    If $M$ is a matroid on a set $E$ and $X\subseteq E$, then $$\mathscr{I}(M|X)=\{I:I
    \subseteq X, I\in \mathscr{I}(M)\}$$ is the class of independent sets of a matroid on $X$. 
\end{theorem}
This matroid is called the \textit{restriction} of $M$ to $X$ and denoted by $M|X$ or $M\setminus E-X$.

\begin{theorem}
    If $M$ is a matroid on a set $E$ and $X\subseteq E$, then $$\mathscr{I}(M/X)=\{I'\subseteq E-X:I'\cup B\in \mathscr{I}(M) \mbox{ for some }B\in \mathscr{B}(M|X)\}$$ is the class of independent sets of a matroid on $E-X$.
\end{theorem}
This matroid is called the \textit{contraction} of $X$ from $M$ and denoted by $M/X$.

The operations of deletion and contraction on a graphic matroid $M(G)$ have nice correspondences in $G$. Namely, $M(G)|X=M(G\setminus E-X)$ and $M(G)/X=M(G/X)$ \cite[Prop.3.2.1]{Oxley}, where $X$ is any subset of the edge set $E$ of $G$. Figure \ref{fig:deletion, contraction} is a simple illustration of a case of deletion. 

\begin{figure}[h]
\begin{subfigure}{0.5\textwidth}
\includegraphics[width=0.7\linewidth, height=5cm]{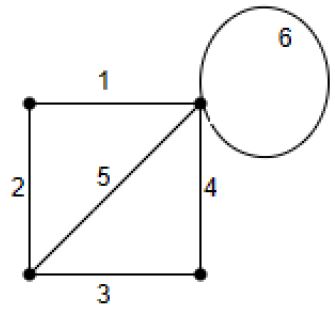} 
\caption*{$G$}
\end{subfigure}
\begin{subfigure}{0.5\textwidth}
\includegraphics[width=0.7\linewidth, height=5cm]{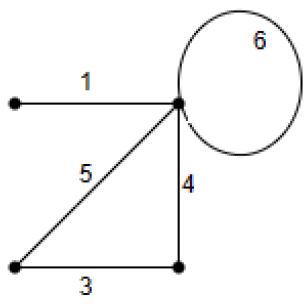}
\caption*{$G\setminus 2$}
\end{subfigure}
\begin{subfigure}{0.5 \textwidth}
    
    \includegraphics[width=0.9\linewidth]{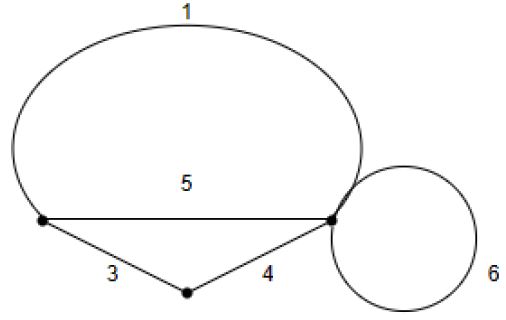}
    \caption*{$G/2$}
    \end{subfigure}
\caption{Operations of deletion and contraction on a graph.}
\label{fig:deletion, contraction}
\end{figure}

\newpage
\noindent$\mathscr{B}(M(G))=\{\{1,2,3\},\{1,2,4\},\{1,3,4\},\{1,3,5\},\{1,4,5\},\{2,3,4\},\{2,3,5\},\{2,4,5\}\}$, $\mathscr{B}(M(G\setminus2))=\{\{1,3,4\},\{1,4,5\},\{1,3,5\}\}$.  
\\ \noindent $\mathscr{B}(M(G/2))=\{\{1,3\},\{1,4\},\{3,4\},\{3,5\},\{4,5\}\}$.

Note that the operations of deletion and contraction on a matroid are commutative  \cite[Prop. 3.1.25]{Oxley}. A matroid $N$ obtained by performing a sequence of these two operations on a matroid $M$ is called a \textit{minor} of $M$.

\vspace{6mm}
Let $I$ be a finite index set and $\mathscr{A}=(A_i:i\in I)$ be a collection of subsets (they need not be distinct) of a finite $E$. A \textit{transversal} or \textit{system of distinct representatives} of $\mathscr{A}$ is a subset of $E$, $\{e_i:e_i\in A_i,\forall i\in I\}$ and $e_i\neq e_j$ if $i\neq j$. A set is called a \textit{partial transversal} of $\mathscr{A}$ if it is a transversal of some subcollection of $\mathscr{A}$.

\begin{theorem}[{\cite[Edmonds and Fulkerson]{Edmonds}}]
    Let $\mathscr{A}$ be a collection of subsets of a finite set $E$ and $\mathcal{I}$ be the set of partial transversals of $\mathscr{A}$. Then $\mathcal{I}$ is the collection of independent sets of a matroid on $E$.
\end{theorem}

\begin{definition} \label{def: transversal matroids}
    Let $\mathscr{A}$ be a collection of subsets of a finite set $E$. The matroid obtained from the set of partial transversals of $\mathscr{A}$ is denoted by $M(\mathscr{A})$. Such a matroid is called a transversal matroid and $\mathscr{A}$ a presentation of $M(\mathscr{A})$. Furthermore, a presentation $\mathscr{A}$ of a transversal matroid $M$ is a minimal presentation if deleting any member of $\mathscr{A}$, or deleting any element of any member results in a presentation of a different matroid. 
\end{definition}

\subsection{Transversal matroids and bipartite graphs} \label{G_N,G_M}

Given a (finite) bipartite graph, we may consider a transversal matroid induced by the bipartite graph. But there is a possible ambiguity. Given a (finite) bipartite graph $G$ with vertex bipartition $G_V=V_1\cup V_2$, let $A_v$ be the set of vertices of $G$ adjacent to $v\in G_V$. The collection $\{A_v:v\in V_1\}$ is a collection of subsets of $V_2$. If we consider the collection of partial transversals of $\{A_v:v\in V_1\}$, we get a matroid on $V_2$. Similarly, if we consider the collection of partial transversals of $\{A_v:v\in V_2\}$ we get a matroid on $V_1$. To bypass this ambiguity, we define a \textit{pointed bipartite graph} to be a (finite) bipartite graph $G$ with vertex bipartition $G_V=G_N\cup G_M$ such that when we mention the transversal matroid induced by $G$ (denoted by $T(G)$), we mean the matroid on $G_N$ with $\mathscr{I}(T(G))$ equal to the collection of partial transversals of $\{A_v:v\in G_M\}$.  

Throughout this dissertation, whenever we mention a bipartite graph, we mean a pointed bipartite graph.

Given a collection $\mathscr{A}=\{A_1,...,A_n\}$ of subsets of a finite set $E$, we may use a bipartite graph $G_{\mathscr{A}}$ to represent $\mathscr{A}$ as follows: Let the set of vertices of $G_{\mathscr{A}}$ be $(G_{\mathscr{A}})_N\cup (G_{\mathscr{A}})_M$ where $(G_{\mathscr{A}})_N:=E$ and $(G_{\mathscr{A}})_M:=\{1,...,n\}$ such that $i\in (G_{\mathscr{A}})_M$ is adjacent to all elements in $A_i$. By what we have discussed, $T(G_{\mathscr{A}})$ is the matroid on $E$ and $\mathscr{I}(T(G_{\mathscr{A}}))$ is the collection of partial transversals of $\mathscr{A}=\{A_1,..,A_n\}.$

A \textit{matching} in a graph is a set of edges in the graph such that no two edges in the set share a common endpoint. Then given a collection $\mathscr{A}$ of subsets of a finite set $E$, let $G$ be the bipartite graph $G$ representing $\mathscr{A}$ as above. By construction, $X\subseteq G_N$ is a partial transversal of $\mathscr{A}$ ($X$ is independent in $T(G)$) if and only if there is a matching in $G$ such that every edge of it has one endpoint in $X$, and one endpoint in $G_M$. Denote the set of those endpoints in $G_M$ by $J$. Then we say that $X$ is matched to/onto $J$, or vice versa. 

With the convention of this subsection in mind, when we say a presentation of a transversal matroid $M$, we can mean a collection of subsets $\mathscr{A}$ of $E(M)$ such that $M(\mathscr{A})=M$, or a bipartite graph $G$ such that $T(G)=M$. It should be clear from the content which one we are referring to.

\section{Preliminaries in probability}
Readers can find the contents of this section in references such as Robert B. Ash \cite[Chapter 1,5]{Ash}, Donald L. Cohn \cite[Chapter 1,2,10]{Cohn}, H.L. Royden and P.M. Fitzpatrick \cite[Chapter 2,3]{Royden}.
\begin{definition}
    Let $\Omega$ be an arbitrary set. A collection $\mathscr{A}$ of subsets of $\Omega$ is a $\sigma$-algebra on $\Omega$ if it satisfies the following conditions:
    \begin{enumerate}
        \item $\Omega\in \mathscr{A}$.
        \item for each set $A\in \mathscr{A}$, $A^c\in \mathscr{A}$.
        \item for each infinite sequence $\{A_i\}$ of sets that belong to $\mathscr{A}$, $\bigcup\limits_{i=1}^{\infty}A_i\in \mathscr{A}$.
        \item for each infinite sequence $\{A_i\}$ of sets that belong to $\mathscr{A}$, $\bigcap\limits_{i=1}^{\infty}A_i\in \mathscr{A}$.
    \end{enumerate}
\end{definition}

\begin{definition}
Let $\mathscr{A}$ be a $\sigma$-algebra on a set $\Omega$. An extended real-valued function $\mu$ on $\mathscr{A}$ is a measure on $\mathscr{A}$ if $\mu(\emptyset)=0$ and it satisfies the countably additive condition:
\begin{equation*}        \mu\left(\bigcup\limits_{i}A_i\right)=\sum\limits_{i}\mu(A_i)
\end{equation*}
whenever $A_1,A_2,...$ is a finite or countably infinite collection of disjoint sets in $\mathscr{A}$. 
\end{definition}

\begin{definition}\label{def:Borel inter.}
    The smallest $\sigma$-algebra of subsets of $\mathbb{R}$ containing all intervals $(a,b]$ with $a,b\in \mathbb{R}$ (technically, the intersection of all such $\sigma$-algebras) is called the Borel $\sigma$-algebra, or Borel set of $\mathbb{R}$, denoted by $\mathscr{B}(\mathbb{R})$.
\end{definition}

\begin{remark}\label{remark:inter. type}
    The intervals in the previous definition can be replaced by other classes of intervals of $\mathbb{R}$, for instance, 
    \\The set of all open intervals,
    \\The set of all closed intervals,
    \\The set of all intervals $[a,b),a,b\in \mathbb{R}$,
    \\The set of all intervals $(a,\infty),a\in \mathbb{R}$,
    \\The set of all intervals $[a,\infty),a\in \mathbb{R}$,
    \\The set of all intervals $(-\infty,a),a\in \mathbb{R}$,
    \\The set of all intervals $(-\infty,a],a\in \mathbb{R}$.
    \\That means the smallest $\sigma$-algebra containing all intervals of any particular kind mentioned above still gives $\mathscr{B}(\mathbb{R})$. And as a result, $\mathscr{B}(\mathbb{R})$ contains all intervals of the form above. 
\end{remark}

\begin{definition}
Let $\mathscr{A}_1$,$\mathscr{A}_2$ be $\sigma$-algebras on sets $\Omega_1,\Omega_2$, respectively. A function $f:\Omega_1\to \Omega_2$ is measurable relative to the $\sigma$-algebras $\mathscr{A}_i,i=1,2$ if for every $A\in \mathscr{A}_2,f^{-1}(A)\in \mathscr{A}_1$. In particular, the codomain $\Omega_2$ of the measurable functions we deal with later in this dissertation is $\mathbb{R}$ equipped with $\sigma$-algebra $\mathscr{B}(\mathbb{R})$. And we simply call them measurable functions (given that $\mathscr{A}_1$ is clear from the context).
\end{definition}

\begin{definition}\label{PS}
    A probability space is a triple $(\Omega,\mathscr{A},P)$, where $\mathscr{A}$ is a $\sigma$-algebra on a set $\Omega$ and $P$ is a non-negative measure on $\mathscr{A}$ such that $P(\Omega)=1$.
\end{definition}

\begin{definition}
    A real-valued random variable $X$ is a measurable real-valued function on a probability space $(\Omega,\mathscr{A},P)$, $X:\Omega\to \mathbb{R}$.
\end{definition}

\begin{remark}
    Given a measurable function $f$ on a probability space $(\Omega,\mathscr{A},P)$, and $I$ an interval of a kind described in Def. \ref{def:Borel inter.} or Remark \ref{remark:inter. type}, from now on if we write $P(f\in I)$, we mean $P\{x\in\Omega|f(x)\in I \}.$
\end{remark}

\begin{definition}
    If $X$ is a real-valued random variable on a probability space $(\Omega,\mathscr{A},P)$, the expectation of $X$ is 
    $$\mathbb{E}(X)=\int_{\Omega}XdP.$$
In particular, if $X$ only takes on finitely many values, $x_1,...,x_n\in \mathbb{R}$, then $$\mathbb{E}(X)=\sum_{i=1}^nx_iP(X=x_i).$$
\end{definition}

\begin{definition}
    If $X$ is a real-valued random variable on a probability space, the variance of $X$, denoted by $\operatorname{Var}(X)$, is 
  $$\operatorname{Var}(X):=\mathbb{E}\{(X-\mathbb{E}(X))^2\}=\mathbb{E}(X^2)-(\mathbb{E}(X))^2.$$
\end{definition}

The following result is well known in measure theory. For instance, see \cite[Prop.2.3.10]{Cohn} for reference.

\begin{proposition}
    Let $(\Omega,\mathscr{A},P)$ be a measure space, and let $X$ be a $[0,+\infty]$-valued $\mathscr{A}$-measurable function on $\Omega$. If $t$ is a positive real number and $A_t:=\{y\in \Omega:X(y)\geq t\}$, then $$P(A_t)\leq \frac{1}{t}\int_{A_t}XdP\leq \frac{1}{t}\int_{\Omega}XdP.$$
\end{proposition}

In particular, if $X$ is a real-valued random variable on a probability space $(\Omega,\mathscr{A},P)$, we get \textit{Markov's inequality}:
$$P(X\geq t)\leq \frac{\mathbb{E}(X)}{t}.$$

As a special case, if $X$ is integer-valued, then letting $t=1$ gives:
\begin{equation} \label{eq:1.1}
    P(X>0)\leq \mathbb{E}(X).
\end{equation}
Applying Markov's inequality, if $d>0$, we have \textit{Chebyshev's inequality}:
$$\mathbb{E}\{(X-\mathbb{E}(X))^2\}\geq P(|X-\mathbb{E}(X)|^2\geq d^2)d^2=P(|X-\mathbb{E}(X)|\geq d)d^2$$
In particular, if $\mathbb{E}(X)>0$, then
\begin{equation}\label{eq:1.2}
    P(X=0)\leq P(|X-\mathbb{E}(X)|\geq \mathbb{E}(X))\leq \frac{\operatorname{Var}(X)}{(\mathbb{E}(X))^2}.
\end{equation}

\section{Random graphs}
We begin with setting some notations used in this dissertation thereafter. Given an undirected, not necessarily simple graph $G$, let $G_V$ and $G_E$ denote the set of vertices and the set of edges of $G$, respectively. A graph $H$ is a \textit{subgraph} of a graph $G$, denoted by $H\subseteq G$, if $H_V\subseteq G_V,H_E\subseteq G_E$, and for every $e\in H_E$, the two endpoints of $e$ must also belong to $H_V$. Given $V'\subseteq G_V$ of a graph $G$, the \textit{subgraph induced} by $V'$ in $G$ is a subgraph $K\subseteq G$ such that $K_V=V'$ and $K_E$ consists of all edges in $G_E$ having both vertices in $V'$ as endpoints. 

The theory of random graphs was founded by Paul Erdős and Alfréd Rényi \cite{Erdos1,Erdos2,Erdos3,Erdos4} after Erdős \cite{Erdos}. There are two basic models of random graphs, the uniform model and the binomial model. The uniform model is defined as follows: Given a vertex set $V$, with $|V|=n,n\in \mathbb{N}$ and an integer $M$ such that $0\leq M\leq \binom{n}{2}$, take $\Omega$ to be the family of all graphs on $V$ with exactly $M$ edges, thus 
$$P(G)=\binom{\binom{n}{2}}{M}^{-1}, G\in \Omega.$$
In this dissertation, we restrict our attention to the binomial model. Readers who want to explore various topics in classical random graph theory can consult references such as Béla Bollobás \cite{Bollobas}, Svante Janson et al. \cite{Janson}, Alan Frieze and Michał Karoński \cite{Frieze}. In particular, \cite[Section 1.4]{Janson} gives a summary of conditions under which the two models are equivalent asymptotically.
Given a set $V=[n]=\{1,2,...,n\}$ and a real number $p$, $0\leq p\leq 1$, the binomial model $\mathscr{G}\{n,P(edge)=p\}$ consists of all 
simple graphs with vertex set $V$ with $|V|=n, n\in \mathbb{N}$ in which the edges are chosen independently with probability $p$. In the language of Definition \ref{PS}, we are working on a probability space $(\Omega,\mathscr{A},P)$ where $\Omega$ is the set of all simple graphs on $V$, denoted by $\mathscr{G}^n$. $\mathscr{A}$ is $\mathscr{P}(\mathscr{G}^n)$, and if $G$ is a simple graph on $V$ with $M$ edges, 
\renewcommand{\arraystretch}{0.6}
$0\leq M \leq \binom{n}{2}$, then
$$P(G)=p^M(1-p)^{\binom{n}{2}-M}$$
and if $A\in \mathscr{A}=\mathscr{P}(\mathscr{G}^n)$, then $A$ is a subset of $\Omega=\mathscr{G}^n$, say $A=\{G_1,...,G_r\}$, then $P(A)=P(G_1)+...+P(G_r)$. By Binomial Theorem, $P(\mathscr{G}^n)=1$.

This random-graph model is denoted by $G(n;p)$. Therefore, if $p=\frac{1}{2}$, then each graph is equally likely to appear. Very often, $p=p(n)$ is a function of $n$ and $p(n)\rightarrow 0$ as $n\rightarrow \infty$.

A subset $Q$ of $\mathscr{G}^n$ is called a \textit{property of graphs of order $n$} if it is closed under isomorphism. We may simply call it a \textit{property} should there be no ambiguity. When we say a graph $G$ has the property $Q$ we mean $G\in Q$. A property $Q$ is called \textit{monotone increasing} if whenever $G\in Q$ and $G\subseteq H$ then $H\in Q$. A property is \textit{monotone decreasing} if its complement is monotone increasing. For two functions $f(n),g(n)$ of $n$, we write $f\ll g$ and $f\gg g$ if $f/g\rightarrow 0$ and $f/g\rightarrow \infty$ as $n\rightarrow \infty$, respectively.

In $G(n;p)$, for a monotone increasing property $Q$, $\hat{p}=\hat{p}(n)$ is called a \textit{threshold} if 

\begin{equation*}
    P(Q)\rightarrow \left\{\begin{array}{ccccc}
         0 &\mbox{if}& p\ll \hat{p}  \\
         1 &\mbox{if}& p\gg \hat{p}
    \end{array}. \right. 
\end{equation*}

Lastly, we introduce the following standard notation for the asymptotic behavior of the comparison of the order of magnitude of two sequences $a_n,b_n$ as $n\rightarrow \infty$. Assume $b_n>0$ for sufficiently large $n$.
\begin{itemize}
    \item $a_n=o(b_n)$ if $a_n/b_n\rightarrow 0$ as $n\rightarrow \infty$.
    \item $a_n\ll b_n$ or $b_n\gg a_n$ if $a_n\geq 0$ and $a_n=o(b_n)$.
    \item $a_n=O(b_n)$ if there exist $C,n_0>0$ such that $|a_n|\leq Cb_n$ for $n\geq n_0$ as $n\rightarrow \infty$.
    \item $a_n=\Omega(b_n)$ if there exist $c,n_0>0$ such that $a_n\geq cb_n$ for $n\geq n_0$ as $n\rightarrow \infty$.
    \item $a_n=\Theta(b_n)$ if $a_n=O(b_n)$ and $a_n=\Omega(b_n)$ as $n\rightarrow \infty$. In other words, if there exist $C,c,n_0>0$ such that $cb_n\leq a_n\leq Cb_n$ for $n\geq n_0$.
\end{itemize}
    
\chapter{Randomly Generated Matroid}
In this Chapter, we introduce a way of generating matroids by randomly picking indicator vectors. In many circumstances, we get uniform matroids. Since we see that uniform matroids are transversal, we are going to write down the number of minimal presentations of uniform matroids, so that we will have a sense of how many transversal matroids are in the collection of matroids we generate by this model.
\section{A model of randomly generated matroids using indicator vectors}
We try to mimic the binomial model of random graphs mentioned in section 1.4 to create a model of randomly generated matroids on a finite set. Given a natural number $n$, we let $E_n$ be the set $\{x_1,...,x_n\}$. An \textit{indicator vector} for $E_n$ is a vector $v=(v_1,...,v_n)\in (\mathbb{F}_2)^n$. $v_1,...,v_n$ are called \textit{coordinates} of $v$. When we choose a set $V$ of indicator vectors for $E_n$, we always choose in such a way that all vectors in $V$ have the same number of 1's in their coordinates. This is because, in a matroid, all bases have the same cardinality. If $v$ is an indicator vector for $E_n$, let $v(1)=\{x_i\mid v_i=1,1\leq i\leq n\}$. For a set $V$ of indicator vectors for $E_n$, we let $H(V)=\{v(1)\mid v\in V\}$. Moreover, if $V=\emptyset$, let $H(V)=\{\emptyset\}$.

Given a real number $p$, $0\leq p\leq 1$ and $m\leq n$, let each indicator vector for $E_n$ with $m$ ones in their coordinates be chosen independently with probability $p$. Thus, we are working in a probability space $(\Omega,\mathscr{A},P)$, where $\Omega$ is the power set of all indicator vectors with $m$ ones of $E_n$, $\mathscr{A}$ is $\mathscr{P}(\Omega)$. For each collection of vectors $\{v_1,...,v_k\}$,  
$$P(\{v_1,...,v_k\})=p^k(1-p)^{\binom{n}{m}-k}.$$

Clearly, $H(V)$ may not form a collection of independent sets of a matroid. We make the following design.

\begin{definition}\label{def:exchange}
    Let $v,w$ be two indicator vectors for $E_n$ (with the same number of 1's in their coordinates), $i,j\in [n]$ are a compatible pair of indices for $v,w$ if $v_i=1, w_i=0$ and $v_j=0,w_j=1$ and let $v_w^{(i,j)}$ denote the indicator vector for $E_n$ with
    $$ \left\{ \begin{array}{ccc}
         (v_w^{(i,j)})_k=0  & \mbox{for} & k=i
        \\(v_w^{(i,j)})_k=1 &\mbox{for}&   k=j
        \\ (v_w^{(i,j)})_k=v_k &\mbox{for}& 1\leq k\leq n, k\neq i,k\neq j.
                      \end{array} \right. $$
Moreover, let $$E_{v,w}=\{v_w^{(i,j)}\mid i,j\in [n] \, \mbox{is a compatible pair of indices for}\, v,w \}.$$
If $V$ is a set of indicator vectors for $E_n$, we let $V_U=V\cup\displaystyle\bigcup_{v,w\in V} (E_{v,w}\cup E_{w,v}).$
\end{definition}

Intuitively, for every two distinct indicator vectors $v,w$ (with the same number of 1's in their coordinates), there must be $i,j\in [n]$ such that $v_i=1,w_i=0$ and $v_j=0,w_j=1$. $v_w^{(i,j)}$ is the indicator vector for $E_n$ obtained by replacing the 1 at the $i$-th position of $v$ with 0 in exchange for the 1 at the $j$-th position from $w$. $E_{v,w}$ collects all vectors obtained in this way. And we also want those vectors obtained in the same way with the roles of $w$ and $v$ swapped, that is the set $E_{w,v}$. Taking union of all those vectors together with $V$, we get $V_U$.

For every $h,k\in V_U$, if we also want $h,k$ together with $E_{h,k}$ and $E_{k,h}$, we can repeat the process in the preceding definition to get $(V_U)_U$. Continuing this process we get a chain:
\begin{equation}\label{eq:chain V_U}
    V\subseteq V_U\subseteq (V_U)_U\subseteq ((V_U)_U)_U\subseteq ...\subseteq
\end{equation}
Since there are finitely many indicator vectors for $E_n$ for a given $n\in \mathbb{N}$, \eqref{eq:chain V_U} must stabilize. 

\begin{definition}\label{M(V)}
    Given a collection $V$ of indicator vectors (with same number of 1's in their coordinates) for $E_n$, where $n$ is a given natural number, let $V_{stab}$ denote the set $((V_U)...)_U$ in chain \eqref{eq:chain V_U} where \eqref{eq:chain V_U} stabilizes, and $M(V)$ be the matroid on $E_n$ with $\mathscr{B}(M(V))=H(V_{stab}).$
\end{definition}    
    It is quite straightforward to see that $\mathscr{B}(M(V))$ is a collection of bases of a matroid. To be precise, if $A,B\in H(V_{stab})$, with $x_i\in A\setminus B$, there are indicator vectors $u=\{u_1,...,u_n\}, v=\{v_1,...,v_n\}$ for $E_n$ such that $H(u)=A$ and $H(v)=B$. Since $x_i\in A\setminus B$, $u_i=1$ while $v_i=0$, there must be an index $j$ so that $v_j=1$ and $u_j=0$ with $u_v^{(i,j)}\in V_{stab}$. Therefore, $H(u_v^{(i,j)})$ is $(A\setminus \{x_i\})\cup \{x_j\}$. So we see that $H(V_{stab})$ satisfies the conditions in Theorem \ref{BE}. This is a very canonical, and the most obvious way of generating a matroid by a collection of chosen indicator vectors because we force every vector resulting from an arbitrary exchange from vectors in $V_{stab}$ to be included.  

\begin{example}
    Let $V$ be the collection of indicator vectors $\{(1,1,0,0),(0,0,1,1)\}$ for $(x_1,x_2,x_3,x_4)$. Then \\$V_{stab}=\{(1,1,0,0),(0,0,1,1),(0,1,1,0),(0,1,0,1),(1,0,1,0),(1,0,0,1)\}$. \\$M(V)=(\{x_1,x_2,x_3,x_4\},\{\{x_1,x_2\},\{x_3,x_4\},\{x_2,x_3\},\{x_2,x_4\},\{x_1,x_3\},\{x_1,x_4\}\})$.
\end{example}

\section{Monotone properties of $M(V)$}
If $m\leq n$, a \textit{property $Q$ of random matroids of order $m,n$} is a class of $M(V)$, which is closed under isomorphism, where $V$'s are sets of indicator vectors for $E_n$ with $m$ ones in their coordinates. We may simply call it a \textit{ property} should there be no ambiguity. When we say $M(V)$ has a property $Q$, we mean $M(V)\in Q$. A property $Q$ is \textit{monotone increasing} if whenever $V_1\subseteq V_2$ and $M(V_1)\in Q$, then $M(V_2)\in Q$. 

\begin{lemma}\label{BF}
    If $v_1,v_2$ are indicator vectors with $m$ ones for $E_n$, $m\leq n$, then $\mathscr{B}(M(v_1,v_2))$ consists of all subsets of $E_n=\{x_1,...,x_n\}$ of the form:
    \begin{equation}\label{eq:(v_1(1) u v_2(1)) u X}
    (v_1(1)\cap v_2(1))\cup X
    \end{equation}
    where $$X\in \left \{Y\in \mathscr{P}((v_1(1)\cup v_2(1))\setminus(v_1(1)\cap v_2(1)))\Bigm   ||Y| =m-|v_1(1)\cap v_2(1)|\right \}.$$
    \end{lemma}
\begin{proof}
    If $m=n$, then there is nothing to prove. Suppose $m<n$, and there is a subset $J$ of $E_n$ in the form of \eqref{eq:(v_1(1) u v_2(1)) u X} such that $J\notin \mathscr{B}(M(v_1,v_2))$. Let $h(u,v)$ be the Hamming distance between two 0-1 vectors having $n$ coordinates. And let $y$ be the 0-1 vector such that $y(1)=Y.$ There must be coordinates $i,j$ such that $v_{1i}=1$, $v_{2i}=y_i=0$ and $v_{1j}=0$, $v_{2j}=y_j=1$. Such a pair of coordinates must exist, otherwise $y$ would be $v_1$ or $v_2$. Perform the exchange described in Definition \ref{def:exchange} to get $v_{1\hspace{.5mm} v_2}^{(i,j)}.$ Then $h(y, v_{1\hspace{.5mm} v_2}^{(i,j)})=h(v_1,y)-2$. Repeat this exchange process until we get a vector having Hamming distance 0 with $y$. That means we get $y$ from the exchange process. That is a contradiction.
\end{proof}

If there is a $v\in V$ such that $x_i\in v(1)$, $x_i$ is not a loop of $M(V)$. Similarly, if $x_i \notin v(1)$, $x_i$ is not a coloop of $M(V)$. Therefore, $M(V)$ containing a non-loop/non-coloop, and not containing a loop/coloop are monotone increasing properties.   

\begin{theorem}
    Let $m\leq n$. Suppose the probability of choosing each 0-1 indicator vector for $E_n$ with $m$ ones is $0\leq p\leq 1$. Let $V$ denote the set of vectors being chosen, $p^*(n)=\binom{n-1}{m-1}^{-1}$ is the threshold for $x_i$ to be a non-loop in $M(V),1\leq i\leq n$.
\end{theorem}
\begin{proof}
    Let $\mathscr{U}$ be the set of all 0-1 indicator vectors for $E_n$, there are $\binom{n-1}{m-1}$ indicator vectors in $\mathscr{U}$ having 1 at coordinate $x_i$. Denote them by $v_1,v_2,...,v_{\binom{n-1}{m-1}}$.
    Let $I_j:\mathscr{P(U)}\rightarrow \{0,1\},1\leq j\leq \binom{n-1}{m-1}$ be given by
   $$I_j(V)= \left\{ \begin{array}{ccc}
         1  & \mbox{if} & v_j\in V
        \\0 &\mbox{if}&   v_j\notin V
                      \end{array} \right.$$
Set $$X=I_1+...+I_{\binom{n-1}{m-1}}$$ then $X$ is a random variable counting the number of elements of $\{v_1,...,v_{\binom{n-1}{m-1}}\}$ in $V$. Note that the random variable $I_j$, where $1\leq j\leq \binom{n-1}{m-1}$ are pairwise independent random variables since indicator vectors are chosen independently in our model. The probability  $P(I_i=1)$ is by definition, $p$. Therefore,

$$\mathbb{E}(X)=\sum_{j=1}^{\binom{n-1}{m-1}}1\cdot P(I_j=1)=\binom{n-1}{m-1}p.$$
When $p\ll\binom{n-1}{m-1}^{-1}$, $\mathbb{E}(X)\rightarrow 0$ as $n\rightarrow \infty$. And by \eqref{eq:1.1}, $P(X>0)\rightarrow 0$. Now, note that 
$$\operatorname{Var}(X)=\mathbb{E}(\sum_i I_i^2+\sum_{i\neq j}I_iI_j)-(\sum_i\mathbb{E}I_i)^2$$
$$=\sum_i\mathbb{E}(I_i^2)+\sum_{i\neq j}\mathbb{E}(I_iI_j)-\sum (\mathbb{E}(I_i))^2-\sum_{i\neq j}(\mathbb{E}(I_i))(\mathbb{E}(I_j))$$
Since every indicator vector is chosen independently with probability $p$, $I_j$'s are pairwise independent variables, so 
$$\sum_{i\neq j}\mathbb{E}(I_iI_j)-\sum_{i\neq j}(\mathbb{E}(I_i))(\mathbb{E}(I_j))=0.$$

As a result, we get
\begin{dmath}\label{eq:Var of Berno}
 \operatorname{Var}(X)=\sum_i 1\cdot P({I_i^2=1})-1\cdot P({I_i=1})^2=\sum_{i}1\cdot P({I_i=1})-1\cdot P({I_i=1})^2=\sum_i p-p^2=\binom{n-1}{m-1}p(1-p)=\mathbb{E}(X)(1-p).
 \end{dmath}
Then by \eqref{eq:1.2}, $P(X=0)\leq \frac{\operatorname{Var}(X)}{(\mathbb{E}(X))^2}=\frac{1-p}{\mathbb{E}(X)}$, but when instead $p\gg\binom{n-1}{m-1}^{-1}$, $\mathbb{E}(X)\rightarrow \infty$. Therefore, $P(X=0)\rightarrow 0$.
\end{proof}

\begin{proposition}
    For any natural numbers $h,q$ with $h\leq q$, the property ``$M(V)$ contains $U_{h,q}$ as a minor'' is monotone increasing.
\end{proposition}
\begin{proof}
    Let $m\leq n$ and $V=\{v_1,\ldots,v_k\}$ be a collection of 0-1 vectors with $m$'s 1 for $E_n$, $0\leq k\leq \binom{n}{m}$. If $m$ or $k$ equals 0, then there is nothing to prove. Thus, we assume both of them are bigger than or equal to 1. If $m=n$, then the whole ground set of $M(V)$ is a basis, adding more indicator vectors does not change anything. So we assume that $m<n$.  Suppose $M(V)$ contains $U_{h,q}$ as a minor, $h\leq m$, $q\leq n$, then there exists a $q$-subset of $\{x_1,\ldots,x_n\}$ such that each of its $h$-subsets is contained in a $v(1)$ for some $v\in V_{stab}$. Without loss of generality, we assume that $q$-subset to be $\{x_1,\ldots,x_q\}$. Note that if $h=q$, no matter what indicator vector we add to $V$, $M(V\cup w)|\{x_1,\ldots,x_q\}$ gives us $U_{h,q}$. Thus, we further assume that $h<q$. Now, let's add a new set of indicator vectors $V'$ with $m$ ones to $V$, if every vector in $U_{V\cup V'}$ has less than or equal to $h$ ones in 1st to $q$th coordinates , $M(V\cup V')|\{x_1,\ldots,x_q\}$ gives us $U_{h,q}$. Suppose that there is some vector in $V'$ having strictly more than $h$ ones in  coordinates $\{x_1,\ldots,x_q\}$. Pick a set $W$ of $\binom{q}{h}$ vectors from $V_{stab}$ such that for each $h$-subset of $\{x_1,\ldots,x_q\}$, there is a vector from $W$ such that it has ones corresponding to that $h$-subset. Since all vectors in $(V\cup V')_{stab}$ have $m$ ones, each vector in $W$ has $m-h$ ones in $q+1$th to $n$th coordinates. By using the exchange process in Definition \ref{M(V)} and Lemma \ref{BF}, we can get a set of indicator vectors $W_E$ such that each of them agrees with exactly one vector from $W$ in 1st to $q$th coordinates, and all of them have ones in $q+1$th to $q+m-h$ th coordinates, then $(M(V\cup V')\setminus\{x_{q+m-h+1},\ldots,x_n\})/\{x_{q+1},...,x_{q+m-h}\}\cong U_{h,q}$. 
    \end{proof}

\section{Uniform matroid as transversal matroid}
\begin{lemma}\label{typical mini presentation}
    The collection of subsets $\mathscr{A}=\{A_1,...,A_m\}$ of $[n]$, where $A_i=\{i,m+1,...,n-1,n\}$, $1\leq i\leq m$, is a minimal presentation for the uniform matroid of rank $m$ on $[n]$.
\end{lemma}
\begin{proof}
    First, we show that $\mathscr{A}$ is a presentation for $U_{m,n}$. $\{1,...,m\}$ is indeed a transversal of $\mathscr{A}$. If $\{x_1,...,x_m\}\subseteq [n]\setminus [m]$, it is still a transversal of $\mathscr{A}$ since each $A_i$ contains all of $[n]\setminus [m]$. Now, given $\{y_1,...,y_k,y_{k+1},...,y_m\}$ such that $\{y_1,...,y_k\}\subseteq [m]$, each $y_j$ belongs to a unique $A_{i_j}$, $1\leq j\leq k$. And since every $y_l$, $k+1\leq l\leq m$, belongs to every member of $\mathscr{A}\setminus \{A_{i_1},...,A_{i_k}\}$ by design, $\{y_1,...,y_m\}$ is a transversal of $\mathscr{A}$. Now we want to show that the presentation $\mathscr{A}$ is minimal. Clearly, removing any $i$ from $A_i$ makes $\{1,...,m\}$ a non-transversal. If instead, we only remove $m+r$, $1\leq r\leq n-m$ from $A_1$, without loss of generality, then  $\{m+r,2,...,m\}$ becomes a non-transversal.
\end{proof}

If there is another minimal presentation for $U_{m,n}$, is it possible for it to have strictly fewer subsets of $[n]$ than $\mathscr{A}$? The next lemma shows that it is impossible.

\begin{lemma}\label{mini presentation n=m}
    If $\mathscr{A}=\{A_1,...,A_n\}$ is a minimal presentation for $M(\mathscr{A})$, then the rank of $M(\mathscr{A})$ is $n$.
\end{lemma}
\begin{proof}
    Suppose the rank of $M(\mathscr{A})$ is $m$, $m<n$. Note that by hypothesis each $A_i$, $1\leq i\leq n$, is non-empty. Then there must be a partial transversal of $\mathscr{A}$ of maximal cardinality $m$. Without loss of generality, we let that be $\{e_1,...,e_m\}$, where $e_i\in A_i$, $1\leq i\leq m$. Observe that:
    \begin{enumerate}
    \item Any set in $\{A_{m+1},...,A_n\}$ can contain only elements from $\{e_1,...,e_m\}$. 
    \item If a set $A_j$ contains $e_i$, $m+1\leq j\leq n$, $1\leq i\leq m$, then $A_i$ can only contain elements from $\{e_1,...,e_m\}$.
    \end{enumerate}
    Otherwise, we get a partial transversal of $\mathscr{A}$ of cardinality strictly bigger than $m$, contradicting our hypothesis. Then, we claim that erasing $A_{m+1},...,A_n$ from $\mathscr{A}$ still gives us the same collection of partial transversals. The proof is then complete since we derive a contradiction to our assumption that $\mathscr{A}$ is a minimal presentation for $M(\mathscr{A})$. To justify our claim, we need to show that any partial transversal of cardinality $m$ of $\{A_1,...,A_m,...,A_n\}$ must be a partial transversal of $\{A_1,...,A_m\}$. Certainly, if the partial transversal is $\{e_1,...,e_m\}$, we are done. If it is $\{e_1',...,e_m'\}$ and is disjoint from $\{e_1,...,e_m\}$, by (1) above, we are also done. So we suppose that $\{A_1,...,A_n\}$ admits a partial transversal whose intersection with $\{e_1,...,e_m\}$ has cardinality $k$, where $1\leq k<m$. Let us write it as $\{e_{r_1},...,e_{r_k},e_{k+1}',...,e_{m}'\}$, where $\{r_1,...,r_k\}\subseteq \{1,...,m\}$. If some element in $\{e_{r_1},...,e_{r_k}\}$ belongs to some $A_j$ for $m+1\leq j\leq n$, say $e_{r_1}$. Note that, by (1), $\{e_{k+1}',...,e_m'\}$ is a partial transversal of $\{A_{i_{k+1}},...,A_{i_{m}}\}$, where $\{i_{k+1},...,i_m\}\subseteq \{1,...,m\}$. But our hypothesis says that $e_{r_1}$ belongs to some $A_i$, $1\leq i\leq m$, and by 2, $e_{r_1}$ does not belong to any set in the collection $\{A_{i_{k+1}},...,A_{i_{m}}\}$. This argument implies that every element in $\{e_{r_1},...,e_{r_k}\}$ belongs to a set in $\{A_1,...,A_m\}\setminus \{A_{i_{k+1}},...,A_{i_{m}}\}$. But from our hypothesis, $\{e_{r_1},...,e_{r_k}\}$ is a partial transversal of $\{A_1,...,A_m\}$, this forces it to be a partial transversal of $\{A_1,...,A_m\}\setminus \{A_{i_{k+1}},...,A_{i_{m}}\}$.     
\end{proof}

Does every member of a minimal presentation of $U_{m,n}$ share the same cardinality? We need more tools to tackle this. Below is an obvious lemma:

\begin{lemma}
    Given a presentation $\mathscr{A}=\{A_1,...,A_n\}$ of $M(\mathscr{A})$ on a finite set $E$, $A_i\subseteq E$ for all $1\leq i\leq n$, $E\setminus A_i$ is a flat of $M(\mathscr{A})$.
\end{lemma}
\begin{proof}
     Obviously, if $A_i$ is empty, the statement is clearly true. Suppose $A_i$ is non-empty and the rank of $E\setminus A_i$ is $m$ in $M(\mathscr{A})$, then there is a set $\{e_1,...,e_m\}\subseteq E\setminus A_i$ such that $e_j\in A_{i_j}$ for $1\leq j\leq m$, and $\{A_{i_1},...,A_{i_m}\}\subseteq \{A_1,...,A_{i-1},A_{i+1},...,A_n\}$. Then $\{e_1,...,e_m,e'\}$ is still a partial transversal for $\mathscr{A}$ for any $e'\in A_i$.
\end{proof}

\begin{corollary}
    If $\mathscr{A}=\{A_1,...,A_m\}$ is a minimal presentation for the uniform matroid $U_{m,n}$ on a finite set $E$, then $A_i$ must have cardinality at least $n-m+1$ for all $1\leq i\leq m$.
\end{corollary}
\begin{proof}
    Suppose there is an $A_i$ with cardinality $k<n-m+1$, then by our previous Lemma, $E\setminus A_i$ is a flat for $U_{m,n}$ with cardinality $n-k>m-1$. That means this flat of $U_{m,n}$ is of cardinality at least $m$, which is impossible, since all $m$-element subsets of $E$ are bases in $U_{m,n}$ and adding elements to $E\setminus A_i$ does not increase the rank of it, which is a contradiction.
\end{proof}

But is it possible for the cardinality of an $A_i$ to be strictly bigger than $n-m+1$? We need the following lemma from Bondy and Welsh:

\begin{lemma}{\cite[Bondy and Welsh]{Bondy}}
    Let $\mathscr{A}=\{A_1,...,A_m\}$ be a presentation of a matroid $M$ with $\rho(M)=n$. Fix $j\in [m]$, if $P$ is a transversal of $\{A_1,...,A_{j-1},A_{j+1},...,A_m\}$ for which $|P\cap A_j|$ is minimal, then
    \begin{enumerate}
        \item $A_j-P$ is a cocircuit of $M$, and
        \item replacing $A_j$ by $A_j-P$ in $\mathscr{A}$ yields a presentation of $M$.
    \end{enumerate}
\end{lemma}

\begin{corollary}\label{n-m+1}
If $\mathscr{A}=\{A_1,...,A_m\}$ is a minimal presentation for $U_{m,n}$, then $|A_i|=n-m+1$ for all $1\leq i\leq m$.
\end{corollary}
\begin{proof}
    Suppose $|A_j|>n-m+1$, pick any transversal $P={e_1,...,e_{j-1},e_{j+1},...,e_m}$ for $\{A_1,...,A_{j-1},A_{j+1},...,A_m\}$. By the pigeonhole principle, $|P\cap A_j|\geq 1$. And by 2. of previous lemma, replacing $A_j$ by $A_j-P$ in $\mathscr{A}$ yields a presentation for $U_{m,n}$, contradicting the assumption that $\mathscr{A}$ is a minimal presentation for $U_{m,n}$. 
\end{proof}

It is worth pointing out that, in \cite{Bonin} Bonin and de Mier proved the following result, which can imply our preceding corollary. The rank of $(U_{m,n})^*$ is $n-m$; moreover, every subset of $E((U_{m,n})^*)$ of cardinality $n-m+1$ is a circuit of $U_{m,n}^*$.
\begin{proposition}[{\cite[Bonin and de Mier]{Bonin}}]
    A presentation $\mathscr{A}$ of a transversal matroid $M$ is minimal if and only if each member of $\mathscr{A}$ is a cocircuit of $M$.
\end{proposition}

The following theorem is essential for us to count the number of distinct minimal presentations of a uniform matroid.

\begin{theorem}[{\cite[Mirsky and Perfect]{Mirsky}}]
    Let $\mathscr{A}=\{A_i:i\in I\}$ and $\mathscr{B}=\{B_i:i\in I\}$ be two finite families of subsets of $E$. Then $\mathscr{A}$ and $\mathscr{B}$ have a common transversal which contains a prescribed set $C$ if and only if for all subsets $J,K$ of $I$,
    \begin{equation}\label{eq:Transversal}
                |(A(J)\cap B(K))\setminus C|+|A(J)\cap C|+|B(K)\cap C|\;\geq |J|+|K|+|C|-|I| 
    \end{equation}
    where $A(J)=\bigcup_{i\in J}A_i$, $B(K)=\bigcup_{i\in K}B_i$.
\end{theorem}

If we want another family of sets to induce the same matroid as the one in Lemma \ref{typical mini presentation}, then they have to satisfy \eqref{eq:Transversal} for all $m$-element subsets $C$. As a result, we have:

\begin{corollary}
    Let $\mathscr{A}$ be the family of sets in Lemma \ref{typical mini presentation}, then a family of sets $\mathscr{B}=\{B_i:i\in I\}$ induces the same matroid as $\mathscr{A}$, namely, $U_{m,n}$, if and only if for all $m$-element subsets $C$, all $K\subseteq I$
    \begin{equation}\label{eq:uniform}
        |B(K)\cap C|\;\geq |K|.
    \end{equation}
\end{corollary}
\begin{proof}
    Note that $|C|-|I|\;=0$, the only if part follows from taking $J=\emptyset$. For the converse, note that if $J\subseteq I$, then $|A(J)|\;\ge n-m+|J|$. By the pigeonhole principle, we have $|A(J)\cap C|\;\geq |J|$. Since $|A(J)\cap C|\;\geq |J|$ for all $J$, \eqref{eq:uniform} is clearly sufficient for \eqref{eq:Transversal} to hold.
\end{proof}

If there is $K\subseteq I$ such that $|B(K)\cap C|\;<|K|\;=t$ for some $m$-subset $C$, then we say that there are $t$ members of $\mathscr{B}$ violating \eqref{eq:uniform}. It is important to note that in such a case, it is possible to have $K'\subset K$ such that $|B(K')\cap C|\;<|K'|$ for some $m$-subset $C$. 

\subsection{Number of minimal presentations for uniform matroid} \label{subsec:no. of mini. pre.}
Earlier in this section, we saw that every uniform matroid is transversal. Before we move to the study of transversal matroids, we want to have an idea of how many uniform matroids there are in the class of transversal matroids. Because we want to take the fact ``a transversal matroid can have many different presentations'' into consideration, it is natural to ask how many minimal presentations there are for a uniform matroid. 

In view of Lemma \ref{mini presentation n=m} and Corollary \ref{n-m+1}, apart from satisfying \eqref{eq:uniform}, there must be exactly $m$ distinct members in $\mathscr{B}=\{B_i:i\in I\}$ and each member must be of cardinality $n-m+1$ for $\mathscr{B}$ to be a minimal presentation. Note that $|B_i\cap C|\;\geq 1$ for all $i$ and every $m$-element subset $C$ due to the pigeonhole principle. And given any $K\subseteq I$ with $|K|=2$, $B(K)$ satisfies \eqref{eq:uniform}, because $B(K)$ must contain at least one element which is not in $B_i$ or $B_j$, $K=\{i,j\}$. Thus, the minimum cardinality of $K$ such that $\mathscr{B}$ violates \eqref{eq:uniform} is 3. That is, a violation occurs if there exists $i,j,k\in I$ such that $|B_i\cup B_j\cup B_k|=n-m+2$. Let us begin with a simple example to illustrate our strategy. 

\begin{example}
    Find the number of distinct minimal presentations of $U_{3,4}$. What we try to find is the number of ways of choosing $\mathscr{B}=\{B_1,B_2,B_3\}$ such that $B_i$'s are pairwise distinct subsets of $[4]$ with cardinality 2, and $\mathscr{B}$ satisfies \eqref{eq:uniform}. First, ignoring \eqref{eq:uniform}, the number of ways of choosing $\mathscr{B}$ is $\displaystyle\binom{\binom{4}{2}}{3}$. Now we want to subtract the cases violating \eqref{eq:uniform} from it. By what we observe in the above paragraph, the minimal number of members that can violate \eqref{eq:uniform} is 3. By the pigeonhole principle, the cases such that $|B_1\cup B_2 \cup B_3|=3$ are the cases we need to abandon. There are $\binom{4}{3}\binom{\binom{3}{2}}{3}=4$ of them. Explicitly, they are
    \begin{align*}
    &\{\{2,3\},\{2,4\},\{3,4\}\},\{\{1,3\},\{1,4\},\{3,4\}\},\{\{1,2\},\{1,4\},\{2,4\}\},\\&\{\{1,2\},\{1,3\},\{2,3\}\}.
    \end{align*}
    Therefore, the number of distinct minimal presentations for $U_{3,4}$ is $\displaystyle\binom{\binom{4}{2}}{3}-4$.
\end{example}

\begin{example}
    We compute the number of minimal presentations for $U_{3,n}$ for any natural number $n\ge 3.$ We want to find the number of ways to choose a collection $\mathscr{B}=\{B_1,B_2,B_3\}$ of subsets of $[n]$ such that $\mathscr{B}$ satisfies equation \eqref{eq:uniform}, and each member of $\mathscr{B}$ has cardinality $n-2$. Due to the pigeonhole principle, any $\mathscr{B}$ we choose should satisfy:
    $$\left|\bigcup_{i\in I} B_i  \right|\ge n-3+|I|, I\subseteq [3]$$
    Equivalently, we can compute the number of ways to choose a collection $\mathscr{D}=\{D_1,D_2,D_3\}$ such that $D_i=[n]\setminus B_{\sigma(i)}$ for some $\{B_1,B_2,B_3\}$, $\sigma\in S_3.$ It is straightforward to see that $\mathscr{D}$ is a collection we want if and only if 
    $$\left|\bigcap_{i\in I}D_i\right|\leq 3-|I|, I\subseteq [3].$$
    In particular, $|D_i|\;=2, \mbox{ for all } i\in [3].$ If we think of each $D_i$ as an edge of the complete graph $K_n$ of order $n$, then our goal is to find the number of ways to choose three edges from $K_n$ such that they do not share a common vertex. This number is equal to
    $$\binom{\binom{n}{2}}{3}-n\binom{n-1}{3}.$$
\end{example}
  
We may write down a formula for the general case $U_{m,n}$, and the inspiration of this approach came from ChatGPT 5 pro. Our goal is to find the number of collections $\mathscr{B}=\{B_i:i\in [m]\}$ of subsets of $[n]$ of cardinality $n-m+1$ so that 
$$\left|\bigcup_{i\in I}B_i\right|\geq n-m+|I|, I\subseteq [m].$$
Similar to our preceding example, it is equivalent to finding the number of ways to choose $\mathscr{D}=\{D_i:i\in [m]\}$ such that 
$$\left|\bigcap_{i\in I}D_i\right|\leq m-|I|, I\subseteq [m].$$
In particular, $|D_i|\;=m-1 \mbox{ for all } i\in [m].$ For $S\subseteq [m]$, let $c_S$ be the number of elements in $[n]$ such that the members of $\mathscr{D}$ to which they belong are exactly $\{D_j:j\in S\}$. We want a collection $\mathscr{D}$ such that 
\begin{equation}\label{eq:Sum c_S=n}
\sum_{S\subseteq [m]}c_S=n
\end{equation}
This is because each element in $[n]$ must be counted once for exactly one subset $S\subseteq [m].$ Moreover, we should have 
\begin{equation}\label{eq:sum c_S m-1}
    \sum_{S\subseteq [m],i\in S} c_S=m-1, \mbox{ for every } i\in [m]
\end{equation}
and
\begin{equation}\label{eq:sum c_S m-|I|}
    \sum_{S\subseteq [m],I\subseteq S} c_S\leq m-|I|, \mbox{ for every } I\subseteq [m]
\end{equation}
Then the number of ways to choose $\mathscr{D}$ is 
\begin{equation*}
    \frac{1}{m!} \sum_{\{c_S:S\subseteq [m]\}} \frac{n!}{\prod_{S\subseteq [m]}c_S!}
\end{equation*}
where the sum is taken over all $\{c_S:S\subseteq [m]\}$ satisfying Equations \eqref{eq:Sum c_S=n}, \eqref{eq:sum c_S m-1} and \eqref{eq:sum c_S m-|I|}. 

\chapter{Monotonicity}
Recall that whenever we mention a bipartite graph $G$, we adopt the convention in Subsection \ref{G_N,G_M}. Remember that it gives us a collection of subsets of $G_N$. The transversal matroid on $G_N$ induced by this collection of subsets is denoted by $T(G)$. Unless otherwise specified, in this chapter, we assume that $|G_N|\;=|G_M|$. This assumption allows us to investigate the behavior of random bipartite graphs as $|G_N|\;=|G_M|\rightarrow \infty$. If $|G_N|\;\neq |G_M|$ then we need to let them go to infinity respectively, and that may cause further complications. The vertices in $G_M,G_N$ which are connected to some edges of $G$ are denoted by $G_M^c,G_N^c$, respectively. In other words, $G_M\setminus G_M^c$ and $G_N\setminus G_N^c$ are isolated vertices of $G$. 

 In this chapter, we will be investigating the threshold for a fixed subgraph to appear in random bipartite graphs, since we will see in Theorem \ref{thm:truncation}, such bipartite graphs correspond to transversal matroids containing a fixed matroid as a submatroid. We will also investigate the threshold for randomly generated transversal matroids to be simple. Any bipartite graph $G$ such that $T(G)$ is simple has to satisfy the conditions that $G_N$ has no isolated vertices and any subset of $G_N$ of cardinality two has a matching. 

A \textit{property} (of order $n$) $\mathscr{E}$ of transversal matroids is a class of transversal matroids on a finite set (of cardinality $n$) closed under isomorphism, for some nonnegative integer $n$. Whenever we discuss a property, we always mean a property of a given order $n$. $\mathscr{E}$ is \textit{monotone increasing} if whenever $M\in \mathscr{E}$ and $G\subseteq G'$ for some bipartite graph $G$ such that $T(G)=M$, then $T(G')\in \mathscr{E}$. $\mathscr{E}$ is monotone decreasing if $\mathscr{E}^c$ is monotone increasing. $\mathscr{E}$ is \textit{convex} if, whenever $M,M'\in \mathscr{E}$, and $G_1\subseteq G_2\subseteq G_3$ for some bipartite graphs $G_1,G_3$ such that $T(G_1)=M, T(G_3)=M'$, then $T(G_2)\in \mathscr{E}$. Clearly, every monotone property is convex but not vice versa.
Given a graph $G$, the class of graphs containing $G$ as a subgraph is clearly a monotone increasing property. It is a classical topic in random graph theory to investigate the thresholds for the graph to appear in $G(n;p)$ (The model of random graphs in which every edge is chosen independently with probability $p$). In Section \ref{sec:fixed subgraph} of this chapter, we are also going to investigate that in the bipartite-graph setting. And we show that, under certain conditions, $H\subseteq G$ implies that $T(G)$ contains $T(H)$ as a submatroid. In Section \ref{sec:simple T(G)}, we show that ``$T(G)$ is simple" is a monotone increasing property and we find its threshold.

\section{Bipartite graphs containing a fixed subgraph}\label{sec:fixed subgraph}
\begin{definition}
    A homomorphism from a bipartite graph $G_1$ to another bipartite graph $G_2$ is a pair of functions $f_N \colon G_{1N}\rightarrow G_{2N}$ and $f_M\colon G_{1M}\rightarrow G_{2M}$ such that whenever $(n,m)$ is an edge of $G_1$, where $n\in G_{1N}$ and $m\in G_{1M}$, we have $(f_N(n),f_M(m))$ is an edge of $G_2$. If both $f_N,f_M$ are bijective and the pair of their inverses is also a homomorphism from $G_2$ to $G_1$, then it is an isomorphism between $G_1$ and $G_2$. An automorphism of a bipartite graph is an isomorphism of the bipartite graph to itself.
\end{definition}

\begin{lemma}
    Given a bipartite graph $G$, there are $$\binom{n}{|G_N|}\binom{n}{|G_M|}|G_N|!|G_M|!\Big / |Aut(G)| $$ copies of $G$ in the complete bipartite graph $K_{n,n}$.
\end{lemma}
\begin{proof}
    There are $\binom{n}{|G_N|}\binom{n}{|G_M|}$ ways of choosing the vertices of $G$. And there are $|G_N|!|G_M|!$ pairs of maps, each pair maps $G_N$ to $(K_{|G_N|,|G_M|})_N$ and $G_M$ to $(K_{|G_N|,|G_M|})_M$. Note that each pair of maps is a homomorphism from $G$ to $K_{|G_N|,|G_M|}$. Each homomorphism induces a copy of $G$ and two pairs of maps $(f_N,f_M),(f'_N,f'_M)$ induce the same copy if and only if $(f'_N,f'_M)\circ(f^{-1}_N,f^{-1}_M)$ is an automorphism of $G$.
\end{proof}

Now we want to find a threshold for a fixed subgraph to appear in a random bipartite graph (thus, a random transversal matroid). For this purpose, we need to set up a probability space to work on. A random bipartite graph $G$ in model $G(n,n;p)$ arises as follows: Let $G_N,G_M$ be two disjoint sets of vertices of cardinality $n$. The probability that an edge between any vertex in $G_N$ and any vertex in $G_M$ appears is $0\leq p\leq 1$. As before, $p$ is a function of $n$. Now we are working on a probability space $(\Omega,\mathscr{A},P)$, in which $\Omega$ is the set of all bipartite graphs on the set of vertices $G_N\cup G_M$, $\mathscr{A}$ is the power set of $\Omega$, and $$P(G_0)=p^k(1-p)^{n^2-k}$$ for a bipartite graph $G_0$ with $k$ edges. This model is denoted by $G(n,n;p)$. Let $G_N,G_M$ play the roles described in Subsection \ref{G_N,G_M}. So every graph appearing in $G(n,n;p)$ corresponds to a collection of subsets of $G_N$, and thus, a transversal matroid on $G_N$.

Given a bipartite graph $G$, for each copy $G'$ of $G$ in $K_{n,n}$, let $I_{G'}$ be the indicator random variable. That is, $I_{G'}$ is a function from $\Omega$ to $\{0,1\}$. Given $G_0$ arising in $G(n,n;p)$, if $G\subseteq G_0$ then $I_{G}(G_0)=1$ and the output is 0 otherwise. And let $X_G:=\sum_{G'}I_{G'}$ where the sum is taken over all copies $G'$ of $G$ in $K_{n,n}$. Note that $I_{G'}$ are not independent variables but we do not need independence of our variables in our argument. Simply put, $X_G$ is the random variable counting the number of copies of $G$ in $G_0$. Recall that $a_n=\Theta(b_n)$ if there exists $C,c,n_0>0$ such that whenever $n>n_0$, we have $cb_n\leq a_n\leq Cb_n$.

Then 
\begin{align}\label{eq:E(X_G)}
\mathbb{E}( X_G)&=\left(\left(\begin{array}{cc}
     n  \\
     |G_N| 
\end{array}\right)\left(\begin{array}{cc}
     n  \\
     |G_M| 
\end{array}\right)|G_N|!|G_M|!\Big /|Aut(G)| \right)p^{|G_E|}\nonumber \\&=\Theta(n^{|G_V|}p^{|G_E|}).
\end{align}

This converges to 0  if $p\ll n^{-|G_V|/|G_E|}$ and converges to $\infty$ if $p\gg n^{-|G_V|/|G_E|}$. The conclusion of \eqref{eq:E(X_G)} is still true in the context of general random simple graphs, namely, $G(n;p)$. The difference is, in $G(n;p)$, the term following the first equality becomes
$$\left(\begin{array}{cc}
     n  \\
     |G_V| 
\end{array}\right)|G_V|!\Big /|Aut(G)|,$$ which does not alter what follows the second inequality.

Now let
\begin{equation}\label{eq:m(G)}
 m(G)=\max\left\{\frac{|H_E|}{|H_V|}:H\subseteq G,|H_E|>0\right\}
\end{equation}

\begin{theorem}[{\cite[p.56]{Janson}}]
    Given a bipartite graph $G$ with at least one edge,
    $$\lim_{n\rightarrow \infty}P(G\subseteq G(n,n;p))=\left\{\begin{array}{rcl}
         0 &\mbox{if}& p\ll n^{-1/m(G)} \\
         1 &\mbox{if}& p\gg n^{-1/m(G)}
    \end{array}   \right.$$
\end{theorem}
The preceding theorem is a classical result in $G(n;p)$, whose proof relies only on the conclusion of \eqref{eq:E(X_G)}. Thus it stays true in $G(n,n;p)$. And the proof is identical to that in \cite[Theorem 3.4]{Janson}.

\begin{remark}\label{Re:label free}
    In this section, we are concerned about certain subgraphs in the complete bipartite graph because we are concerned about certain submatroids induced by the subgraph. The isolated vertices of $G_M$ do not contribute to the collection of subsets of $G_N$ induced by $G$. Thus, we can erase those isolated vertices. Now, given a bipartite graph $G$, the number of copies of $G$ in the complete bipartite graph $K_{n,n}$ becomes:
    $$\left(\begin{array}{cc}
         n  \\
         |G_N| 
    \end{array}\right)\left(\begin{array}{cc}
         n  \\
         |G_M^c| 
    \end{array}\right)|G_N|!|G_M^c|!\Big /|Aut(G)|.$$ That does not change what follows the second inequality of equation \eqref{eq:E(X_G)}. And note that in equation \eqref{eq:m(G)}, the subgraph $H\subseteq G$ inducing the maximum of $m(G)$ must not include any vertex of $G_M\setminus G_M^c$. Therefore, the above theorem is still valid. 
\end{remark}     

The reason why we are concerned about the thresholds for a given graph to appear in $G(n,n;p)$ is that if $T(H)$ is a uniform matroid and $H\subseteq G$, then we can perform restriction followed by truncation on $T(G)$ to get $T(H)$, where truncation is introduced in the following theorem:

\begin{theorem}[{\cite[p.59]{Welsh}}]
    Let $M$ be a matroid on a finite set $E$ and $k\leq \rho(M)$ be a positive integer. Let $$\mathscr{I}_k:=\{X\in \mathscr{I}(M):|X|\leq k\}.$$ Then $\mathscr{I}_k$ is the collection of independent sets of a matroid $M_k$ on $E$.
\end{theorem}
\noindent $M_k$ is called the \textit{truncation} of $M$ at $k$.

\begin{theorem}\label{thm:truncation}
    Let $T(H)$ be a uniform matroid. If $H\subseteq G$ in a bipartite graph $G$, the matroid obtained from restricting $T(G)$ to $H_N$ and performing truncation at $\rho(T(H))$ is $T(H)$.
\end{theorem}
\begin{proof}
    Clearly, $T(G)|\{H_N\}$ is a matroid on the set $H_N$, and $\rho(T(G)|\{H_N\})\geq \rho(T(H))$. Then the collection of bases of the truncation of $T(G)|\{H_N\}$ at $\rho(T(H))$ is the same as that of $T(H)$.
\end{proof}

Unfortunately, if $H\subseteq G$ but $T(H)$ is not a uniform matroid, as the following example suggests, the conclusion of the above theorem may not hold.

\begin{figure}[h]
    \centering
    \includegraphics[width=0.2\linewidth]{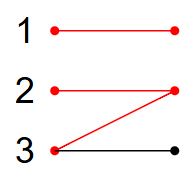}
    \caption{We may not get $T(H)$ through truncation from $T(G)$ even if $H$ is a subgraph of $G$.}
    \label{fig:no truncation}
\end{figure}

\begin{example}
    Let $G$ be the graph in Figure \ref{fig:no truncation}, $G_N=\{1,2,3\}$, and $H$ be the subgraph labeled red. It is easily verified that no truncation of $T(G)$ can result in $T(H)$.
\end{example}

\begin{figure}[h]
    \centering
    \includegraphics[width=0.6\linewidth]{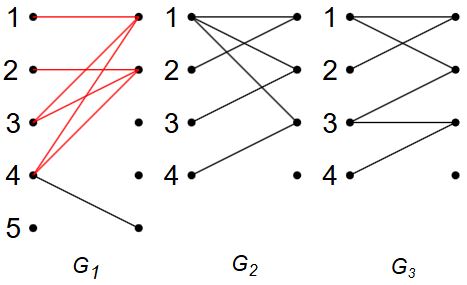}
    \caption{``$G$ contains $H$ as a subgraph'' is neither sufficient nor necessary for $T(G)$ to contain $T(H)$ as a minor.}
        \label{fig:subgraph minors iff nonexample}
\end{figure}

It is natural to ask whether a bipartite graph $G$ contains $H$ as a subgraph implies that $T(G)$ contain $T(H)$ as a minor? Figure \ref{fig:no truncation} and $G_1$ in Figure \ref{fig:subgraph minors iff nonexample} provide counterexamples. Denote the red subgraph of $G_1$ in Figure \ref{fig:subgraph minors iff nonexample} as $H$. It is straightforward to see that $T(H)$ is $U_{2,4}$, and $H$ is the unique (up to isomorphism) minimal presentation of $U_{2,4}$. But $T(G_1)$ does not contain $U_{2,4}$ as a minor.  But is it necessary that $G$ contains a subgraph $H$ for $T(G)$ to contain $T(H)$ as a minor? $G_2$ and $G_3$ in Figure \ref{fig:subgraph minors iff nonexample} provide a counterexample. $G_2$ and $G_3$ are both minimal presentations for $U_{3,4}$, and certainly they both do not contain each other as a subgraph.

We give sufficient conditions for $T(G)$ to contain $T(K)$ as a minor, $K\subseteq G$, in the following theorem:

\begin{theorem}\label{thm:sufficient U_2,4}
    If $G$ satisfies the following conditions, then $T(G)$ contains $T(K)$ as a minor:
    \begin{enumerate}
        \item $G$ contains $K$ as a subgraph;
        \item Every vertex in $G_M\setminus K_{M}$ which is adjacent to some vertex in $K_N$ is also adjacent to at least one vertex in $G_N\setminus K_N$;
        \item  If $\{x_1,\ldots,x_k\}$ is the complete list of vertices in $G_M\setminus K_M$ that are adjacent to some vertex (not necessarily the same) of $K_N$, then there exist vertices $\{y_1,\ldots,y_k\}$ in $G_N\setminus K_N$ such that $y_i$ is adjacent to $x_i$ but to no other vertices in $G_M$.
    \end{enumerate}
\end{theorem}
\begin{proof}
    Restrict $T(G)$ to $K_N\cup\{y_1,\ldots,y_k\}$. Note that $\{k_1,\ldots,k_r\}\subseteq E(T(K))$ is independent in $T(K)$ if and only if $\{k_1,\ldots,k_r\}$ can be matched by some subset of $K_M$. By condition 3, $\{k_1,\ldots,k_r\}\cup \{y_1,\ldots,y_k\}$ is independent in $T(G)$. Then perform contraction of $\{y_1,\ldots,y_k\}$ from $T(G)$. The result is $T(K)$. 
\end{proof}

\begin{figure}[h]
    \centering
    \includegraphics[width=0.24\linewidth]{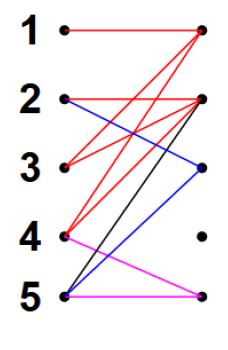}
    \caption{$T(G)$ does not contain $U_{2,4}$ as a minor.}
    \label{U_2,4 minor non example}
\end{figure}

\noindent Figure \ref{U_2,4 minor non example} serves as a nonexample. Let $G$ be the graph in the figure, and $G_N=\{1,2,3,4,5\}$. The red subgraph is $K$, $T(K)=U_{2,4}$, and $U_{2,4}$ is not a minor of $T(G)$. $G$ violates condition (3) of the previous theorem. 

\begin{example}
    We give an example to illustrate the process performed in the proof of Theorem \ref{thm:sufficient U_2,4}. In Figure \ref{fig:U_2,4 minor process}, $G$ is the leftmost subfigure with $G_N=[8]$, $K$ is the subgraph marked as red, and $\{5,6,7\}$ is $\{y_1,y_2,y_3\}$ in the theorem. The middle subfigure represents $T(G)|[7]$, and the last one represents $T(G)|[7]/\{5,6,7\}$.

\end{example}

\begin{figure}[h]
    \centering
    \includegraphics[width=0.8\linewidth]{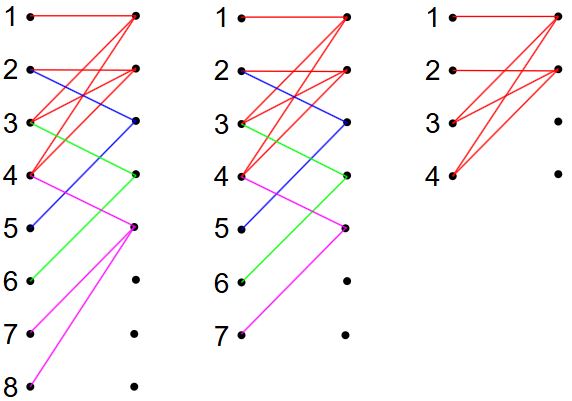}
    \caption{An example illustrating Theorem \ref{thm:sufficient U_2,4}.}
    \label{fig:U_2,4 minor process}
\end{figure}

Given a transversal matroid $T$, is there a fixed bipartite graph $H_T$ such that whenever a bipartite graph $G$ contains $H_T$ as a subgraph, then $T(G)$ contains $T$ (or does not contain $T$) as a minor? We answer these questions in the following two propositions:

\begin{proposition}
    Given a transversal matroid $\mathscr{T}$ which is not free\footnote{A matroid is called free if every subset of its ground set is independent.}, there is no fixed bipartite graph $H$ such that whenever a bipartite graph $G$ contains $H'$ as a subgraph, where $H'\cong H$, then $T(G)$ contains $\mathscr{T}$ as a minor. 
\end{proposition}
\begin{proof}
    Suppose there exists such a bipartite graph $H$ for a transversal matroid $\mathscr{T}$ and $G$ contains $H'$ as a subgraph. Then we may form a bipartite graph $G'$ by adding edges to $G$ such that $T(G')$ is a free matroid. But any minor of a free matroid is also free. Therefore it is impossible for $T(G')$ to contain $\mathscr{T}$ as a minor. That is a contradiction.
\end{proof}

\begin{proposition}
    Given a transversal matroid $\mathscr{T}$, there is no fixed bipartite graph $H$ such that whenever a bipartite graph $G$ contains $H$ as a subgraph, then $T(G)$ does not contain $\mathscr{T}$ as a minor. 
\end{proposition}
\begin{proof}
    Suppose there exists a bipartite graph $H$ for a transversal matroid $\mathscr{T}$. Since $\mathscr{T}$ is transversal, there exists a bipartite graph $K$ such that $T(K)=\mathscr{T}$. Form a bipartite graph $G$ containing $H$ as a subgraph as follows: let $|G_N|=|G_M|$ be large enough such that $G_N$ contains $H_N,K_N$ as subsets with $H_N\cap K_N=\emptyset$ in $G_N$, and $G_M$ contains $H_M,K_M$ as subsets with $H_M\cap K_M=\emptyset$ in $G_M$. Moreover, let $G$ inherit all the adjacencies of $H$ and $K$. By construction, $T(G)|K_N=\mathscr{T}$, a contradiction.
\end{proof}

\begin{proposition}
    Given a transversal matroid $\mathscr{T}$, there is no fixed bipartite graph $H$ such that whenever a bipartite graph $G$ contains $H'$ as a subgraph, where $H'\cong H$, then $T(G)$ does not contain $\mathscr{T}$ as a minor. 
\end{proposition}
\begin{proof}
    Suppose there exists a bipartite graph $H$ for a transversal matroid $\mathscr{T}$. Since $\mathscr{T}$ is transversal, there exists a bipartite graph $K$ such that $T(K)=\mathscr{T}$. Form a bipartite graph $G$ containing $H$ as a subgraph as follows: let $|G_N|=|G_M|$ be large enough such that $G_N$ contains $H_N,K_N$ as subsets with $H_N\cap K_N=\emptyset$ in $G_N$, and $G_M$ contains $H_M,K_M$ as subsets with $H_M\cap K_M=\emptyset$ in $G_M$. Moreover, let $G$ inherit all the adjacencies of $H$ and $K$. By construction, $T(G)|K_N=\mathscr{T}$, a contradiction.
\end{proof}

Matroids containing $U_{2,4}$ as a minor are of particular significance\footnote{Let $A$ be a $m\times n$ matrix over a field $\mathbb{F}$, the collection of subsets of the column vectors of $A$ each of which is an independent set in $V(m,\mathbb{F})$ forms the independent sets of a matroid, called the vector matroid of $A$.A matroid is called binary if it is a vector matroid of a matrix over $GF(2)$. A matroid is binary iff it has no $U_{2,4}$ minor\cite{Tutte1}}. So we are going to say some more about that. Let us bring the concept of lattice up to our mind. A \textit{lattice} is a partially ordered set $(\mathscr{L},\leqslant)$ such that for any $x,y\in \mathscr{L}$, there exists a unique greatest lower bound $x\wedge y\leqslant x,y$ (\textit{meet} of $x,y$), and a unique least upper bound $x\vee y\geqslant x,y$ (\textit{join} of $x,y$), such that if $z\leqslant x$ and $z\leqslant y$ then $z\leqslant x\wedge y$, and if $x\leqslant z$ and $y\leqslant z$ then $x\vee y\leqslant z$. $x$ is said to \textit{cover} $y$ in $\mathscr{L}$ if $y< x$ and whenever $y< z< x$ then $z=x$. In a lattice $\mathscr{L}$ with a universal lower bound $\mathscr{L}_O$, the elements covering $\mathscr{L}_{\mathscr{O}}$ are called \textit{atoms} or \textit{points}. In a lattice $\mathscr{L}$ with a universal upper bound $\mathscr{L}_I$, the elements covered by $\mathscr{L}_I$ are called \textit{coatoms} or \textit{copoints}. A finite lattice $\mathscr{L}$ is \textit{semimodular} if for all $x,y\in \mathscr{L}$, whenever $x$ and $y$ cover $x\wedge y$, then $x\vee y$ covers $x$ and $y$. A finite lattice is called \textit{geometric} if it is semimodular and every element is the join of atoms. Moreover, if a partially ordered set has a universal lower bound $\mathscr{L}_{\mathscr{O}}$, the \textit{height} $h(x)$ of $x\in P$ is the least upper bound of lengths of chains $\mathscr{L}_{\mathscr{O}}=x_0<x_1<...<x_k=x$ between $\mathscr{L}_{\mathscr{O}}$ and $x$. And a finite poset $P$ is said to satisfy the \textit{Jordan-Dedekind chain condition} if given arbitrary $x,y\in P$, all maximal chains between them have the same length. Note that a semimodular lattice satisfies the Jordan-Dedekind chain condition and its height $h$ satisfies $h(x\vee y)+h(x\wedge y)\leq h(x)+h(y)$.

Now recall Definition \ref{def:flat} for the concepts of flats in a matroid. Declare a partial order on the collection of all flats of a matroid by set inclusion, then the collection forms a geometric lattice (see \cite[p.50]{Welsh} for an argument). The lattice of flats of a matroid $M$ is denoted by $\mathscr{L}(M)$. 

\begin{theorem}[{\cite[p.169]{Welsh}}]
    A geometric lattice $\mathscr{L}$ of height $r$ is the lattice of a matroid that does not contain $U_{2,4}$ as a minor if and only if every element of height $r-2$ is contained in at most 3 copoints.
\end{theorem}

\begin{example}
    \begin{figure}[h]
        \centering
        \includegraphics[width=0.8\linewidth]{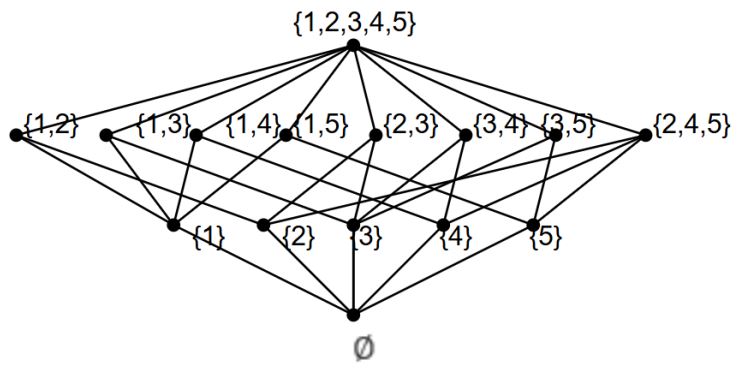}
        \caption{Hasse diagram of the lattice of flats of a matroid.}
        \label{lattice of flats example}
    \end{figure}
    Figure \ref{lattice of flats example} is the Hasse diagram of the lattice of flats of $T(G)$, where $G$ is from the figure below. From the previous theorem, $T(G)$ contains $U_{2,4}$ as a minor because in its lattice of flats, there is an atom which is contained in more than three copoints. 
\end{example}

\begin{figure}[h]
    \centering
    \includegraphics[width=0.2\linewidth]{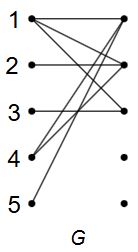}
    \caption{$T(G)$ contains $U_{2,4}$ as a minor.}
    \label{fig:contains U_2,4 as minor serve example for lattice of flats}
\end{figure}

\begin{figure}[h]
    \centering
    \includegraphics[width=0.65\linewidth]{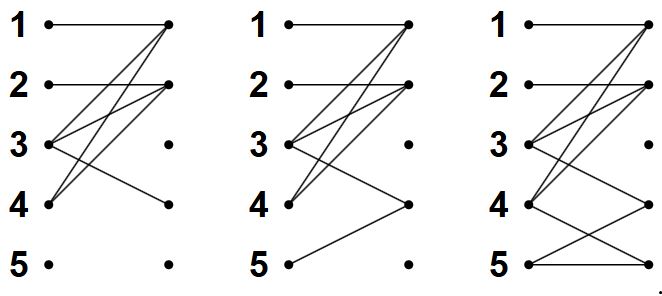}
    \caption{``$T(G)$ is binary'' is not a convex property.}
    \label{fig:binary not convex}
\end{figure}

\begin{remark}
Figure \ref{fig:binary not convex} shows that ``$T(G)$ is binary" is not even a convex property. Let the graphs in the figure be $G_1\subseteq G_2\subseteq G_3$. Both $T(G_1)$ and $T(G_3)$ do not possess $U_{2,4}$ as a minor, but $T(G_2)/5\cong U_{2,4}$.
\end{remark}

\section{Simple $T(G)$}\label{sec:simple T(G)}
We bear this stipulation in mind: Whenever we have a bipartite graph $G$, if we add an edge to $G$, we keep the bipartition of $G$. That is, the new graph formed by adding edges to $G$ has the same bipartition as $G$.

Recall that a matroid is called simple if it has no loops or parallel edges (a parallel edge means a circuit of cardinality two). In the context of transversal matroids, $T(G)$ has no loops if and only if $G_N$ has no isolated vertex, and this is clearly a monotone increasing property. And $T(G)$ has no parallel edges if and only if there are no two elements $x,y\in G_N$ such that both of them are adjacent to some $z$ and $z$ is the only vertex they are adjacent to. This is not a monotone increasing property. The problem arises when there is an $x\in G_N$ which is adjacent to only one vertex $z$, and there is an isolated vertex $y\in G_N$. After adding an edge joining $y$ to $z$, $\{x,y\}$ becomes a parallel edge. Clearly, if we require that $G_N$ have no isolated vertices, this situation does not happen. Therefore, ``$T(G)$ is simple" is a monotone increasing property.

\begin{theorem}
    In $G(n,n;p), \frac{\log n}{n}$ is a threshold for $T(G)$ to have no loop.
\end{theorem}
\begin{proof}
    Let $I_x,x\in G_N$ be independent random variables such that $I_x(G)=1$ if $x$ is a loop in $T(G)$ and $I_x(G)=0$ otherwise. Let $X$ be the random variable $\sum_{x\in G_N}I_x$. $P(I_x(G)=1)=(1-p)^n$ for all $x\in G_N$. We have $$\mathbb{E}(X)=n(1-p)^n$$ If $p\gg \frac{\log n}{n}$, then there is a sequence $\{s_n\}\rightarrow \infty$ such that $p(n)=\frac{\log n+s_n}{n}$ and since $0\leq p\leq 1$, we have $1-p\leq e^{-p}$ as $n\rightarrow \infty$
    $$\mathbb{E}(X)\leq ne^{-pn}\leq ne^{-\log n-s_n}=e^{-s_n}=o(1).$$ On the other hand, $p\ll \frac{\log n}{n}$ implies the existence of a sequence $\{s_n\}\rightarrow \infty$ such that $p(n)=\frac{\log n-s_n}{n}$, 
    $$\mathbb{E}(X)=n(1-p)^n=n \left(1+\frac{s_n-\log n}{n}\right)^n=ne^{s_n-\log n+o(1)}=e^{s_n+o(1)}$$ as $n\rightarrow \infty$.

    Note that 
    \begin{dmath*}
        \operatorname{Var}(X)=\mathbb{E}(X^2)-(\mathbb{E}(X))^2=n(1-p)^n+n(n-1)(1-p)^{2n}-n^2(1-p)^{2n}=\mathbb{E}(X)(1+(n-1)(1-p)^n-n(1-p)^2)=\mathbb{E}(X)(1-(1-p)^n).
    \end{dmath*}
    By inequality \eqref{eq:1.2}, $P(X>0)\geq 1-\frac{\operatorname{Var}X}{(\mathbb{E}(X))^2}=1-\frac{1-(1-p)^n}{\mathbb{E}(X)}\rightarrow 1$ as $n\rightarrow \infty$.
\end{proof}
The following result is analogous to the classical result that $\frac{\log n}{n}$ is a threshold for random graphs to have no isolated vertex in $G(n,p)$. It can also be found in \cite{Erdos5}.

\begin{proposition}
    In $G(n,n;p)$, the probability that $T(G)$ contains a parallel edge goes to 0 if $p\gg \frac{\log n}{n}$. 
\end{proposition}
\begin{proof}
    Let $I_{x,y},x,y\in G_N$ be the indicator random variable such that $I_{x,y}(G)=1$ if $\{x,y\}$ is a circuit of $T(G)$ and $I_{x,y}(G)=0$ otherwise. Let $X=\sum_{\{x,y\}\in G_N}I_{x,y}$, and $s_n$ be as in the preceding theorem. 
    \begin{dmath*}
        \mathbb{E}(X)=\left(\begin{array}{cc}
             n  \\
             2 
        \end{array}\right)np^2(1-p)^{2(n-1)}\leq \left(\begin{array}{cc}
             n  \\
             2 
        \end{array}\right)np^2e^{-2p(n-1)}=\Theta(\frac{n^3}{2}p^2e^{-2\log n-2s_n})=\Theta(\frac{n}{2}p^2e^{-s_n})\rightarrow 0 
    \end{dmath*}
    as $n\rightarrow \infty$.
\end{proof}

\begin{corollary}
    In $G(n,n;p), p(n)=\frac{\log n}{n}$ is a threshold for $T(G)$ to be simple. 
\end{corollary}
\begin{proof}
    The definition of $T(G)$ being simple is that $T(G)$ has no loops or parallel edges. So the result follows from the preceding two results.
\end{proof}

\chapter{Non Monotonicity}
In this chapter, we look at non-monotone properties of transversal matroids such as how circuits change, the decomposability of the lattice of flats, and the number of circuits. Then we investigate some characterizations of circuits of transversal matroids and how circuits `paste' together in a bipartite graph for a connected transversal matroid. 
\section{Dual of transversal matroids}\label{sec:T(G)^*}
In this section we investigate how the number of circuits of $T(G)^*$ changes as we add an edge to $G$ (note that the number of bases of $T(G)$ is always equal to that of $T(G)^*$). We need to introduce the concept of \textit{strict gammoid}, and we mainly follow the notations in \cite[Chapter 13, Sect. 1,2]{Welsh} in the context in strict gammoids. 

Let $D=(V,E)$ be a directed graph (abbrev. as digraph) with vertex set $V$ and edge set $E$. A \textit{path} in $D$ is a sequence $(v_0,...,v_k)$ of distinct vertices with $k\geq 0$, and $(v_{i-1},v_i)$ is an edge of $E$, $1\leq i\leq k$. Two paths in $D$ are \textit{disjoint} if their vertex sets are disjoint. If $A$and $B$ are subsets of $V$ and there exists a bijection $\gamma:A\rightarrow B$ so that we can find disjoint paths $\{P_v:v\in A\}$ in $D$ such that $P_v$ has initial vertex $v$ and terminal vertex $\gamma(v)$, then we say that $A$ can be linked onto $B$. And we say that $A$ can be linked into $C$ if $A$ can be linked onto a subset of $C$. For $Z\subseteq V$, let $\tilde Z:=Z\cup\{v\in V:(z,v)\in E \mbox{ for some } z\in Z\}$. And for $v\in V$, let $\mathscr{A}_v:=\widetilde{\{v\}}$, $\mathscr{A}$ be the family of sets $\{\mathscr{A}_v:v\in V\}$ and $\mathscr{A}_X$ be the subfamily $\{\mathscr{A}_v:v\in X\},X\subseteq V$. Choose a subset $B\subseteq V$, let $L(D,B)$ denote the collection of subsets of $V$ which can be linked into $B$.

The following lemma is called the \textit{linkage lemma} in \cite[Chapter 13]{Welsh}. Essentially, the content of the Lemma is due to Ingleton and Piff \cite{Ingleton}. D.J.Welsh also proved it in \cite[p.218]{Welsh}. 

\begin{lemma}
    Let $X,Y$ be subsets of the vertex set $V$ of a digraph $D$. Then $X$ can be linked onto $Y$ in $D$ if and only if $V\setminus X$ is a transversal of the family of sets $\mathscr{A}_{V\setminus Y}$.
\end{lemma}

The following theorem is due to H. Perfect. D. J. Welsh also proved it using the linkage lemma in \cite[p.220]{Welsh}. We copy the proof from Welsh here because we will have something important to note from there along our development below, and his proof also proves Theorem \ref{thm: strict gam. M*} stated below which is due to A. A. W. Ingleton and M. J. Piff.
\begin{theorem}{\cite{Perfect}}\label{thm:L(D,B)}
    $L(D,B)$ is the collection of independent sets of a matroid on $V$, where $V$ is the vertex set of a digraph $D$ and $B$ is a fixed subset of $V$.
\end{theorem}
\begin{proof}
    By the linkage lemma, if $X\subseteq V$, $X$ can be linked into $B$ if and only if $X$ can be linked onto $B'\subseteq B$, this is to say that $V\setminus X$ is a transversal of the family $\mathscr{A}_{V\setminus B'}$. Thus $X\in L(D,B)$ if and only if $V\setminus X$ contains a transversal of $\mathscr{A}_{V\setminus B}$. This is to say that $X\in L(D,B)$ if and only if $V\setminus X$ is spanning in the transversal matroid $M[\mathscr{A}_{V\setminus B}]$. Note that the complement of a spanning set of a matroid $M$ is independent of $M^*$. 
\end{proof}

A matroid that can be obtained from a digraph $D$ and some subset $B$ of the vertex set of $D$ in the above way is called a \textit{strict gammoid}.

\noindent As we said, the proof above also proves the following theorem which is due to Ingleton and Piff:

\begin{theorem}{\cite{Ingleton}}\label{thm: strict gam. M*}
    A matroid $M$ is a strict gammoid if and only if its dual matroid $M^*$ is transversal.
\end{theorem}

Based on the proof of theorem \ref{thm:L(D,B)}, given a transversal matroid $M[\mathscr{A}]$, D. J. Welsh wrote a procedure to construct $(M[\mathscr{A}])^*$. It is described as follows:

\begin{construction}\label{con:dual of trans.}
Given $M[\mathscr{A}]$ on a set $V$, select a basis $B$ of $M[\mathscr{A}]$, say $B=\{x_1,\ldots,x_r\}$. Recall from subsection \ref{G_N,G_M}, $\mathscr{A}$ corresponds to a bipartite graph $G$ such that $G_N=V$, $m\in G_M$ is adjacent to all vertices in $A_m$. Pick $J=\{j_1,\ldots,j_r\}\subseteq G_M$ such that $j_i$ is matched to $x_i$, $1\leq i\leq r$. Then form a digraph $D^J_B(M[\mathscr{A}])$ by joining $x_i$ to all $y\in A_{j_i}\setminus x_i$ with vertex set $V$. And $L(D^J_B(M[\mathscr{A}]),V\setminus B)=(M[\mathscr{A}])^*$.
\end{construction}
Now there may be an ambiguity coming up to our mind. According to Theorem \ref{thm:L(D,B)}, in construction \ref{con:dual of trans.}, $L(D^J_B(M[\mathscr{A}]),V\setminus B)=(M[\mathscr{A}_B])^*$. But the digraph $\mathscr{A}_B$ depends on the choice of $B$ and $J$. How can we be sure that it is always equal to $(M[\mathscr{A}])^*$? The following lemma clears up this mist:

\begin{lemma}[{\cite[Brualdi]{Brualdi}}]\label{lemma: exists |G_M|=|base|}
    Let $\mathscr{A}=\{A_i:i\in [r]\}$ be a presentation of a transversal matroid $M$. If $\phi:B\rightarrow [r]$ is a matching of a basis $B$ of $M$ into $\mathscr{A}$, then $\{A_i:i\in \phi(B)\}$ is also a presentation of $M$.
\end{lemma}

This implies that $M[\mathscr{A}_B]=M[\mathscr{A}]$ in construction \ref{con:dual of trans.} for every choice of basis $B$ and matching $J$.

\begin{lemma}
    Let $G$ be a bipartite graph, and $G'$ be a bipartite graph with the same bipartition as $G$, obtained by adding one edge to $G$. If $\rho (T(G'))=\rho (T(G))+1$ and the new edge connects $m\in G_M$ to $n\in G_N$, then there exists a basis $B$ of $T(G)$ such that $n\notin B$. 
\end{lemma}
\begin{proof}
    Suppose $B_1$ is a basis of $T(G)$ such that $n\in B_1$. Then there is a $J\subseteq G_M$ matching to $B_1$ in $G$. Now, since $\rho(T(G'))=\rho(T(G))+1$, there is an element $s\in G_N$ such that $B_1\cup \{s\}$ is a basis for $T(G')$. That means $B_1\cup \{s\}$ can be matched to $J\cup \{r\}$ for some $r\in G_M$. But all edges connecting $s$ in $G'$ are also in $G$. Therefore, $(B_1\cup \{s\})\setminus \{n\}$ can also be matched into $J\cup \{r\}$ in $G$. That means $(B_1\cup \{s\})\setminus \{n\}$ is also independent, and thus a basis, in $T(G)$.
\end{proof}

\begin{theorem}
    Let $G$ be a bipartite graph, and $G'$ be a bipartite graph obtained by adding one edge to $G$ connecting $m\in G_M$ to $n\in G_N$. If $\rho (T(G'))=\rho (T(G))+1$, let $B$ be a basis of $T(G)$ that satisfies the previous lemma and $J$ a matching of $B$ in $G$ such that $m\notin J$, then we have the following:
    \begin{enumerate}[label=(\roman*)]
        \item If $V$ is a circuit of $T(G)^*$ and $v\in V$ can be linked only to $n$ in $D_B^J(T(G))$, then both $v$ and $V\setminus \{v\}$ are circuits of $T(G')^*$.
        \item If $V$ is a circuit of $T(G)^*$ and there exists a $v\in V$ such that $v$ can be linked to $n$ in $D_B^J(T(G))$, then $V\setminus \{r\}$ is dependent in $T(G')^*$ for all $r\in V$.
        \item If $V$ is a circuit of $T(G)^*$ and no element of $V$ can be linked to $n$ in $D_B^J(T(G))$, then $V$ stays as a circuit in $T(G')^*$.
    \end{enumerate}
\end{theorem}
\begin{proof}
    Pick a basis $B=\{n_1,\ldots,n_t\}$ of $T(G)$ such that $n\notin B$. First, we show that we may pick a $J$ specified in the assumption of our theorem, that is, $\{m_1,...,m_t\}=J\subseteq G_M$ matching to $B$ and $m\notin J$, say $m_i$ is matched to $n_i,1\leq i\leq t$ (Such $J$ exists since $B\cup \{n\}$ is a basis in $T(G')$, there is a $J'\subseteq G_M'$ matching to $B\cup \{n\}$. Note that $m\in J'$, otherwise $\rho(T(G))=t+1$, a contradiction. In this matching in $G'$, $m$ has to be matched to $n$. Otherwise, this matching still works in $G$, again that would make $\rho(T(G))=t+1$. $J'\setminus \{m\}$ is the $J$ we want). Note that in the particular case $G_N\setminus B=\{n\}$, we get a stronger result, namely, every single element of $G_N$ becomes a circuit in $T(G')^*$. Indeed, in this case, if $V$ is a circuit of $T(G)^*$ that satisfies the assumption of $(i)$ or $(ii)$, $V$ must be of cardinality 2, then clearly the conclusions of $(i)$ or $(ii)$ hold, respectively. If $V$ is a circuit of $T(G)^*$ that satisfies the assumption of $(iii)$, $V$ must be a singleton and the conclusion of $(iii)$ obviously holds. 
    
    Now we return to the general case, note that $m$ cannot be connected to any vertex in $G_N\setminus B$, since that would again make $\rho(T(G))=t+1$. For the same reason, there can not be a path $(m,n_{i_1},m_{i_1},\ldots,n_{i_k},m_{i_k},n')$, where $\{n_{i_1},\ldots,n_{i_k}\}$ (resp. $\{m_{i_1},\ldots,m_{i_k}\}$) is a subset of $\{n_1,\ldots,n_t\}$ (resp. $\{m_1,\ldots,m_t\}$), and $n'\in G_N\setminus B$ (a path in a graph is a sequence of distinct vertices such that there is an edge in the graph joining any two consecutive vertices).

    $J'=J\cup \{m\}$ is a matching for $B\cup \{n\}$ in $G'$.  Now we compare $D^{J'}_{B\cup \{n\}}(T(G'))$ to $D^J_B(T(G))$. The only difference between them is that $n$ can become a source vertex in $D^{J'}_{B\cup \{n\}}(T(G'))$. But by our discussion above, that cannot create a path from $n$ to any vertex in $G_N\setminus B\cup \{n\}$ in $D^{J'}_{B\cup\{n\}}(T(G'))$. Suppose $V$ is a circuit in $T(G)^*$, that means $V$ is minimal with respect to the fact that $V$ cannot be linked into $G_N\setminus B$. 
    
    Statement $(i)$ is straightforward, since $V$ cannot be linked to $G_N\setminus B$ in $D_B^J(T(G))$, $V$ is not a singleton if it satisfies the assumption of $(i)$, so $V\setminus\{v\}$ is nonempty, and it cannot be linked to $G_N\setminus B\cup \{n\}$ in $D_B^J(T(G))$, and thus in $D_{B\cup \{n\}}^{J'}(T(G'))$. To show $(ii)$, suppose $v\in V$ can be linked to $n$ in $D_B^J(T(G))$, then again $V\setminus \{v\}$ is nonempty and cannot be linked to $G_N\setminus B\cup \{n\}$ in $D_B^J(T(G))$, and thus it cannot be linked to $G_N\setminus B\cup \{n\}$ in $D_{B\cup \{n\}}^{J'}(T(G'))$. If $w\in V$ such that $w$ cannot be linked to $n$ in $D_B^J(T(G))$, we want to show that $V\setminus \{w\}$ cannot be linked to $G_N\setminus B\cup \{n\}$ in $D_B^J(T(G))$, and thus also in $D_{B\cup \{n\}}^{J'}(T(G'))$. Suppose otherwise, there is a linking from $V\setminus \{w\}=\{v,x_1,\ldots,x_m\}$ to $\{y_v,y_1,\ldots,y_m\}\subseteq G_N\setminus B\cup \{n\}$ in $D_B^J(T(G))$, where $v$ is linked to $y$, $x_i$ is linked to $y_i,1\leq i\leq m$, and $n\notin \{y_v,y_1,\ldots,y_m\}$. Observe that $w$ cannot be linked to any element in $G_N\setminus B\cup \{n\}$ other than an element in $\{y_v,y_1,\ldots,y_m\}$. If such a case occurs, then $V$ can be linked to $G_N\setminus B\cup \{n\}$ and thus to $G_N\setminus B$, contradicting our assumption. But $w$ can be linked to some element in $G_N\setminus B\cup \{n\}$, since otherwise $w$ itself would be a circuit of $T(G)^*$, which is absurd. Thus, $w$ has to be linked to at least one of $\{y_v,y_1,\ldots,y_m\}$:
    \\Case 1. If $w$ is linked to $y_v$, since $v$ is linked to $n$, $V$ can be linked to $\{n,y_v,y_1,\ldots,y_m\}\subseteq G_N\setminus B$. This is a contradiction. 
    \\Case 2. If $w$ cannot be linked to $y_v$, then without loss of generality, we may assume that $w$ is linked to $y_m$. Again, $x_m$ cannot be linked to any element in $G_N\setminus B\cup \{n\}$ other than the elements in $\{y_v,y_1,\ldots,y_m\}$. But if $w$ and $x_m$ can be linked only to $y_m$, then $\{x_m,w\}$ would be a circuit of $T(G)^*$, contradicting our assumption. Again, without loss of generality, assume $x_m$ is linked to $y_{m-1}$. Continuing this process, it can be assumed that $w$ is linked to $y_m$, $x_{i}$ is linked to $y_{i-1},1<i\leq m$, $x_1$ is linked to $y_v$ in $D_B^J(T(G))$. But recall that $v$ can be linked to $n$, thus $V$ can be linked to $G_N\setminus B$ in $D_B^J(T(G))$. This is a contradiction. 

Statement $(iii)$ is obvious.
\end{proof}

\begin{theorem}
    Let $G$ be a bipartite graph, and $G'$ be a bipartite graph with the same bipartition as $G$, obtained by adding one edge to $G$. If $\rho(T(G'))=\rho(T(G))$, suppose the new edge joins $m\in G_M$ to $n\in G_N$, and $V$ is a circuit of $T(G)^*$, then either $V$ is a circuit of $T(G')^*$, or $U\cup \{n\}$ is a circuit of $T(G')^*$ for some $U\subseteq V$.
\end{theorem}
\begin{proof}
    Suppose the new edge joins $m\in G_M$ to $n\in G_N$. If $\mathscr{B}(T(G))=\mathscr{B}(T(G'))$, then there is nothing to prove. Thus, we assume that their sets of bases are not equal. Then, since any basis of $T(G)$ stays as a basis in $T(G')$, there must be some basis of $T(G')$ that is not a basis of $T(G)$. First we show that there must be a basis of $T(G)$ not containing $n$. Suppose otherwise, every basis of $T(G)$ contains $n$. Pick a basis $B$ of $T(G)$, it is also a basis of $T(G')$. And pick $J\in G_M$ matching $B$ in $G$, note that $m$ cannot be matched to $n$ in $G$. Then we can use the same matching $J$ for $B$ in $G'$, where still, $m$ cannot be matched to $n$. Now we compare $D^J_B(T(G'))$ with $D^J_B(T(G))$. If $m$ is not adjacent to any element in $B$ in $G$, then $D^J_B(T(G'))=D^J_B(T(G))$ and therefore, $T(G')^*$ and $T(G)^*$ share the same set of circuits. But that means they are the same matroids, which is absurd since the set of bases of $T(G)$ is assumed to be a proper subset of that of $T(G')$, and thus the set of bases of $T(G)^*$ is a proper subset of $T(G')^*$. Now, if $m$ is adjacent to some element of $B$ in $G$, then the only possible new directed edge (if there is one) in $D^J_B(T(G'))$ is from some element of $B$ to $n$. Because we assumed that every basis of $T(G)$ contains $n$, that implies $n$ is a coloop of $T(G)$ (a loop of $T(G)^*$). That means there is no path from $n$ to any vertex in $G_N\setminus B$ in $D^J_B(T(G))$. And there is still no path from $n$ to any vertex in $G'_N\setminus B$. That means a set is dependent in $T(G)^*$ if and only if it is dependent in $T(G')^*$ and if $V$ is a circuit of $T(G)^*$, then $V$ is still a circuit in $T(G')^*$. If $V'$ is a circuit of $T(G')^*$, clearly $V'$ is dependent in $T(G)^*$. But if $V'\setminus U$ is a circuit of $T(G)^*$, where $U$ is a proper subset of $V'$, then $V'\setminus U$ is dependent in $T(G')^*$, which is absurd. So, again, we reach the conclusion that $T(G)^*=T(G')^*$, which is a contradiction.

    Therefore, there exists a basis $B$ of $T(G)$ such that $n\notin B$ with $J\subseteq G_M$ matched to $B$. Note that $m$ must be in $J$, otherwise $\rho(T(G'))=\rho(T(G))+1$, contradicting our hypothesis. Again, we compare $D^J_B(T(G'))$ with $D^J_B(T(G))$. Since $n\notin B$, it is not a source vertex in both $D^J_B(T(G))$ and $D^J_B(T(G'))$. The only difference between these two directed graphs is that there is a new directed edge from some element of $B$ to $n$ in $D^J_B(T(G'))$, where $n$ is only a terminal vertex. 
    
    \noindent Case 1) If $V$ is a circuit of $T(G)^*$, $n\in V$, then $V$ is still a circuit in $T(G')^*$ since in $D_B^J(T(G'))$, there is no new path from any vertex in $G_N$ to $G_N\setminus B$ when compared to $D_B^J(T(G))$.

    \noindent Case 2) If $V$ is a circuit of $T(G)^*$, $n\notin V$. If there is no path from any element of $V$ to $n$ in $D^J_B(T(G'))$, then $V$ remains a circuit of $T(G')^*$. $V$ is certainly dependent in $T(G')^*$. If $U\subsetneq V$ were dependent in $T(G')^*$, then it cannot be matched to $G'_N\setminus B$, and it also cannot be matched to $G_N\setminus B$, a contradiction. 

    If there is a path from some element of $V$ to $n$ in $D_B^J(T(G'))$, then there are two subcases:
    \\ Subcase 1) If $V$ is still dependent in $T(G')^*$, then $V$ stays a circuit in $T(G')^*$. Since again, if $U\subsetneq V$ is dependent in $T(G')^*$, then $U$ is dependent in $T(G)^*$. That is absurd.
    \\Subcase 2) If $V$ is independent in $T(G')^*$, then any matching of $V$ to $G'_N\setminus B$ in $D^J_B(T(G'))$ must be as follows: $v$ is matched to $n$ for some $v\in V$, and $V\setminus \{v\}$ is matched to $G_N\setminus (B\setminus \{n\})$. Then $V\cup \{n\}$ is dependent in $T(G')^*$. And any circuit of $T(G')^*$ contained in $V\cup \{n\}$ must contain $n$. 
\end{proof}

\section{Decomposability of lattice of flats}\label{sec:L(T(G))}
Recall that, from the previous chapter, we mentioned that the lattice of flats of a matroid $M$ is a geometric lattice and it is denoted by $\mathscr{L}(M)$. If $(\mathscr{L}_1,\leqslant_1),(\mathscr{L}_2,\leqslant_2)$ are lattices, their \textit{direct product} is the lattice $\mathscr{L}_1\times \mathscr{L}_2=\{a,b:a\in \mathscr{L}_1,b\in \mathscr{L}_2\}$ and the order is given by $(a,b)\leqslant (c,d)$ if $a\leqslant_1 c$ and $b\leqslant_2 d$. Note that direct product of geometric lattices is still a geometric lattice and a component in a direct product decomposition of a geometric lattice is also a geometric lattice (see for instance, \cite[p.60-61]{Aigner}).

\begin{figure}[h]
    \centering
    \includegraphics[width=1.1\linewidth]{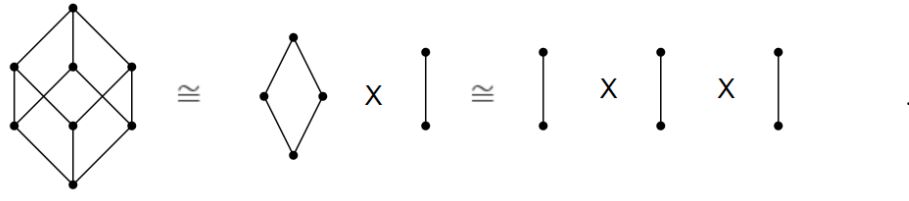}
    \caption{Decomposition of a geometric lattice.}
    \label{fig:decom. geo.}
\end{figure}

Figure \ref{fig:decom. geo.} shows a decomposition of a geometric lattice.
\noindent The following theorem characterizes the decomposability of a geometric lattice in terms of the complements of its atoms. (In a lattice $\mathscr{L}$ with $\mathscr{L}_O,\mathscr{L}_I$, $c\in \mathscr{L}$ is called a complement of $x\in \mathscr{L}$ if $c\wedge x=\mathscr{L}_O$ and $c\vee x=\mathscr{L}_I$).

\begin{theorem}[{\cite[Maeda]{Maeda}}]
    A geometric lattice $\mathscr{L}$ is directly indecomposable if and only if any two atoms share a common complement.
\end{theorem}

\begin{theorem}
    If $G$ is a bipartite graph such that $T(G)$ is simple, let $G'$ be a bipartite graph with the same bipartition as $G$, obtained by adding an edge to $G$, then
    \begin{enumerate}
        \item If $\mathscr{L}(T(G))$ is directly indecomposable and $\rho(T(G'))=\rho(T(G))$, then $\mathscr{L}(T(G'))$ is directly indecomposable.
        \item If $\mathscr{L}(T(G))$ is directly decomposable and $\rho(T(G'))=\rho(T(G))+1$, then $\mathscr{L}(T(G'))$ is directly decomposable.
    \end{enumerate}
\end{theorem}
\begin{proof}
    The height function for a geometric lattice coincides with the rank of the flats in the simple matroid corresponding to the lattice. Therefore, $\rho$ is used in place of $h$ for the height function. Since $\mathscr{L}(T(G))$ is semimodular, we have $\rho(x)+\rho(y)\geq \rho(x\vee y)+\rho(x\wedge y)$ for all $x,y\in \mathscr{L}(T(G))$. Thus, a complement of an atom must be a copoint in $\mathscr{L}(T(G))$. 

    \noindent For (1), 
    we claim that, for every hyperplane $H$ of $T(G)$, a subset of $H$ is a hyperplane of  $T(G')$. Let $\rho$ and $\rho'$ denote the rank function of $T(G)$ and $T(G')$, respectively. If $\rho'(H)=\rho (H)$, it means $H$ stays as a hyperplane in $T(G')$, we are done. If $\rho'(H)=\rho(H)+1$, just pick a subset $H'\subseteq H$ such that $\rho'(H')=\rho(H)$, then $H'$ is a hyperplane of $T(G')$. The desired result follows directly from our claim. 

    \noindent For 2) 
   Since $T(G)$ is directly decomposable, there exist atoms $p$ and $q$ in $\mathscr{L}(T(G))$ such that they do not have any common complement in $\mathscr{L}(T(G)).$ We get our desired conclusion if we can show that $p,q$ do not share any common complement in $\mathscr{L}(T(G'))$. We will assume that they share a common complement in $\mathscr{L}(T(G'))$ and derive a contradiction. This assumption is equivalent to saying that there is a hyperplane $H'$ such that $H'$ contains neither $p$ nor $q$. Comparing $\rho'(H')$ with $\rho(H')$, we may further split our investigation to two cases:

   \noindent Case i) $\rho'(H')=\rho(H')$
   \\ We have $\rho(H')=\rho'(H')=\rho'(T(G'))-1=\rho(T(G))$, and $H'\subsetneq G'_N=G_N$. If $n\notin H'$, we claim that $G_N\setminus H'=n$. Suppose, contrary to our claim, there is an $s\neq n$ such that $s\in G_N\setminus H'$. Due to $\rho'(s\cup H')=\rho'(T(G'))=\rho(T(G))+1$ and the only added edge to $G$ is adjacent to $n$, we have $\rho(s\cup H')=\rho'(s\cup H')=\rho (T(G))+1.$ This is impossible. So $n$ is the only element in $G_N$ with $n\notin H'$. But this contradicts our assumption that $H'$ contains neither $p$ nor $q$. Therefore, we have $n\in H'$. Let $\rho(T(G))=r$. $H'$ contains a base $\{h'_1,\ldots,h'_{r-1},n\}$ of $T(G).$ Since $H'\setminus n$ is of rank $r-1$ in $T(G)$, it is contained in some hyperplane $H$ of $T(G)$. We claim that $H=H'$. Suppose, contrary to our claim, there is a $s\in H\setminus H'$. Then we have $\rho(\{h'_1,\ldots,h'_{r-1},s\})=r-1$, and $\rho'(\{h'_1,\ldots,h'_{r-1},s\})=r$. This is not possible since the only added edge to $G$ is incident to $n$. Thus, we have to reach our claim $H=H'$. But this is again a contradiction because $H$ has to contain $p$ or $q$.

   \noindent Case ii) $\rho'(H')=\rho(H')+1$
   \\ In this case, $n$ must be in $H'$. We have 
      $$\rho(H')=\rho'(H')-1=\rho'(T(G'))-2=\rho(T(G))-1.$$
      So $H'$ is contained in some hyperplane $H$ of $T(G).$ We claim $H'=H$. Suppose, contrary to our claim, there is a $h\in H\setminus H'$. We have $\rho'(H'\cup h)=\rho'(H')+1$. Because every base in $T(G')$ must contain $n$ and any matching of the base must employ the newly added edge to match $n$, there must be an edge different from the newly added one to match $h$. That means
      $\rho (H'\cup h)=\rho(H')+1$, which is absurd since $H'$ is of the same rank as the hyperplane $H$ in $T(G)$ and $h\in H.$ Thus, we have to reach our claim $H'=H$. Again, this is still a contradiction.
    \end{proof}

\begin{figure}[h]
    \centering
    \includegraphics[width=0.68\linewidth]{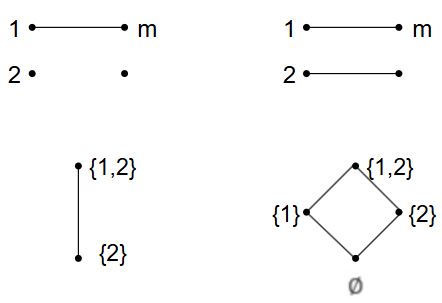}
    \caption{The preceding theorem can fail if $T(G)$ is not simple.}
    \label{fig:nonsimple indecom. to decom.}
\end{figure}
Figure \ref{fig:nonsimple indecom. to decom.} shows how the preceding theorem can fail without the assumption of simplicity of $T(G)$. $T(G)$ is a matroid on the set $\{1,2\}$ where $\{1\}$ is a coloop of $T(G)$, and $\{2\}$ is a loop of $T(G)$. $G'$ is obtained from $G$ by adding the edge $(2,m)$. But $\mathscr{L}(T(G))$ is indecomposable while $\mathscr{L}(T(G'))$ is decomposable.

\section{Dynamics of circuits}
In this section, we see how the number of circuits changes as we add an edge to a bipartite graph. 
\begin{figure}[h]
    \centering
    \includegraphics[width=0.5\linewidth]{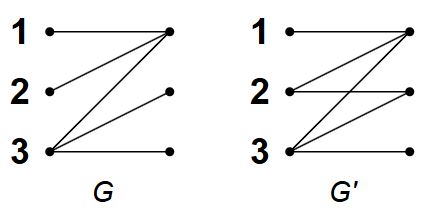}
    \caption{Number of circuits strictly decreases.}
    \label{fig:circuit decrease}
\end{figure}

\begin{figure}[h]
    \centering
    \includegraphics[width=0.45\linewidth]{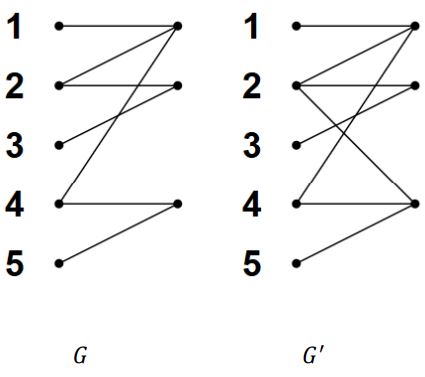}
    \caption{Number of circuits strictly increases.}
    \label{fig:circuit increase}
\end{figure}

\begin{example}
    In Figure \ref{fig:circuit decrease}, $T(G)$ has one circuit $\{1,2\}$. $G'$ is obtained by adding one edge to $G$, and $T(G')$ is free. In Figure \ref{fig:circuit increase}, $T(G)$ has $\{1,2,3\},\{1,2,4,5\},\{2,3,4,5\},\{1,3,4,5\}$ as circuits, while $T(G')$ has $\{1,2,3,4\},\{1,2,3,5\},\{1,2,4,5\},\{2,3,4,5\},\{1,3,4,5\}$ as circuits.
\end{example}

\begin{lemma}\label{lemma:circuit connected}
    Let $G_1,G_2$ be two connected components of a bipartite graph $G$, then $T(G)$ cannot have a circuit containing elements from both $G_{1N}$ and $G_{2N}$.
\end{lemma}
\begin{proof}
    Suppose $T(G)$ possesses a circuit $S\cup Q$ where $S,Q$ are nonempty subsets of $G_{1N},G_{2N}$, respectively. Since $S$ and $Q$ are independent in $T(G)$, this implies that there are $J_S\subseteq G_{1M},J_Q\subseteq G_{2M}$ matching to $S,Q$, respectively. But that means $J_S\cup J_Q\subseteq G_M$ is a matching for $S\cup Q$, a contradiction. 
\end{proof}

\begin{theorem}
    Let $G$ be a bipartite graph with at least two connected components, and $G'$ be obtained by adding an edge connecting two components of $G$ so that the bipartition of $G$ is preserved and $\rho(T(G'))=\rho(T(G))$. Then there is an injective function $f\colon \mathscr{C}(T(G))\rightarrow \mathscr{C}(T(G')).$ 
\end{theorem}
\begin{proof}
   Suppose the newly added edge connects $m\in G_{2M}$ and $n\in G_{1N}$, where $G_1,G_2$ are two components of $G$. If $n$ is a coloop of $T(G)$, we can let $f$ be the identity function on $\mathscr{C}(T(G)).$ Suppose $n$ is not a coloop of $T(G)$. Note that in any matching of an arbitrary base of $T(G)$, $m$ must be matched to some element from the base. Otherwise, we run into a contradiction to our assumption $\rho(T(G))=\rho(T(G'))$. Due to Lemma \ref{lemma: exists |G_M|=|base|}, we may assume $|G_{2M}|$ to be the cardinality of the maximal independent set in $G_{2N}$. Given $C_1\in \mathscr{C}(T(G))$ with $n\in C_1$, we have $|C_1\cup G_{2N}|\geq |(C_{1\mathscr{A}})_M\cup G_{2M}|$, where $\mathscr{A}$ is the presentation of $T(G)$ determined by the bipartite graph $G$. Due to $|(C_{1\mathscr{A}})_M\cup G_{2M}|\;=|(C_{1\mathscr{A'}})_M\cup G_{2M}|$, where $\mathscr{A'}$ is the presentation of $T(G')$ determined by $G'$, we have $C_1\cup G_{2N}$ is dependent in $T(G')$. Therefore, there is $X\subseteq G_{2N}$ such that $C_1\cup X\in \mathscr{C}(T(G')).$ Fix the set $X$. We claim that $C\cup X\in \mathscr{C}(T(G'))$ for all $C\in \mathscr{C}(T(G))$ with $n\in C$. This claim can be easily observed to be true by noticing that $m$ is the only element in $G_{2M}$ that is adjacent to an element, namely $n$, in any $C\in \mathscr{C}$ with $n\in C$. To be precise, we know that any subset of $C_1\cup X$ with $|C_1\cup X|-1$ can be matched in $G'$. Let that subset be $(C_1\cup X)\setminus z$ with $z\in X$. Then for an arbitrary $C\in \mathscr{C}(T(G))$, $(C\cup X)\setminus z$ can also be matched in $G'$ because $C_1$ (resp. $C$) can be matched by $(C_1)_M\cup \{m\}$ (resp. $(C)_M\cup \{m\}$) and $X\setminus z$ can be matched by some subset in $G_{2M}\setminus m$. On the other hand, if we take an element $s$ out of $C\cup X$ with $s\in C$, then observe that $C\setminus s$ can be matched by $(C\setminus s)_M$ and $X$ can be matched by some subset of $G_{2M}.$
   
   Thus, we may define $f$ to be 
   $$f(C) = \left\{
\begin{array}{ll}
    C\cup X & \text{if } C\in \mathscr{C}(T(G)), n\in C \\
    C  & \text{if } C\in \mathscr{C}(T(G)), n\notin C
\end{array}
\right.$$
\end{proof}

We also have a result on how the number of bases changes.

\begin{proposition}
    If $G'$ is a bipartite graph obtained by adding an edge to a bipartite graph $G$, and $G'$ has the same bipartition as $G$. If $T(G)$ is simple, and $\rho(T(G'))=\rho(T(G))+1$, then the number of bases of $T(G')$ is less than or equal to that of $T(G)$. 
\end{proposition}
\begin{proof}
     Suppose the newly added edge is incident to $n\in G_N$, and $B'$ is a basis of $T(G')$, then $B'\setminus \{n\}$ is a basis of $T(G)$. It is clear that if $B',B''$ are two distinct new bases created, then $B'\setminus \{n\},B''\setminus \{n\}$ are two distinct bases in $T(G)$ eliminated in the process from $T(G)$ to $T(G')$.
\end{proof}

\section{Some characterizations of circuits}
Let $\mathscr{A}$ be a presentation for a transversal matroid $M$, with $S\subseteq E(M)$. Then we let $S_{\mathscr{A}}$ denote the induced subgraph of $S\cup X$ in $G$, where $X$ is the set consisting of all the elements of $(G_{\mathscr{A}})_M$ which are adjacent to at least one element in $S$.  

\begin{remark}\label{re:(G)_M-(S)_M NO S}
   Let $(S_{\mathscr{A}})_M$ be the set $X$ we discussed above. We emphasize that since $X$ is defined to collect all elements of $(G_{\mathscr{A}})_M$, no element in $(G_{\mathscr{A}})_M\setminus (S_{\mathscr{A}})_M$ can be adjacent to any element in $S$. 
\end{remark}
The following proposition is a consequence of Hall's theorem (\cite{Hall},\cite[p.98]{Welsh}). We include a proof that is adapted to our content.

\begin{proposition}\label{prop: cir. N=M+1} 
    Let $\mathscr{A}$ be a presentation for a transversal matroid $M$. If $C\in \mathscr{C}(M)$, then $|(C_{\mathscr{A}})_N|=|(C_{\mathscr{A}})_M|+1$.
\end{proposition}
\begin{proof}
    Pick a $C'=\{c_1,\ldots,c_{|C|-1}\}\subseteq C$ of cardinality $|C|-1$. Then $C'$ must be independent in $M$. Thus, there exists a subset $J=\{m_1,\ldots,m_{|C|-1}\}\subseteq(G_{\mathscr{A}})_M$ such that $m_i$ is matched to $c_i$ in $G_{\mathscr{A}}, 1\leq i\leq |C|-1.$ $\{c_{|C|}\}=C\setminus C'$ cannot be isolated and it cannot be adjacent to vertices in $(G_{\mathscr{A}})_M\setminus J$ (otherwise $C$ becomes independent), so it must be adjacent to at least one of $J$. Our goal is to show that none of the vertices in $C$ can be adjacent to any vertex in $(G_{\mathscr{A}})_M\setminus J$. 
    
     Suppose $c_{|C|}$ is adjacent to $m_1$, then no vertex in $\{c_{|C|},c_1\}$ can be adjacent to any vertex in $(G_{\mathscr{A}})_M\setminus J$. Otherwise, if one of $\{c_{|C|},c_1\}$ is matched to an element $q\in (G_{\mathscr{A}})_M\setminus J$, then $\{c_{|C|},c_1\}$ can be matched by $\{q,m_1\}$. But $C\setminus \{c_{|C|},c_1\}$ can be matched by $J\setminus \{m_1\}$, so $C$ is independent. This is absurd. Now, since $\{c_{|C|},c_1\}$ is not a circuit of $M$, so at least one of them must be adjacent to some vertex in $J\setminus \{m_1\}$ so that $\{c_{|C|},c_1\}$ is independent. Now for $\{c_{|C|},c_1,\ldots,c_k\}\subsetneq C$,$1\leq k <|C|-1$, since at least one of them should be adjacent to at least one vertex in $J\setminus \{m_1,\ldots,m_k\}$, without loss of generality, we always assume it is $m_{k+1}$. Now we proceed with an inductive argument. Our induction hypothesis is to suppose that if $\{c_{|C|},c_1,\ldots,c_k\},1\leq k< |C|-1$ is a set such that if any one vertex in $(G_{\mathscr{A}})_M\setminus \{m_1,\ldots,m_k\}$ is adjacent to any one vertex in $\{c_{|C|},c_1,...,c_k\}$, then $\{c_{|C|},c_1,\ldots,c_k\}$ is independent in $M$. Then we claim that if any one vertex in $\{c_{|C|},c_1,\ldots,c_k,c_{k+1}\}$ is adjacent to any one vertex in $(G_{\mathscr{A}})_M\setminus \{m_1,\ldots,m_k,m_{k+1}\}$, $\{c_{|C|},c_1,\ldots,c_k,c_{k+1}\}$ is independent in $M$. 
    
    To verify the claim, if some $m_0\in (G_{\mathscr{A}})_M\setminus \{m_1,\ldots,m_{k+1}\}$ is matched to $c_{k+1}$, since $m_{k+1}$ is adjacent to at least one of $\{c_{|C|},c_1,\ldots,c_k\}$, by induction hypothesis the claim follows. Or if $m_0$ is adjacent to any one of $\{c_{|C|},c_1,\ldots,c_k\}$, again by induction hypothesis and the assumption that $m_{k+1}$ is matched to $c_{k+1}$, the claim still follows. And our desired result is a consequence of our claim. When $k+1$ reaches $|C|-1$, our claim shows that if any element in $C=\{c_{|C|},c_1,\ldots,c_{|C|-1}\}$ is adjacent to any element in $(G_{\mathscr{A}})_M\setminus \{J\}$, then $C$ becomes independent, which is a contradiction.
\end{proof}

\begin{corollary}\label{coro:(C_A)_M matched}
    If $\mathscr{A}$ is a presentation for a transversal matroid $M$, and $C\in \mathscr{C}(M)$, then every subset of $(G_{\mathscr{A}})_M$ can be matched in $C_{\mathscr{A}}$. 
\end{corollary}
\begin{proof}
    Every subset of $C$ of cardinality $|C|-1$ is independent and thus can be matched in $G_{\mathscr{A}}$. By Remark \ref{re:(G)_M-(S)_M NO S}, they can be matched in $C_{\mathscr{A}}$. By Proposition \ref{prop: cir. N=M+1}, $|(C_{\mathscr{A}})_M|\;=|C|-1$, so $|(C_{\mathscr{A}})_M|$ can be matched.  
\end{proof}

\begin{definition}
    A subgraph $H$ of a bipartite graph $G$ is called type Z if the induced subgraph of $H_V$ in $G$ is of the form $n_1m_1...n_km_k$ (expression $n_1m_1...n_km_k$ means that two vertices are adjacent in the subgraph if they are adjacent in $n_1m_1...n_km_k$) such that $m_k$ is not adjacent to $n_1$, and also $m_k$ is not adjacent to any vertex in $G_N\setminus \{n_1,\ldots,n_k\}$, where $\{n_1,
    \ldots,n_k\}= H_N,\{m_1,\ldots,m_k\}=H_M$.
\end{definition}

\begin{figure}[h]
\begin{subfigure}{0.53\textwidth}
\includegraphics[width=0.5\linewidth, height=6cm]{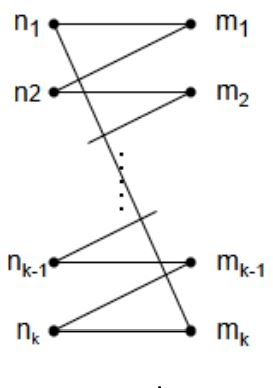} 
\caption{}
\end{subfigure}
\begin{subfigure}{0.53\textwidth}
\includegraphics[width=0.5\linewidth, height=6cm]{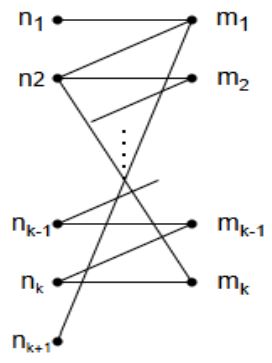}
\caption{}
\end{subfigure}
\caption{An illustration of the subgraph $H$.}
\label{fig:cir.curcuit}
\end{figure}

\begin{lemma}\label{lemma:cir.}
    Suppose $G$ is a bipartite graph and it contains a subgraph $H$ which can be read as $n_1m_1n_2m_2\ldots n_km_kn_1$ such that $H_N=\{n_1,\ldots,n_k\},H_M=\{m_1,\ldots,m_k\}$, if $H'$ is obtained by connecting a vertex $n_{k+1}\in G_N$ to any vertex of $\{m_1,\ldots,m_k\}$, then $\{n_1,\ldots,n_k,n_{k+1}\}$ is a circuit in $T(H')$.  
\end{lemma}
\begin{proof}
    Figure \ref{fig:cir.curcuit} (a) is a visualization of such a $H$. Now suppose we form $H'$ by connecting a vertex $n_{k+1}$ to a $m_j,1\leq j\leq k$. Note that $T(H')$ is a circuit if any subset of cardinality $k$ in $\{n_1,\ldots,n_{k+1}\}$ is independent. Clearly, $\{n_1,\ldots,n_k\}$ is independent. Now, we want to show that any cardinality-$k$ subset $P\subseteq \{n_1,\ldots,n_k,n_{k+1}\}$ such that $n_{k+1}\in P$ can be matched by $\{m_1,\ldots,m_k\}$. In all three cases below, $n_{k+1}$ is matched by $m_j$.
    \\ Case 1. $P=\{n_1,\ldots,n_k,n_{k+1}\}\setminus \{n_i\}$ for some $1\leq i<j$
    \\ $n_h$ is matched by $m_h$ for $1\leq h\leq i-1$ or $j<h\leq k$;
    \\ $n_h$ is matched by $m_{h-1}$ for $i+1\leq h\leq j.$
    \\Case 2. $P=\{n_1,\ldots,n_k,n_{k+1}\}\setminus \{n_i\}$ for some $j<i\leq k$
    \\ $n_h$ is matched by $m_{h-1}$ for $2\leq h\leq j$ or $i+1\leq h\leq k$;
    \\ $n_1$ is matched by $m_k$;
    \\ $n_h$ is matched by $m_{h}$ for $j<h<i$.
    \\Case 3. $P=\{n_1,\ldots,n_k,n_{k+1}\}\setminus\{n_i\}, i=j$
    \\ $n_h$ is matched by $m_h$ for $1\leq h\leq k$ such that $h\neq j$.

    Note that adding more edges inside $H'$ does not affect the result. Any set $P$ of cardinality $k$ in $\{n_1,\ldots,n_{k+1}\}$ can be matched, and it still can be matched if we add any edge whose endpoints are vertices in $H'_V$. However, $\{n_1,\ldots,n_k,n_{k+1}\}$ still cannot be matched by $\{m_1,\ldots,m_k\}$. 
\end{proof}

Figure \ref{fig:cir.curcuit} (b) suggests that it is necessary for $n_1$ to be adjacent to $m_k$ in the previous lemma, as $\{n_1,n_{k+1}\}$ becomes a circuit in $T(H')$.

To make our proof of the following theorem smoother, we set some more terminology: A \textit{Z-chain} in a bipartite graph $G$ is a subgraph of the chain form $a_1b_1\ldots a_ib_i$, where $\{a_1,\ldots,a_i\}\subseteq G_N$, and $\{b_1,\ldots,b_i\}\subseteq G_M$, such that two vertices are adjacent in the subgraph if (not necessarily only if) they are adjacent in the chain. 

If we have a Z-chain $a_1b_1\ldots a_ib_i$ in a bipartite graph $G$, and there are $a'\in G_N$, $b'\in G_M$ such that $a'\notin \{a_1,\ldots,a_i\}$, $b'\notin \{b_1,\ldots,b_i\}$, so that $a_1b_1\ldots a_ib_ia'b'$ is a Z-chain in $G$, we say the Z-chain $a_1b_1\ldots a_ib_i$ is extended to $b'$ through $a'$. We also need the following straightforward lemma:
\begin{lemma}\label{lemma: extended Z chain}
    Let $G$ be a bipartite graph without type-Z subgraph, and $G$ has a Z-chain $n_1m_1\ldots n_sm_s$ such that $m_s$ is not adjacent to $n_1$. Moreover, there are $\{n_{t_1},\ldots,n_{t_r}\}\subseteq G_N$, $\{m_{k_1},\ldots,m_{k_r}\}\subseteq G_M$ such that $n_{t_i}$ is adjacent to $m_{t_i}$ for all $1\leq i\leq r$, with $\{n_{t_1},\ldots,n_{t_r}\}\cap \{n_1,\ldots,n_s\}=\emptyset$, $\{m_{k_1},\ldots,m_{k_r}\}\cap \{m_1,\ldots,m_s\}=\emptyset$, then at least one of the following occurs:
    \begin{itemize}
        \item $m_s$ is adjacent to some $n\notin \{n_1,\ldots,n_s,n_{t_1},\ldots,n_{t_r}\}.$
        \item There is an $m_{k_i}\in \{m_{k_1},\ldots,m_{k_r}\}$ such that $m_{k_i}$ is an element to which the Z-chain $n_1m_1\ldots n_sm_s$ can be successively extended through an element in $\{n_{t_1},\ldots,n_{t_r}\}$, and $m_{k_i}$ is adjacent to $n_1$, or to $n\notin \{n_1,\ldots,n_s,n_{t_1},\ldots,n_{t_r}\}.$
    \end{itemize}
\end{lemma}
\begin{proof}
    If $m_s$ is not adjacent to any element in $\{n_{t_1},\ldots,n_{t_r}\}$, due to our assumption that $G$ contains no type Z subgraph, $m_s$ must be adjacent to $n_1$, or to some $n\notin \{n_1,\ldots,n_s,n_{t_1},\ldots,n_{t_r}\}$. If $m_s$ is adjacent to some element in $\{n_{t_1},\ldots,n_{t_r}\}$, without loss of generality, suppose $m_s$ is adjacent to $n_{t_1}$. Then the Z chain $n_1m_1\ldots n_sm_s$ is extended to $m_{k_1}$ through $n_{t_1}$. If $m_{k_1}$ is adjacent to $n_{t_2}$, then the extended Z-chain can be further extended to $m_{k_2}$ through $n_{t_2}$. Continuing to extend our Z chain through elements in $\{n_{t_1},\ldots,n_{t_r}\}$, there is an $m_{k_i}\in \{m_{k_1},\ldots,m_{k_r}\}$ such that our Z- chain can be successively extended as $n_1m_1\ldots n_sm_sn_{t_1}m_{k_1}\ldots n_{t_i}m_{k_i}$ and it cannot be further extended to an element in $\{m_{k_1},\ldots,m_{k_r}\}$ through an element in $\{n_{t_1},\ldots,n_{t_r}\}$. Then to avoid a type Z subgraph, $m_{k_i}$ must be adjacent to $n_1$, or to some $n\notin \{n_1,\ldots,n_s,n_{t_1},\ldots,n_{t_r}\}.$ 
\end{proof}

\begin{theorem}\label{thm:type Z Thm}
    If $G$ is a bipartite graph, $|G_N|=|G_M|\;>2$, then $G_N$ is a circuit in $T(G)$ if $|G_N|\;=|G^c_M|+1$, $G$ is connected except $|G_M\setminus G^c_M|$ isolated vertices and $G$ contains no subgraph of type Z.
\end{theorem}
\begin{proof}
    Suppose $G$ satisfies the assumptions of the theorem, but $G_N$ is not a circuit of $T(G)$. Our goal is to show that $G$ contains a type Z subgraph. First, note that $G$ is a connected graph except $|G_M\setminus G_M^c|\;=1$ isolated vertex. 

    Due to $|G_N|\;=|G^c_M|+1$, $G_N$ is dependent in $T(G)$. Thus, $G_N$ properly contains some $C\in \mathscr{C}(T(G))$. Let $\mathscr{A}$ be the presentation for $T(G)$ corresponding to $G$. Pick a non-isolated $m_s\in G_M\setminus (C_{\mathscr{A}})_M$. In view of Remark \ref{re:(G)_M-(S)_M NO S}, elements in $G_M\setminus (C_{\mathscr{A}})_M$ cannot be adjacent to any element in $C$. Because of the connectedness of $G$, there must be a path from $m_s$ to some element $m_c\in (C_{\mathscr{A}})_M$, and $m_c$ has to be adjacent to some $c\in C$. Without loss of generality, we may assume that the path is $m_sn_s\ldots m_1n_1m_c$. So we have a Z-chain $cm_cn_1m_1\ldots n_sm_s$. If $m_s$ is not adjacent to $G_N\setminus (C\cup \{n_1,\ldots,n_s\})$, we get a type Z subgraph. Otherwise, $m_s$ has to be adjacent to some $n'\notin C\cup\{n_1,\ldots,n_s\}$. See Figure \ref{fig:type Z Thm illu} for an illustration. The part colored in blue corresponds to $C_{\mathscr{A}}$.

    Note that $$|C\cup \{n_1,\ldots,n_s,n'\}|\;=|(C_{\mathscr{A}})_M\cup \{m_1,\ldots,m_s\}|+2.$$
    Together with the assumption $|G_N|\;=|G_M^c|+1$, we know that there must be a non-isolated element $m'\in G_M\setminus ((C_{\mathscr{A}})_M\cup \{m_1,\ldots,m_s\}).$ Then we ask the following question:
    \begin{equation}\label{eq:type Z question}
        \mbox{Is that non-isolated element adjacent to some element in } \{n_1,\ldots,n_s,n'\}?
    \end{equation}

    \begin{figure}[h]
    \centering
         \includegraphics[width=0.22\linewidth,height=9cm]{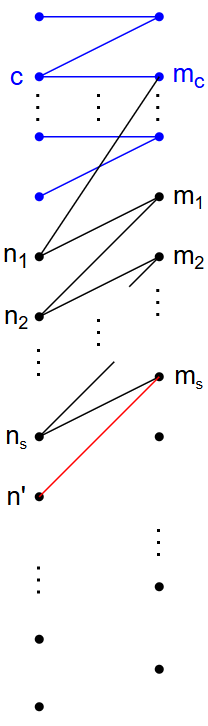}
         \caption{$m_s$ is adjacent to some $n'\in G_N\setminus( C\cup \{n_1,
         \ldots,n_s\})$.}
        \label{fig:type Z Thm illu}
    \end{figure}

\noindent Case 1) The answer to \eqref{eq:type Z question} is yes. 

The non-isolated element $m'$ is adjacent to some $n_i\in \{n_1,\ldots,n_s,n'\}$. See Figure \ref{fig:type Z Thm illu 1} for an illustration. Then we get a new Z-chain $cm_c\ldots n_im'$, and it either results in a type Z subgraph and we are done, or by Lemma \ref{lemma: extended Z chain}, the non-isolated $m'$ or an element in $\{m_1,\ldots,m_s\}$ has to be adjacent to an element in $G_N\setminus (C\cup \{n_1,\ldots,n_s,n'\}).$

\vspace{5mm}
\noindent Case 2) The answer to \eqref{eq:type Z question} is no. 

Let $m'$ be $m_{s+w}$. We employ the same technique in the first paragraph of our proof. That is, there must be a path from $m_{s+w}$ to some element in $(C_{\mathscr{A}})_M\cup \{m_1,\ldots,m_s\}$. Let's say there is a path from $m_{s+w}$ to $m_1$ (we choose $m_1$ because it is convenient to illustrate this with a picture, and this is not going to affect the argument). Let the path be $m_{s+w}n_{s+w}\ldots m_{s+1}n_{s+1}m_1$. Then we have a Z chain $cm_cn_1m_1n_{s+1}m_{s+1}\ldots n_{s+w}m_{s+w}.$ See Figure \ref{fig:type Z Thm illu. case 2} for an illustration. Then either we have a type Z subgraph and we are finished, or $m_{s+w}$ is adjacent to some element in $G_N\setminus (C\cup \{n_1,\ldots,n_s,n',n_{s+1},\ldots,n_{s+w}\}).$

\begin{figure}[h]
    \begin{subfigure}{0.53\textwidth}
        \includegraphics[width=0.44\linewidth,height=10.6cm]{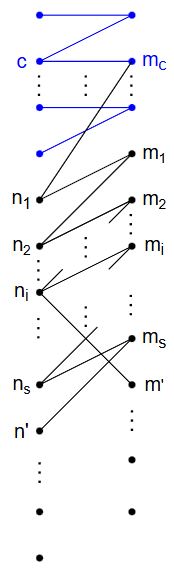}
        \caption{}
        \label{fig:type Z Thm illu 1}
    \end{subfigure}
    \begin{subfigure}{0.53\textwidth}
        \includegraphics[width=0.49\linewidth,height=10.6cm]{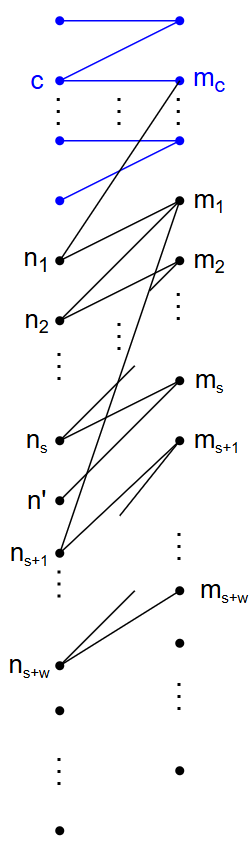}
        \caption{}
        \label{fig:type Z Thm illu. case 2}
    \end{subfigure}
       
        \caption{A non-revealed, non-isolated vertex in $G_M$ is adjacent to either some element in $\{n_1,\ldots,n_s,n'\}$ or not.}
    \label{fig:placeholder}
\end{figure}

Along the course of our proof, the elements from $G_N$ or $G_M$ that are already known to be incident to some known edge are called revealed. For instance, in Figure \ref{fig:type Z Thm illu. case 2}, $C\cup \{n_1,\ldots,n_s,n',n_{s+1},\ldots,n_{s+w}\}$ and $\{(C_{\mathscr{A}})_M\cup \{m_1,\ldots,m_s,m_{s+1},\ldots,m_{s+w}\}\}$ are revealed.

In summary, in order to avoid a type Z subgraph, in case (1), one more non-isolated element in $G_M$ is revealed in $G$,  and one more element in $G_N$ is revealed. In case (2), the number of newly revealed elements in $G_N$ is strictly one more than that of $G_M$. Then again, due to the assumption $|G_N|\;=|G_M^c|+1$, in either case there exists an  element in $G_M$ that should be non-isolated but not yet revealed; we then continue to ask \eqref{eq:type Z question} (with $\{n_1,\ldots,n_s,n'\}$ replaced by the subset of $G_N$ that is already revealed) on that element. Each time we answer that question, the cardinality of elements from $G_N$ being revealed is at least as big as that of elements from $G_M$ being revealed. Since the cardinality of elements that have been revealed strictly increases each time we answer and $|G_N|$ is finite, eventually, we will be running out of elements from $G_N$ to avoid a type Z subgraph.
\end{proof}

Now let us go through the background of connectivity of matroids. There are analogues in the cycle matroid $M(G)$ of the effects of a graph $G$ when we perform deletion or  contraction or taking dual on $G$ (see Figure \ref{fig:deletion, contraction}). It is natural to ask for a matroid concept corresponding to the connectedness of graphs. But every cycle matroid is isomorphic to the cycle matroid of a connected graph (see \cite[Prop.1.2.9]{Oxley} for a proof). So the ordinary notion of connectedness in graphs is not right in the context of matroids. As a result, following \cite[Whitney]{Whitney}, the concept of 2-connectedness is extended from graphs to matroids. 

A subset of the vertex set of a graph $G$ is called a \textit{vertex cut} if the new graph obtained from deleting the subset and its incident edges has strictly more connected components than $G$. For a connected graph $G$ having at least one pair of distinct non-adjacent vertices, the \textit{connectivity} $\kappa(G)$ of $G$ is the smallest integer $m$ such that $G$ has a vertex cut of cardinality $m$. If $G$ is a complete graph, then we take $\kappa(G)=|G_V|-1$. The following three propositions are standard in the field.

\begin{proposition}[{\cite[p.119]{Oxley}}]\label{prop:ref.connected 1}
    Let $G$ be a graph with three or more vertices without loops or isolated vertices. Then $G$ is 2-connected if and only if for every pair of distinct edges of $G$, there is a cycle containing both.
\end{proposition}
Define a relation $\xi$ on $E(M)$ of a matroid $M$ by $a$ $\xi$ $b$ if either $a=b$ or there is a circuit of $M$ containing $a$ and $b$, $a,b\in E(M)$.

\begin{proposition}[{\cite[p.119]{Oxley}}]\label{prop:ref.conected 2}
    The relation $\xi$ is an equivalence relation on $E(M)$ for every matroid $M.$
\end{proposition}
 The $\xi$ equivalence classes are called the connected components of the matroid $M$. A matroid $M$ is called \textit{2-connected} and often simply called \textit{connected} if $E(M)$ is a $\xi$-component. 

Let $M$ be a matroid and $X\subseteq E(M)$, by Theorem \ref{thm:rank function}, we have $\rho(X)+\rho(E(M)\setminus X)\geq \rho(M)$. The following proposition gives a sufficient and necessary condition for connectedness of matroids in terms of their rank functions.

\begin{proposition}[{\cite[p.123]{Oxley}}]\label{prop:ref.connected 3}
    A matroid $M$ is connected if and only if for any nonempty proper subset $T$ of $E(M)$, $$\rho(T)+\rho(E(M)\setminus T)>\rho(M).$$
\end{proposition}
Before we give some features of those bipartite graphs whose corresponding transversal matroid is connected, we want to give some features of those bipartite graphs corresponding to a circuit first. 

Recall that we let $G_V$ denote the set of vertices for a graph $G$. If $H_1,H_2$ are subgraphs of $G$, let $H_1\cup H_2$ be the subgraph of $G$ such that $(H_1\cup H_2)_V=(H_1)_V\cup (H_2)_V$ and $(H_1\cup H_2)_E$ is $(H_1)_E\cup (H_2)_E$. $H_1\cap H_2$ is the subgraph of $G$ defined analogously, in the obvious way.   

Given a graph $G=(V,E)$, let $G^c=(V^c,E^c)$ denote the graph whose vertex set $V^c$ is the subset of non-isolated vertices in $V$, and edge set $E^c=E$.

\begin{lemma}(Circuits Pasting)\label{lemma:cir. past.}
    If $G$ is a bipartite graph such that $T(G)$ is connected with $|G_N|\geq 2$ and $\mathscr{A}$ is the corresponding presentation for $T(G)$, then there exist distinct $C_1,...,C_n\in \mathscr{C}(T(G))$ such that $G^c=(C_1)_{\mathscr{A}}\cup ...\cup (C_n)_{\mathscr{A}}$, with $(C_m)_{\mathscr{A}}\cap ((C_{m-1})_{\mathscr{A}}\cup \ldots\cup (C_1)_{\mathscr{A}})\neq \emptyset$ for all $2\leq m\leq n$.
\end{lemma}
\begin{proof}
    Pick two distinct $x,y\in G_N$. Since $T(G)$ is connected, there exists $C_1\in \mathscr{C}(T(G))$ such that $\{x,y\}\subseteq C_1$. If $C_1=E(T(G))$, we are done. If not, pick $z\in E(T(G))\setminus C_1$. By connectedness of $T(G)$ there must be $C_2\in \mathscr{C}(T(G))$ such that $\{y,z\}\subseteq C_2$. Again, if $C_1\cup C_2=E(T(G))$, we are finished. If not, we repeat the above step. Since $|C_1\cup\ldots\cup C_n|$ strictly increases after $n$ steps for every positive integer $n$, and  $E(T(G))$ is finite, we have to stop after finitely many steps. 
\end{proof}

\begin{remark}
    If we have $C_1,\ldots,C_n\in \mathscr{C}(M)$ for a matroid $M$ such that $C_1\cup \ldots\cup C_n=E(M)$, we say $C_1,\ldots,C_n$ cover the matroid $M$. Note that this does not mean that $C_1,\ldots,C_n$ are all the circuits of $M$. 
\end{remark}

Recall that for any bipartite graph $G$, let $G^c_N$ and $G^c_M$ denote the non-isolated vertices of $G_N$ and $G_M$, respectively.
\begin{proposition}
    If $G$ is a bipartite graph such that $|G_N|\;\geq 2$ and $T(G)$ is connected, then $|G_N|\;>|G_M^c|$.
\end{proposition}
\begin{proof}
    Since $|G_N|\;\geq 2$, we may pick distinct $x,y\in G_N$. And due to the connectedness of $T(G)$, there is a circuit $C_1\in \mathscr{C}(T(G))$ containing $x$ and $y$. If $C_1=G_N$, by Proposition \ref{prop: cir. N=M+1}, we are finished. If not, by the circuit pasting, Lemma \ref{lemma:cir. past.}, there exist $C_1,...,C_n\in \mathscr{C}(T(G))$ such that $G^c=(C_1)_{\mathscr{A}}\cup ...\cup (C_n)_{\mathscr{A}}$ and $(C_m)_{\mathscr{A}}\cap ((C_{m-1})_{\mathscr{A}}\cup ...\cup (C_1)_{\mathscr{A}})\neq \emptyset$ for all $2\leq m\leq n$. Now we set our induction hypothesis:
    $$|C_1\cup \ldots\cup C_k|\;>|((C_1\cup\ldots\cup C_k)_{\mathscr{A}})_M|,$$
    for some integer $k$, $1\leq k<n$. Since $C_{k+1}$ is a circuit of $T(G)$, $C_{k+1}\cap (C_1\cup\ldots\cup C_k)$, as a proper subset of $C_{k+1}$, is independent in $T(G)$, and thus can be matched by some $J\subseteq ((C_{k+1})_{\mathscr{A}})_M$. By Remark \ref{re:(G)_M-(S)_M NO S}, $J\subseteq ((C_{k+1})_{\mathscr{A}})_M \cap ((C_1\cup\ldots\cup C_k)_{\mathscr{A}})_M)$. Therefore, we have 
    $$|((C_1\cup \ldots\cup C_k)_{\mathscr{A}})_M \cap ((C_{k+1})_{\mathscr{A}})_M|\;\geq |C_{k+1}\cap (C_1\cup \ldots\cup C_k)|.$$
    Combining this with Proposition \ref{prop: cir. N=M+1}, we obtain
    $$
    |C_{k+1}\setminus C_{k+1}\cap (C_1\cup\ldots\cup C_k)|\;>|((C_{k+1})_{\mathscr{A}})_M\setminus (((C_{k+1})_{\mathscr{A}})_M \cap ((C_1\cup \ldots\cup C_k)_{\mathscr{A}})_M)|.
    $$
    And clearly, we can partition $|C_1\cup\ldots\cup C_{k+1}|$ and $|((C_1\cup\ldots\cup C_k)_{\mathscr{A}})_M \cup ((C_{k+1})_{\mathscr{A}})_M|$ as follows:

    $$|C_1\cup\ldots\cup C_{k+1}|\;=|C_1\cup\ldots\cup C_k|+|C_{k+1}\setminus (C_{k+1}\cap (C_1\cup\ldots\cup C_k))|,$$

    and 
    
    \begin{flalign*}    
    &|((C_1\cup\ldots\cup C_k)_{\mathscr{A}})_M\cup ((C_{k+1})_{\mathscr{A}})_M|\\
    &\leq |((C_1\cup\ldots\cup C_k)_{\mathscr{A}})_M|+|((C_{k+1})_{\mathscr{A}})_M\setminus (((C_{k+1})_{\mathscr{A}})_M \cap ((C_1\cup \ldots\cup C_k)_{\mathscr{A}})_M)|.
    \end{flalign*}
    By our induction hypothesis, we have 
    $$|C_1\cup\ldots\cup C_{k+1}|\;>|((C_1\cup\ldots\cup C_{k+1})_{\mathscr{A}})_M|.$$ 
    \end{proof}

\begin{proposition}
    If $T(G)$ is a connected matroid such that $|G_N|\;=0$ or $|G_N|\;\ge 2$, then $\rho(T(G))=|G_M^c|$.
\end{proposition}
\begin{proof}
    A subset of $G_N$ is independent in $T(G)$ if and only if it can be matched in $G$. So the rank of $T(G)$, which is defined to be the cardinality of a maximal independent set, is the cardinality of a maximal subset of $G_N$ that can be matched in $G$. Certainly, $\rho(T(G))\leq |G_M^c|$. And the equality is achieved if we  show that $G_M^c$ can be matched in $G$. 
    
    If $|G_N|\;=0$, then $\rho(T(G))=0=|G_M^c|$. Suppose $|G_N|\;\geq 2$, by circuit pasting, Lemma \ref{lemma:cir. past.}, there exist distinct circuits $C_1,\ldots,C_n$ in $\mathscr{C}(T(G))$ such that $G^c=(C_1)_{\mathscr{A}}\cup\ldots\cup (C_n)_{\mathscr{A}}$ with $(C_m)_{\mathscr{A}}\cap ((C_{m-1})_{\mathscr{A}}\cup\ldots\cup (C_1)_{\mathscr{A}})\neq \emptyset$ for all $2\leq m\leq n$. If $n=1$, every subset of cardinality $|C_1|-1$ is independent. So by Proposition \ref{prop: cir. N=M+1}, we have $\rho(C_1)=|C_1|-1=|((C_1)_{\mathscr{A}})_M|$. Now we set our induction hypothesis:
    $$\rho(C_1\cup\ldots\cup C_k)=|((C_1\cup\ldots\cup C_k)_{\mathscr{A}})_M|,$$
    for an integer $k$ with $1\leq k<n$. Now, we proceed with our inductive step:

    \noindent Case 1) $((C_{k+1})_{\mathscr{A}})_M\subseteq ((C_1\cup\ldots\cup C_k)_{\mathscr{A}})_M$
   
    Our induction hypothesis is equivalent to saying that $((C_1\cup...\cup C_k\cup C_{k+1})_{\mathscr{A}})_M$ has a matching in $(C_1\cup\ldots\cup C_k)_{\mathscr{A}}$, thus in $(C_1\cup...\cup C_k\cup C_{k+1})_{\mathscr{A}}$. 

    \noindent Case 2) $((C_{k+1})_{\mathscr{A}})_M\nsubseteq ((C_1\cup\ldots\cup C_k)_{\mathscr{A}})_M$

    By Corollary \ref{coro:(C_A)_M matched}, $((C_{k+1})_{\mathscr{A}})_M\setminus ((C_1\cup\ldots\cup C_k)_{\mathscr{A}})_M$ can be matched to a subset of $C_{k+1}$. But by Remark \ref{re:(G)_M-(S)_M NO S}, a set being matched to $((C_{k+1})_{\mathscr{A}})_M\setminus (C_1\cup\ldots\cup C_k)_{\mathscr{A}})_M$ must be a subset of $C_{k+1}\setminus(C_1\cup\ldots\cup C_k)$. Together with our induction hypothesis, $((C_{k+1})_{\mathscr{A}})_M\cup ((C_1\cup\ldots\cup C_k)_{\mathscr{A}})_M$ can be matched. 
   \end{proof}

\begin{proposition}
    Let $G$ be a bipartite graph such that $T(G)$ is connected. If $G'$ is obtained by adding an edge to $G$ with endpoints in $G_N,G_M^c$, then $T(G')$ is also connected.  
\end{proposition}
\begin{proof}
    From the preceding proposition, $\rho(T(G))=|G_M^c|$. And $T(G)$ is connected if and only if $$\rho(H)+\rho(E(T(G))\setminus H)>\rho(M)$$ where $H$ is any nonempty proper subset of $E(M)$. In our context, $M=T(G)$, and adding edges as described does not change the value of $\rho(T(G))$. And any subset of $G_N$ that can be matched in $G$ can still be matched in $G'$, therefore, $\rho(H)$ and $\rho(E(T(G))\setminus H)$ do not strictly decrease. So the inequality still holds. 
\end{proof}

\begin{figure}[h]
    \centering
    \includegraphics[width=0.45\linewidth]{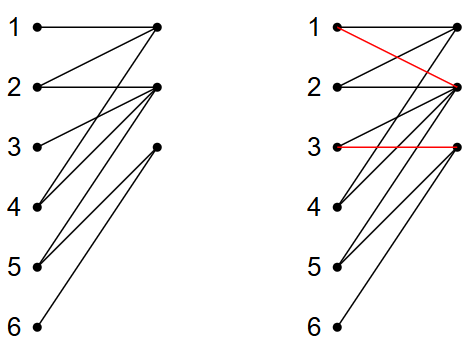}
    \caption{Adding edges between non-isolated vertices in $G$ for a connected $T(G)$ does not alter the connectedness.}
    \label{fig:add.e.not.al.conn.}
\end{figure}

In Figure \ref{fig:add.e.not.al.conn.}, we have an intuitive example of the situation described by the preceding proposition. Edges colored in red represent newly added edges to an existing bipartite graph. That does not affect the connectedness of the corresponding transversal matroids.

\chapter{Acyclic Presentation}
As there can be many presentations for a transversal matroid, it is natural to try to find a simple one. In this chapter, we characterize the bipartite graph $G_{\mathscr{A}}$ when there is a presentation $\mathscr{A}$ for a transversal matroid $M$ such that $G_{\mathscr{A}}$ is a forest. We call such a presentation an acyclic presentation of $M$. We are going to characterize the transversal matroids which admit acyclic presentations. We will also find conditions under which part of a presentation of a transversal matroid can be made acyclic.

Recall that if $\mathscr{A}$ is a presentation for a transversal matroid $M$, and $S\subseteq E(M)$ is an arbitrary subset of the ground set of $M$, we let $S_{\mathscr{A}}$ denote the induced subgraph of $S\cup K$ in $G_{\mathscr{A}}$, where $K$ is the set consisting of every vertex in $(G_{\mathscr{A}})_M$ such that it is adjacent to at least one vertex in $S$. 

\section{Matroids induced by general plabic graphs}
The class of transversal matroids that admit acyclic presentations is already understood, but this is not known to be clearly documented in the literature. The author is grateful for the suggestion from Alex Fink  for the approach in this section. The content in this section was developed through e-mail exchanges between us \cite{Fink}.

The concept of a plabic graph was introduced by Alexander Postnikov in 2006 \cite{Postnikov}. A \textit{planar bicolored graph}, or simply a \textit{plabic graph} is a planar undirected graph $G$ properly embedded into a closed disk. The vertices of $G$ that are on the boundary of the disk are called \textit{boundary vertices}, and the remaining vertices are called \textit{internal vertices}. Every internal vertex of $G$ is colored  black or white. Moreover, each boundary vertex is adjacent to exactly one internal vertex. In this section, we are going to deal with an alternate version that is more general in the sense that we do not require it to be planar. So in our approach, the boundary vertices in Postnikov's sense may not be placed on a disk and it is possible for some of them to be colored black or white.

\begin{definition}
    A generalized plabic graph $P$ is a finite simple graph such that each vertex of $P$ is colored in one of the following three colors: green, black or white. Moreover, a vertex colored green must be a leaf in $P$.
\end{definition}

\begin{remark} 
From now on, we mean a generalized plabic graph whenever we mention plabic graph unless otherwise stated.
\end{remark}

To derive a matroid from a plabic graph $P$, we need to have some orientations of the edges of $P$. Intuitively, an orientation of $P$ is just an assignment of directions to all edges of $P$. To be precise, let $V(P)$ and $E(P)$ denote the set of vertices and the set of edges of $P$, respectively. Let $E_{a,b}=\{(a,b),(b,a)\}$ for $a,b\in V(P), a\neq b.$ An \textit{orientation} $\mathscr{O}$ of $P$ is a function:
$$\mathscr{O} \colon E(P) \to \bigcup_{\substack{a,b\in V(P), \\a\neq b}} E_{a,b}$$
such that 
$$\mathscr{O}\vert_{\{a,b\}\in E(P)} \colon \{a,b\} \to E_{a,b} .$$ 
If $\mathscr{O}(\{a,b\})=(a,b)$, we say that the edge $(a,b)$ is directed out from $a$ and directed towards $b$ in the orientation $\mathscr{O}.$ 

\begin{definition}
An orientation $\mathscr{O}$ of a plabic graph $P$ is called a valid orientation if 
\begin{itemize}
    \item For each black vertex $b$ in $P$, there is exactly one edge directed out from $b$ in $\mathscr{O}$,
    \item For each white vertex $w$ in $P$, there is exactly one edge directed towards $w$ in $\mathscr{O}.$
    \end{itemize}
\end{definition}

The analogues in Postnikov's plabic graph are called \textit{perfect orientations} \cite{Postnikov}. He derives matroids out of perfect orientations of plabic graphs in his sense \cite[p.40]{Postnikov}. We can also derive matroids from plabic graphs out of valid orientations in an analogous way.  

Suppose $\mathscr{O}$ is an orientation of a plabic graph $P$. Let $\mathscr{O}[P]$ denote the set of green leaves of $P$ such that for $l \in \mathscr{O}[P]$, the only edge incident to $l$ is directed out from $l$ in $\mathscr{O}.$

\begin{theorem}
    Let $P$ be a plabic graph. Then $\{\mathscr{O}[P]: \mathscr{O} \mbox{ is a valid orientation of } P\}$ forms a collection of bases of a matroid on the set of green leaves of $P$. This matroid is denoted by $\mathscr{M}(P)$, and is called the matroid induced by the plabic graph $P$. If there is no valid orientation for $P$, we declare the collection to consist of the empty set.
\end{theorem}
\begin{proof}
    Let $X=\mathscr{O}_1[P], Y=\mathscr{O}_2[P]$, where $\mathscr{O}_1$ and $\mathscr{O}_2$ are two distinct valid orientations of $P$. Let $D$ be the subgraph of $P$ such that every edge in $D$ is directed differently in $\mathscr{O}_1$ and $\mathscr{O}_2$. Note that every vertex except a leaf has degree 2 in $D$. Let $x\in X$ such that $x\notin Y$, then $x\in D$. Furthermore, if $D_x$ is the component of $D$ containing $x$, then there must be another leaf $y$ in $D_x$ because the sum of all degrees of vertices in $D_x$ must be even. That means $D_x$ is actually a path graph. Under $\mathscr{O}_1$, $D_x$ is a directed path from $x$ to $y$. Now, we reverse the directions of all edges in $D_x$ and let other edges not in $D_x$ inherit the orientation of $\mathscr{O}_1$. The resulting orientation is still valid, and it corresponds to the base $(X\setminus \{x\})\cup \{y\}.$ This means the bases exchange condition specified in Theorem \ref{BE} holds.   
\end{proof}

\begin{example}
In Figure \ref{fig:$M(P)=U_2,3$}, the plabic graph $P$ we have is a star with three green leaves such that the center is colored black. And all three possible valid orientations of $P$ are displayed. The matroid induced by $P$ is the uniform matroid $U_{2,3}.$
    \begin{figure}[h]
        \centering
        \includegraphics[width=0.5\linewidth]{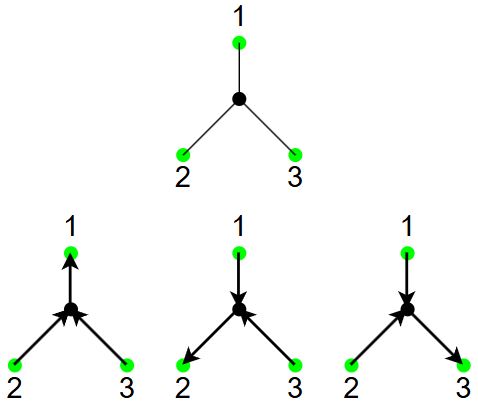}
        \caption{A plabic graph $P$ such that $\mathscr{M}(P)=U_{2,3}.$}
        \label{fig:$M(P)=U_2,3$}
    \end{figure}
\end{example}

In Figure \ref{fig:S(U_k,k-1),S(U_1,k)}, we have a plabic star with $k$ green leaves such that its center is colored black, and a plabic star with $k$ green leaves whose center is colored white. We denote them as $SU_{k-1,k}$ and $SU_{1,k}$, respectively. It is easy to see that $\mathscr{M}(SU_{k-1,k})=U_{k-1,k}$ and $\mathscr{M}(SU_{1,k})=U_{1,k}.$

\begin{figure}[h]
    \centering
    \includegraphics[width=0.4\linewidth]{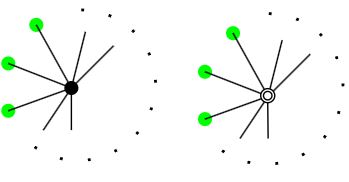}
    \caption{Two plabic graphs $SU_{k-1,k},SU_{1,k}$ such that $\mathscr{M}(SU_{k-1,k})=U_{k-1,k}$ and $\mathscr{M}(SU_{1,k})=U_{1,k}.$}
    \label{fig:S(U_k,k-1),S(U_1,k)}
\end{figure}

Now we need to define an operation on plabic graphs:
\begin{definition}\label{def:2-sum of plabic graphs}
    Let $P,P'$ be two plabic graphs, such that $l$ and $l'$ are green leaves of $P$ and $P'$, respectively. Denote the vertex in $P$ that is adjacent to $l$ as $v$, and the vertex in $P'$ that is adjacent to $l'$ as $v'$. The 2-sum of $P$ and $P'$ along  $l$ and $l'$ is a plabic graph $P \oplus_2^{l,l'} P'$, where $V(P\oplus_2^{l,l'}P')$ is the disjoint union of $V(P)\setminus \{l\}$ and $V(P')\setminus \{l'\}$. All edges in $P$ or $P'$ that are not incident to $l$ or $l'$ stay as edges in $P\oplus_2^{l,l'}P'$, and we replace the edges $\{v,l\}$ and $\{l',v'\}$ by one edge $\{v,v'\}$. This edge $\{v,v'\}$ is called the \textit{bridge} of the 2-sum. If $l$ and $l'$ are clear in the context, we may simply denote the 2-sum as $P\oplus_2 P'.$
\end{definition}

\begin{figure}[h]
    \centering
    \includegraphics[width=0.65\linewidth]{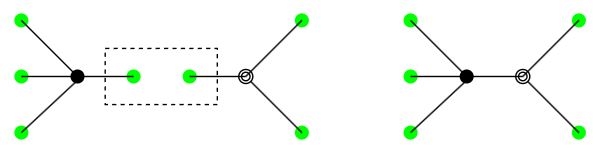}
    \caption{2-sum of two plabic graphs}
    \label{fig:2-sum of plabic graphs}
\end{figure}

In Figure \ref{fig:2-sum of plabic graphs}, we illustrate the operation of 2-sum on two plabic graphs. 

For completeness of information, it is worth pointing out that there is an operation called 2-sum on any two matroids on at least two elements whose ground sets have a common element $p$, such that neither of them has $p$ as a loop/coloop. We are going to see that it is not a coincidence for it to share the same name as the 2-sum of plabic graphs we mention above. Let $M,N$ be two matroids satisfying the specified condition. Then there is a matroid on the ground set $(E(M)\cup E(N))- p$ with set of circuits
$$\mathscr{C}(M\setminus p)\cup \mathscr{C}(N\setminus p)\cup \{(C\cup D)-p: p\in C\in \mathscr{C}(M) \mbox{ and } p\in D\in \mathscr{C}(N)\}.$$
This matroid is called the 2-sum of $M$ and $N$ along $p$ (see \cite[p.308]{Seymour} and \cite[Proposition 7.1.20]{Oxley}), denoted by $M\oplus_2^p N$. It is often simply denoted by $M\oplus_2 N$ if there is no ambiguity. Its set of bases $\mathscr{B}(M\oplus_2 N)$ can be described as 
\begin{align*}
&\{(B_1\setminus \{p\})\cup B_2: p\in B_1\in \mathscr{B}(M)\mbox{ and } p\notin B_2\in \mathscr{B}(N)\}
\\ &\cup \{B_3\cup (B_4\setminus \{p\}): p\notin B_3\in \mathscr{B}(M)\mbox{ and } p\in B_4\in \mathscr{B}(N)\}.
\end{align*}
\begin{proposition}  
Suppose we have two plabic graphs $P$ and $P'$ with green leaves $l$, and $l'$, respectively and $l,l'$ are not loop/coloop of $\mathscr{M}(P),\mathscr{M}(P').$ Relabel $l$ and $l'$ in $\mathscr{M}(P)$ and $\mathscr{M}(P')$ as $p$, then $\mathscr{M}(P\oplus_2 P')\cong \mathscr{M}(P)\oplus_2 \mathscr{M}(P').$
\end{proposition}
\begin{proof}
    If $B$ is a base of $\mathscr{M}(P\oplus_2 P')$, there is a valid orientation $\mathscr{O}$ of $P\oplus_2 P'$ such that every incident edge of a vertex of $B$ in $P\oplus_2 P'$ is directed out from that vertex. Let $v$ (resp. $v'$) be the vertex adjacent to $l$ (resp. $l'$) in $P$ (resp. $P'$). If the bridge is directed as $(v,v')$ in $\mathscr{O}$, orient $\{v,l\}$ as $(v,l)$  and all other edges of $P$ as in $\mathscr{O}$, then we have a valid orientation of $P$ giving $B\cap E(\mathscr{M}(P))$ as a base for $\mathscr{M}(P)$. For $P'$, orient $\{v',l'\}$ as $(l',v')$ and all other edges as in $\mathscr{O}$. That results in a valid orientation of $P'$ giving $(B\cap (\mathscr{M}(P')))\cup \{l'\}$ as a base for $\mathscr{M}(P').$ On the other hand, if the bridge is directed in the other way, the argument is analogous, with the roles of $\{v,l\}$ and $\{v',l'\}$ exchanged.

    For the other direction, let $B$ be a base of $\mathscr{M}(P)\oplus_2 \mathscr{M}(P')$. Without loss of generality, assume that $B=(B_1\setminus \{p\})\cup B_2$ such that $p\in B_1\in \mathscr{M}(P)$ and $p\notin B_2\in \mathscr{M}(P')$. Let $\mathscr{O}_B$ and $\mathscr{O}_B'$ be valid orientations of $P$ and $P'$, corresponding to $B_1$ and $B_2$, respectively. Then we can construct a valid orientation for $P\oplus_2 P'$ as follows: For an edge which is not the bridge of $P\oplus_2 P'$, let it inherit the orientation in $\mathscr{O}_B$ or $\mathscr{O}_B'$ depending on whether it is an edge in $P$ or $P'$. For the bridge, direct it towards $P.$ 
\end{proof}

It is important to emphasize that, in general, when we perform a 2-sum on plabic graphs $P,P'$ along their leaves $l,l'$, we do not always require that neither $l$ nor $l'$ be a loop or a coloop in $\mathscr{M}(P),\mathscr{M}(P').$

Now we proceed by introducing a class of matroids called series-parallel matroids, the definition of which can be found in Minty's paper \cite{Minty} or Welsh \cite[p.235]{Welsh}. After that, we see that there is a relationship between series-parallel matroids and the class of plabic matroids formed by successively 2-summing certain kinds of plabic graphs. 

\begin{definition}\label{def:series and parallel extensions}
Let $M$ be a (not necessarily plabic) matroid on a set $S$. Let $x\in S$ and suppose $y\notin S$. The \textit{series extension} of $M$ at $x$ by $y$ is the matroid on $S\cup \{y\}$ which has its bases the sets of the form (1) or (2)
\begin{enumerate}[label=(\arabic*)]
    \item $B\cup \{y\}; B \mbox{ is a base of }M,$
    \item $B\cup \{x\}; B \mbox{ is a base of } M, x\notin B.$
\end{enumerate}
The \textit{parallel extension} of $M$ at $x$ by $y$ is the matroid on $S\cup \{y\}$ which has its bases the sets of the form (3) or (4)
\begin{enumerate}[label=(\arabic*),start=3]
    \item $B \mbox{ is a base of }M,$
    \item $(B\setminus x)\cup \{y\}; x\in B, B \mbox{ is a base of }M.$
\end{enumerate}
\end{definition}

A \textit{series-parallel matroid}, as defined by Minty \cite{Minty}, is a matroid obtained by successively performing series/parallel extensions on a matroid on one element. By convention, a matroid on one element is also series-parallel. Minty also pointed out that a series-parallel matroid is the cycle matroid of a planar graph. A matroid on one element is either the cycle matroid of an edge incident to two vertices, or an edge incident to only one vertex (that is a loop). Suppose we have a cycle matroid $M$ of a graph $G$. Note that the series extension of $M$ at $x$ by $y$ is the cycle matroid of a graph obtained by subdividing the edge $x$ of $G$ into edges $x$ and $y$. The parallel extension of $M$ at $x$ by $y$ is the cycle matroid of a graph obtained by adding a parallel edge $y$ to the edge $x$. Obviously, a series-parallel matroid need not be connected. We may just begin with an edge incident to two vertices and perform a series extension. Then the cycle matroid induced by the resulting graph is the free matroid on two elements, which is clearly not connected.   

It is important to note that series-parallel matroids in the sense of Oxley \cite{Oxley} must be connected \cite[Corollary 5.4.12]{Oxley}. Before further comparing the definitions of series-parallel matroids we adopt (in Minty's sense) and that of Oxley, we need some more tools:

Let $M_1$ and $M_2$ be matroids on disjoint sets $E_1$ and $E_2$. Let $E=E_1\cup E_2$ and $\mathcal{I}=\{I_1\cup I_2 : I_1\in \mathcal{I}(M_1) \mbox{ and } I_2\in \mathcal{I}(M_2)\}.$ Then $(E,\mathcal{I})$ is a matroid called the direct sum of $M_1$ and $M_2$, and it is denoted by $M_1\oplus M_2$ (see \cite[Proposition 4.2.8]{Oxley} or \cite[p.72]{Welsh}). 

\begin{proposition}[{\cite[Proposition 4.2.7]{Oxley}}] \label{prop:disconnected iff direct sum}
A matroid is disconnected if and only if it is a direct sum of non-empty matroids.
\end{proposition}

A series-parallel matroid in Oxley's sense is a cycle matroid of a graph, which is obtained by successively adding a parallel edge to an edge or subdividing an edge to either of the two single edge graphs. In particular, we are allowed to subdivide an edge only if the edge is not a cut edge (an edge of a graph is called a cut edge if removing the edge increases the number of connected components of the graph) \cite[p.151]{Oxley}. In other words, we can only perform series extension on an edge $x$ if $x$ is not a coloop of $M(G)$. An important characterization of this class of matroids is the following:

\begin{proposition}[\cite{Brylawski}]\label{prop:exclused minor series-parallel}
    A matroid $M$ is a series-parallel matroid in the sense of Oxley if and only if $M$ is connected and does not contain any minor isomorphic to $M(K_4)$ or $U_{2,4}.$
\end{proposition}

Bear in mind that whenever we mention series-parallel matroid, we mean series-parallel matroid in the sense of Minty, unless we say otherwise. Given a series-parallel matroid $M$, we can write $M$ as direct sum of connected components. By the preceding proposition and the fact that both $M(K_4)$ and $U_{2,4}$ are connected matroids, every series-parallel matroid is a direct sum of series-parallel matroids in the sense of Oxley. As a result, the class of direct sums of series-parallel matroids coincides with the class of direct sums of series-parallel matroids in Oxley's sense. This class is called \textit{quasi series-parallel}. We point out this fact because many authors use the definition of Oxley (for instance \cite{Ferroni},\cite{Proudfoot}).

The relations between series or parallel extensions on plabic matroids and 2-summing of plabic graphs are described in the following two lemmas:

\begin{lemma}\label{lemma:series=2 summming star}
    Let $P$ be a plabic graph and $x$ a green leaf of $P$. The series extension of $\mathscr{M}(P)$ at $x$ by $y$ is isomorphic to $\mathscr{M}(P')$, where $P'$ is the plabic graph obtained by 2-summing $P$ with $SU_{2,3}$ along $x$ and one of the leaves of $SU_{2,3}.$ 
\end{lemma}
\begin{proof}
    Denote the center of $S_{2,3}$ by $c$. After 2-summing $P$ with one of the leaves of $S_{2,3}$, label one of the leaves in $P'$ that is adjacent to $c$ as $x$ and another as $y$. Let $B$ be a base of $\mathscr{M}(P)$. There must be a valid orientation $\mathscr{O}_B$ of $P$ in which every edge incident to a vertex $b\in B$ is directed out from $b$. If $x\in B$, then in $\mathscr{O}_B$ the edge incident to $x$ is directed out from $x$. Let all edges of $P'$ not incident to $c$ inherit the orientation of $\mathscr{O}_B$ and direct the bridge (recall Definition \ref{def:2-sum of plabic graphs}) out from $c$. Then we have to direct the edge $\{x,c\}$ as $(x,c)$, $\{y,c\}$ as $(y,c)$. Now we have a valid orientation for $P'$ such that the base of $\mathscr{M}(P')$ corresponding to it is $B\cup \{y\}$. 
    
    On the other hand, if $x\notin B$, there must be a valid orientation $\mathscr{O^*}_B$ in which the edge incident to $x$ is directed towards $x$. Similarly, let all edges of $P'$ not incident to $c$ inherit the orientation of $\mathscr{O^*}_B$ and direct the bridge towards $c$. If we direct the edge $\{x,c\}$ as $(c,x)$ and $\{y,c\}$ as $(y,c)$, then we have a valid orientation for $P'$ such that its corresponding base in $\mathscr{O^*}_B$ is still $B\cup \{y\}$. So we have finished the case (1) in Definition \ref{def:series and parallel extensions}. Moreover, if we direct the edge $\{x,c\}$ as $(x,c)$ and $\{y,c\}$ as $(c,y)$, we still get a valid orientation for $P'$ such that its corresponding base in $\mathscr{O^*}_B$ is $B\cup \{x\}$. This is (2) in Definition \ref{def:series and parallel extensions}.
\end{proof}

The proof of the following lemma is analogous to that in the previous lemma. So we are not going to repeat that.

\begin{lemma}\label{lemma:series/parallel are 2-sum}
    Let $P$ be a plabic graph and $x$ a green leaf of $P$. The parallel extension of $\mathscr{M}(P)$ at $x$ by $y$ is isomorphic to $\mathscr{M}(P')$, where $P'$ is the plabic graph obtained by 2-summing $P$ with $SU_{1,3}$ along $x$ and one of the leaves of $SU_{1,3}.$
\end{lemma}

The following three lemmas describe some operations we may perform on a plabic graph without changing the matroid induced. 

\begin{lemma}\label{lemma:contra edge same color same plabic}
    If we contract an edge between two vertices that are both colored black or white in a plabic graph $P$, then the matroid induced by the new plabic graph is the same as $\mathscr{P}$.
\end{lemma}
\begin{proof}
    Suppose we have two adjacent black vertices $b_1,b_2$ in $F$. Let $F'$ be the plabic graph obtained from contracting the edge between $b_1,b_2$, and denote the vertex in $F'$ obtained by identifying $b_1,b_2$ by $b$. If $\mathscr{O}_F$ is a valid orientation for $F$, just let the edges of $F'$ inherit their directions in $\mathscr{O}_F$. This gives us a valid orientation for $F'$ that gives us the same base. Conversely, if we have a valid orientation $\mathscr{O}_{F'}$, there is exactly one edge directed out from $b$ in $\mathscr{O}_{F'}$. Let that edge be $(b,x)$. Now let all edges in $F$ except $\{b_1,b_2\}$ inherit their directions in $\mathscr{O}_{F'}$. If $x$ is adjacent to $b_1$ in $F$, direct $\{b_1,b_2\}$ as $(b_2,b_1)$ and reverse that direction if $x$ is adjacent to $b_2$ in $F$. Then we get a valid orientation for $F$ that gives us the same base. The argument for two adjacent white vertices is completely analogous. 
\end{proof}

\begin{lemma}\label{lemma: contracting out deg 2 vert.}
    If $v$ is a black/white vertex adjacent to only $v_1$ and $v_2$ in a plabic graph $P$, form a new plabic graph $P'$ by erasing $v$ and replacing the edges $\{v_1,v\}$ and $\{v,v_2\}$ by an edge $\{v_1,v_2\}$. The matroid induced by $P'$ is the same as $\mathscr{M}(P)$. We call this procedure ``contracting out'' the vertex $v$.
\end{lemma}
\begin{proof}
    In any valid orientation $\mathscr{O}$ of $P$ there must be one edge directed towards $v$ and the other directed out from it. For orientations $(v_i,v),(v,v_j)$, direct the edge $\{v_i,v_j\}$ towards $v_j$ with $\{i,j\}=\{1,2\}$ and all the other edges agree with the directions in $\mathscr{O}$.
\end{proof}

\begin{lemma}\label{lemma: black/white leaf loop/co}
    Let $v$ be a black (white) leaf in a plabic graph $P$, and $u$ be a white (black) vertex adjacent to $v$. Moreover, let $\{l_1,\ldots,l_n\}$ be the set of green leaves adjacent to $u$. Form a plabic graph $P'$ by:
    \begin{enumerate}
        \item Remove $\{l_1,\ldots,l_n\}$ and the edges incident to them.
        \item Take the disjoint union of the plabic graph obtained after the above process with $n$ copies of $SU_{0,1}$ (resp. $SU_{1,1}$).
    \end{enumerate}
    Then $\mathscr{M}(P')=\mathscr{M}(P).$
\end{lemma}
\begin{proof}
    Because in any valid orientation of $P$, $v$ must be directed towards (out from) $u$. All edges incident to $\{l_1,\ldots,l_n\}$ must be directed towards (out from) them. So they are loops (coloops) in $\mathscr{M}(P)$. Step (1) removes all of them from $\mathscr{M}(P)$ without affecting other elements in the matroid. Step (2) just adds them back to the matroid, so $\mathscr{M}(P')=\mathscr{M}(P).$
\end{proof}

\begin{proposition}\label{prop:direct sum of series-parallel iff plabic forest}
    A matroid arises from a nonempty plabic forest if and only if it is a quasi series-parallel matroid.
\end{proposition}
\begin{proof}
    We prove the only if part first. 
    
    \noindent Let us tackle some special cases:

    \noindent -If $F$ has a connected component containing no black/white vertex:
    \\ That connected component is a disjoint union of copies of ``two green leaves adjacent to each other". Note that for each copy, we can view it as follows:
    \begin{figure}[h]
        \centering
        \includegraphics[width=0.95\linewidth]{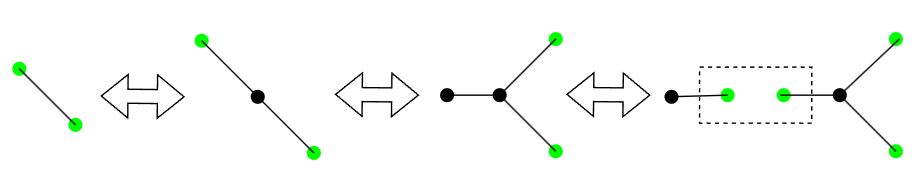}
        \caption{Two green leaves adjacent to each other yield a series-parallel matroid.}
        \label{fig:Only 2 leaves series-parallel}
    \end{figure}
 
  The dotted rectangular box indicates the operation of 2-sum along the leaves inside. The bi-directional arrows mean that the plabic matroids they induce are the same. The matroid induced by them is series-parallel by Lemma \ref{lemma:series=2 summming star}.

\noindent -If there is a black/white leaf in $F$: 
\\ If the connected component $H$ of $F$ containing the black/white vertex does not have any green leaf, then we can forget that component without affecting $\mathscr{M}(P).$ If $H$ only consists of the black/white vertex and a green leaf attached to it, then $H$ corresponds to a loop/coloop in $\mathscr{M}(P)$. Recall that by convention, a loop/coloop is series-parallel. Now, we can assume that there is a black leaf $b$ in $H$. If there is an edge incident to two black or two white vertices, we can make the black/white induced subgraph of $H$ to have alternating vertex colorings. So $b$ is adjacent to some white vertex $w$. If there is a green leaf adjacent to $w$, apply Lemma \ref{lemma: black/white leaf loop/co} to remove it and compensate that by taking disjoint union with $SU_{0,1}$. Thus we may assume that there is no green leaf adjacent to $w$. Observe that in any valid orientation of $F$, there is only one possible orientation for the edge $\{b,w\}$ and all other edges incident to $w$ must be directed out from $w$. So we may erase vertices $b,w$ and all edges incident to $w$ without changing $\mathscr{M}(F).$ We can also get rid of white leaf by a similar argument. Figure \ref{fig:erase black/white leaves} is an illustration:
\begin{figure}[h]
    \centering
    \includegraphics[width=1\linewidth]{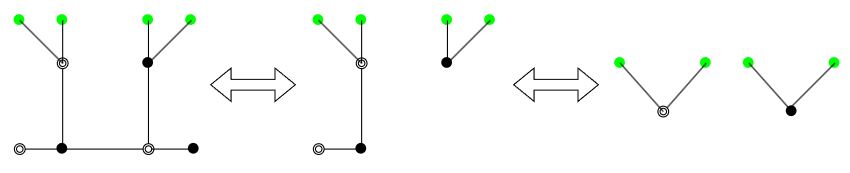}
    \caption{Erase black/white leaves}
    \label{fig:erase black/white leaves}
\end{figure}

Therefore, we restrict our attention to plabic trees that have no black/white leaves. Suppose we have a plabic matroid $\mathscr{M}(F)$ where $F$ is a nonempty plabic forest. Note that if we have a vertex of degree 4, we can perform the process depicted in Figure \ref{fig:focus on internal vertex of dig le 3} due to Lemma \ref{lemma:contra edge same color same plabic}.
    \begin{figure}[h]
        \centering
        \includegraphics[width=0.6\linewidth]{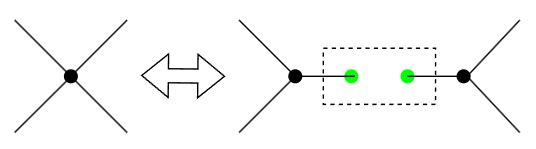}
        \caption{A degree 4 vertex can be decomposed.}
        \label{fig:focus on internal vertex of dig le 3}
    \end{figure}
    
    If we have a vertex of degree strictly bigger than 4 in $F$, we can apply the process described in Figure \ref{fig:focus on internal vertex of dig le 3} successively so we may assume that every vertex in $F$ is of degree less than or equal to three.

Furthermore, by Lemma \ref{lemma: contracting out deg 2 vert.}, we may assume that every black/white vertex is of degree 3. We proceed by induction on the number of internal vertices. The base case is clearly true by the following observation (together with Figure \ref{fig:Only 2 leaves series-parallel}):
\begin{figure}[h]
    \centering
    \includegraphics[width=0.65\linewidth]{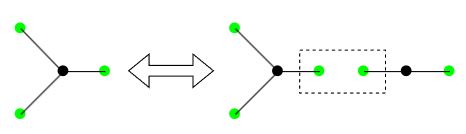}
    \caption{There is only one black vertex.}
    \label{fig:base case plabic tree is series-parallel}
\end{figure}

\noindent Our induction hypothesis is that a plabic tree with $n$ internal vertices induces a series-parallel matroid. Suppose we have a plabic tree $T$ with $n+1$ internal vertices. Note that there must be an internal vertex in $T$ which is adjacent to only one other internal vertex. To see this, just look at a longest path $P$ in $T$. The endpoints in this path must be green leaves and we denote them by $l_1$ and $l_2$. Let $v$ be a vertex in this path that is adjacent to $l_2$. Then $v$ is adjacent to exactly one other internal vertex $v'$. If $v$ were adjacent to another internal vertex $u\neq v'$, then we would get a path $(l_1,\ldots,v',v,u,\ldots)$ that is longer than the original longest path $P=(l_1,\ldots,v',v,l_2)$. This is a contradiction. Now, note that we can view $T$ as a 2-sum of $SU_{2,3}$ or $SU_{1,3}$ (depending on whether $v$ is black or white) with a plabic tree with $n$ internal vertices. Our desired conclusion follows from the induction hypothesis.  

The converse of this Proposition is a direct consequence of Lemma \ref{lemma:series=2 summming star} and Lemma \ref{lemma:series/parallel are 2-sum}.
\end{proof}

We can begin to characterize transversal matroids that admit acyclic presentations in terms of series-parallel matroids. Given a transversal matroid $M$, due to Lemma \ref{lemma: exists |G_M|=|base|}, there must be a bipartite presentation $G$ for $M$ such that $|G_M|=\operatorname{rank}(M).$ We form a plabic presentation for $M$ from $G$ as follows: color vertices in $G_N$ as black and vertices in $G_M$ as white. Attach a green leaf to each black vertex (if there are loops, remove them and compensate that by direct sum). Let $B\subseteq G_N$ be a base of $M$, orient the edges between those green leaves adjacent to vertices in $B$ towards $B$, and orient the matching edges of $B$ to $G_M$ towards $G_M$. For other edges between $G_M$ and $G_N$, orient them towards $G_N$. For other vertices in $G_N$ which are not in $B$, orient the edges between them and green leaves towards the leaves. The resulting orientation is a valid orientation of the plabic graph. It is easy to see that each valid orientation gives us back a base of $M$ in the same way. The following diagram is an illustration. 
\begin{figure}[h]
    \centering
    \includegraphics[width=0.6\linewidth]{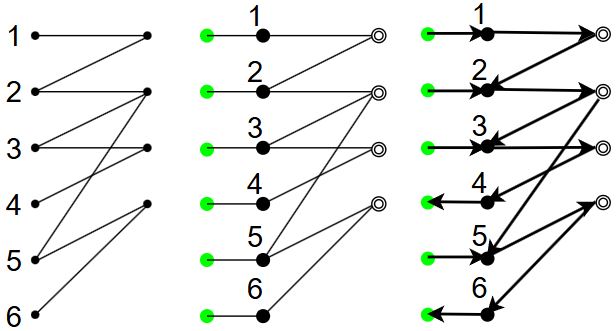}
    \caption{Turn bipartite presentation of $M$ to a plabic presentation.}
    \label{fig:tran. M. as plabic mat.}
\end{figure}

\begin{figure}[h]
    \centering
    \includegraphics[width=0.36\linewidth]{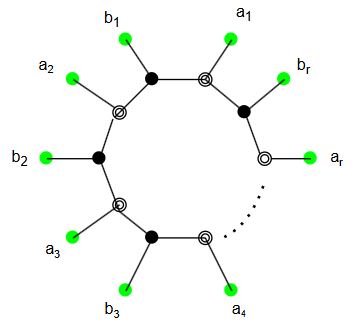}
    \caption{A plabic whirl of order $r$.}
    \label{fig:Whirl as plabic matroid}
\end{figure}

\noindent Therefore, if $M$ has an acyclic (bipartite) presentation, there is a plabic forest inducing $M$. By Proposition \ref{prop:direct sum of series-parallel iff plabic forest}, $M$ is a direct sum of series-parallel matroids. We want to know whether the converse is also true. We call the plabic graph shown in Figure \ref{fig:Whirl as plabic matroid} a \textit{plabic whirl of order $r$}\footnote{We choose this name because the matroid induced by this plabic graph is in fact isomorphic to the graphic matroid that is called the whirl of order $r$ in the literature \cite[Proposition 1.5.14, Example 8.4.3]{Oxley}.}, denoted by $\mathscr{P}_r$.

\newpage
\begin{lemma}\label{lemma:2-plabic whirl is U_{2,4}}
    The matroid induced by the plabic whirl of order 2 is isomorphic to $U_{2,4}$.
\end{lemma}
\begin{proof}
    There is a valid orientation for every two-element set of $[4]$ such that the two-element set is a base. Those valid orientations are listed in Figure \ref{fig:plabic whirl 4 is iso to U2,4}. It is straightforward to see that there is no valid orientation that can give a three-element set as a base.
\end{proof}

\begin{figure}[h]
    \centering
    \includegraphics[width=0.73\linewidth]{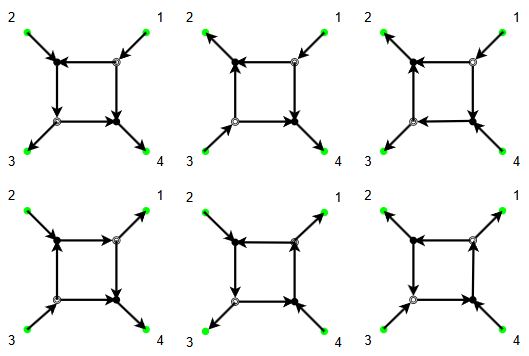}
    \caption{The matroid induced by the plabic whirl of order 2 is isomorphic to $U_{2,4}$.}
    \label{fig:plabic whirl 4 is iso to U2,4}
\end{figure}

\begin{proposition}\label{prop:big whirl has small one as minor}
    $\mathscr{M(P}_r)$ contains $\mathscr{M(P}_{r-1})$ as a minor for all $r\ge 3.$
\end{proposition}
\begin{proof}
    We adopt the notation in Figure \ref{fig:Whirl as plabic matroid}. We claim that $\mathscr{M(P}_r)/b_r\setminus a_r$ is isomorphic to $\mathscr{M(P}_{r-1})$. In $\mathscr{P}_r$, denote the white (black) vertex adjacent to $a_i$ (resp. $b_i$) by $x_i$ (resp. $y_i$) for all $1\leq i\leq r$. Arbitrarily pick a valid orientation $\mathscr{O}$ of $\mathscr{P}_{r-1}$ and let $B_{\mathscr{O}}$ be the base of $\mathscr{M(P}_{r-1})$ corresponding to $\mathscr{O}$. The orientation of the edge $\{x_1,y_{r-1}\}$ is either $(x_1,y_{r-1})$ or the reverse direction. In the former case, we form a valid orientation on $\mathscr{P}_r$ as follows:

    \noindent 
    \begin{enumerate}
    \item Orient $\{x_1,y_r\}$, $\{x_r,y_r\}$ and $\{x_r,y_{r-1}\}$ as $(x_1,y_r)$, $(y_r,x_r)$ and $(x_r,y_{r-1})$, respectively
    \item Orient $\{b_r,y_r\}$ and $\{x_r,a_r\}$ as $(b_r,y_r)$ and $(x_r,a_r)$, respectively.
    \item Let all other edges of $\mathscr{P}_r$ inherit the directions in $\mathscr{O}$. 
    \end{enumerate}
    Then we end up with a valid orientation of $\mathscr{P}_r$. This valid orientation corresponds to the base $B_{\mathscr{O}}\cup \{b_r\}$ of $\mathscr{M(P}_r).$ 

    On the other hand, if in $\mathscr{O}$, the edge $\{x_1,y_{r-1}\}$ is oriented as $(y_{r-1},x_1)$, then just reverse the directions of edges in 1, and keep the directions of edges in 2 and 3. This again results in a valid orientation corresponding to the base $B_{\mathscr{O}}\cup \{b_r\}$ of $\mathscr{M(P}_r).$

    If instead we are given a valid orientation $\mathscr{O'}$ of $\mathscr{M}(\mathscr{P}_r)$ such that the edge incident to $b_r$ (resp. $a_r$) is directed out (resp. towards $a_r$), then the set of green leaves excluding $b_r$ in $\mathscr{P}_r$ such that their incident edges are directed out from them is a base $B_{\mathscr{O'}}$ of $\mathscr{M}(\mathscr{P}_r)/b_r\setminus a_r$. Depending on whether $y_r$ is directed towards $x_r$ or vice versa, we can perform the reverse of the process above. For instance, for the former case, $x_1$ must be directed to $y_r$ and $x_r$ must be directed to $y_{r-1}$ in $\mathscr{O'}$. Then in $\mathscr{P}_{r-1}$, we direct $x_1$ towards $y_{r-1}$ and let other edges inherit the orientations of $\mathscr{O'}$. The resulting orientation is a valid one for $\mathscr{P}_{r-1}$. The base of $\mathscr{M}(\mathscr{P}_{r-1})$ corresponding to that coincides with $B_{\mathscr{O'}}.$ The case when $x_r$ is directed towards $y_r$ in $\mathscr{O'}$ is analogous. 
\end{proof}

\begin{corollary}\label{coro:whirl is not quasi series-parallel}
    $\mathscr{M(P}_r)$ with $r\ge 2$ is not a quasi series-parallel matroid.
\end{corollary}
\begin{proof}
    This follows from Lemma \ref{lemma:2-plabic whirl is U_{2,4}}, Proposition \ref{prop:big whirl has small one as minor}, Proposition \ref{prop:exclused minor series-parallel}, and the fact that $U_{2,4}$ is connected.
\end{proof}

Before we show that if $M$ is a quasi series parallel, transversal matroid then $M$ has an acyclic bipartite presentation, we need some more tools. We say that a bijection $\varphi$ between two matroids $M_1,M_2$ is a \textit{weak map} from $M_1$ to $M_2$ if $\varphi^{-1}(I)$ is independent in $M_1$ for every independent set $I$ of $M_2$ (see \cite[p.279]{Oxley}, \cite{Kung} or \cite{Higgs}). The collection of matroids on a fixed set $E$ can be partially ordered by taking $M_1\ge M_2$ if the identity map is a weak map from $M_1$ to $M_2$. This is called the \textit{weak order} on the collection \cite[p.281]{Oxley}, and $M_2$ is called a \textit{weak image} of $M_1$. Moreover, if the rank of $M_1$ equals that of $M_2$, we say that $M_2$ is a \textit{strong weak image} of $M_1.$

\begin{figure}[h]
    \centering
    \includegraphics[width=0.3\linewidth]{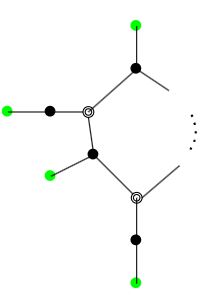}
    \caption{Plabic presentation of $M$ admitting no acyclic presentation contains a plabic whirl $S$ as a subgraph.}
    \label{fig:contains whirl as weak image}
\end{figure}

Given a transversal matroid $M$ with no acyclic bipartite presentation, let $G$ be a bipartite graph representing $M$. Again we assume $|G_M|=\operatorname{rank}(M).$ Color $G_N$ as black and $G_M$ as white, and attach a green leaf to each black vertex. We have a plabic graph $P$ such that $\mathscr{M}(P)=M.$ The adjacent vertices in the black/white induced subgraph of $P$ have alternating vertex colors. By contracting out vertices of degree 2 and edges between two vertices of the same color, we may assume that each internal vertex is of degree at least 3. Then it contains the plabic whirl graph $S$ in Figure \ref{fig:contains whirl as weak image} as a subgraph, where the cycle in the subgraph has length at least 4. Let the matroid induced by $S$ be $W$. 

Now, we want to show that for any valid orientation $\mathscr{O}_S$ in $S$, we can extend the orientation to a valid orientation $\mathscr{O}_S^P$ of $P$. In other words, we want to get a valid orientation $\mathscr{O}$ of $P$ such that when we restrict it to $S$, the orientation coincides with that of $\mathscr{O}_S$. We first let the subgraph $S$ inherit the orientation of $\mathscr{O}_S$. And we show that we may orient all other edges in $P$ not in $S$ so that the whole orientation is valid in $P$. It is easier for readers to follow by using a diagram:
\begin{figure}[h]
    \centering
    \includegraphics[width=0.3\linewidth]{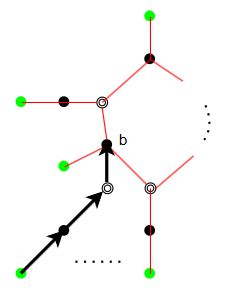}
    \caption{Extend the orientation $\mathscr{O}_S$.}
    \label{fig:extend the orientation of O_I}
\end{figure}

\noindent Figure \ref{fig:extend the orientation of O_I} shows part of $P$. The edges colored in red are edges in the subgraph $S$, and we assume they already inherited the orientation of $\mathscr{O}_S$. Suppose there is a path from a leaf to a black vertex $b$ in $S$, and the whole path except $b$ is not in $S$. We orient all edges in that path towards $b$. This does not affect the validity of the orientation inside $S$. 
\begin{figure}[h]
    \centering
    \includegraphics[width=0.28\linewidth]{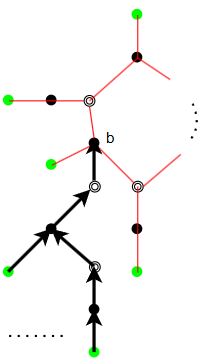}
    \caption{Further extend the orientation $\mathscr{O}_S$.}
    \label{fig:further extend the orientation of O_I}
\end{figure}

If there is any other path branching out, we can orient it in a similar way as shown in Figure \ref{fig:further extend the orientation of O_I}. If the path is branching out from a white vertex, just orient the path out from the white vertex. It is straightforward to see, by induction, that this can be done for any finite number of branching paths. Our argument can extend $\mathscr{O}_S$ to a valid orientation of a spanning set (every vertex of $P$ is in there) of edges of $P$, in which the red cycle is the only cycle. If there is any edge that is not in the validly oriented spanning set, since the edge must be incident to a black and a white vertex, just orient the edge towards the black vertex. Then we get a valid orientation $\mathscr{O}_S^P$ of $P$. 

Let $L$ be the set of green leaves of $P$ not in $S$, $L_{out}^{\mathscr{O}_S^P}$ (resp. $L_{in}^{\mathscr{O}_S^P}$) be the subset of $L$ such that their incident edges are directed out from them (resp. towards them) in $\mathscr{O}_S^P$ for some valid orientation $\mathscr{O}_S$ of $S$. However, note that no matter what valid orientation $\mathscr{O}_S$ we choose for $S$, during the process we construct $\mathscr{O}_S^P$, the orientations of all edges in $P$ that are not in $S$ are fixed. That is because for an edge between an endpoint in $P\setminus S$ and a black (resp. white) vertex in $S$, we always direct the edge towards (resp. out from) the vertex in $S$. 
So $L_{out}^{\mathscr{O}_S^P}$ (resp. $L_{in}^{\mathscr{O}_S^P}$) can be denoted by $L_{out}$ (resp. $L_{in}$).  By construction, every base of $\mathscr{M}(S)$ is a base of $\mathscr{M}(P)/L_{out}\setminus L_{in}.$ In other words, $\mathscr{M}(S)$ is a rank-preserving weak image of $\mathscr{M}(P)/L_{out}\setminus L_{in}$.
  
  Due to Proposition \ref{prop:exclused minor series-parallel} and the fact that both $U_{2,4}$ and $M(K_4)$ are connected, a minor of a quasi series-parallel matroid is quasi series-parallel. Moreover, a rank-preserving weak image of a quasi series-parallel matroid is quasi series-parallel \cite[Theorem 6.5, Theorem 6.10]{Lucas}. Since $M=\mathscr{M}(P)$ is quasi series-parallel, this would imply that $\mathscr{M}(S)$ is quasi series-parallel, a contradiction to Corollary \ref{coro:whirl is not quasi series-parallel}. Thus, $M$ is not a quasi series-parallel matroid. In conclusion, we have proved:

  \begin{theorem}
      A transversal matroid $M$ admits an acyclic bipartite presentation if and only if $M$ is quasi series-parallel and transversal.
  \end{theorem}

Because this is outside the scope of this dissertation, the author is not going to give full details here. One may further characterize which quasi series-parallel matroids are transversal. Given a series-parallel matroid $M$, there is a plabic tree presentation $T$ for $M$. The following moves in $T$ do not affect the matroid induced by the plabic graph:
\begin{enumerate}[(1)]
    \item Contracting edges between two internal vertices of the same color. 
    \item Inserting a black vertex of degree 2 between an edge $\{l,w\}$ where $l$ is a green leaf and $w$ is a white vertex.
    \item If there are green leaves $\{l_1,\ldots,l_n\},n\ge 2$ adjacent to a black vertex $b$, insert a black vertex $b'$ and a white vertex $w'$ of degree 2 so that $\{l_i,b'\},\{b',w'\},\{w',b\}$ become edges. 
\end{enumerate}

\noindent By doing these moves if necessary, we may assume the plabic tree representing $M$ has the following form: the adjacent internal vertices have alternating colors and each green leaf is adjacent to exactly one black vertex. If we can also make each black vertex adjacent to exactly one green leaf, then the black/white induced subgraph is a bipartite presentation of $M$. Thus, $M$ is transversal. On the other hand, suppose there is a black vertex of degree at least 3 and it is not adjacent to any green leaf. If we are not able to get rid of it by all of the above moves, then $M$ is not transversal. Suppose such an $M$ is transversal, then $M$ admits an acyclic bipartite presentation by the preceding Theorem. By attaching one green leaf to each black vertex, we get another plabic tree presentation which has the desired form. One can derive a contradiction by showing that the plabic tree presentation for a matroid is unique up to the moves (1), (2) and (3).

\section{Acyclic presentations for a transversal matroid $M$}
In this section, our main result will let us see how to spot $C\in \mathscr{C}(M)$ in $G_{\mathscr{A}}$ if $\mathscr{A}$ is an acyclic presentation for a transversal matroid $M$. Such circuits can be described as `plantings' of particular kinds of subgraphs in $G_{\mathscr{A}}.$
\vspace{2mm}

\begin{lemma}\label{lemma: cir. u. of cirs. iff no coloop}
    If $C$ is a subset of the ground set of a transversal matroid $M$ such that $|(C_{\mathscr{A}})_N|\;=|(C_{\mathscr{A}})_M|+1$ for some presentation $\mathscr{A}$ of $M$, and $C$ is a union of circuits of $M$, then $C\in \mathscr{C}(M).$
\end{lemma}
\begin{proof}
    A coloop $x$ in $M$ is a loop in $M^*$. Then no basis of $M^*$ contains $x$. That is equivalent to saying that every base of $M$ contains $x$. This means that no circuit of $M$ can contain $x$.
    
     By assumption we may assume $C=\bigcup_{i=0}^n C_i$ such that  $C_m\nsubseteq \bigcup_{i=0}^{m-1} C_i$, for all $0<m\leq n$,  where $C_0,...,C_n \in \mathscr{C}(M)$. We want to derive a contradiction to the assumption $|C|\;=|((C)_{\mathscr{A}})_M|+1$ if $n>0$, and we will proceed by induction. 
    Suppose $n=1$, that is, $C=C_0\cup C_1$, note that 
    $$|C_1\setminus (C_0\cap C_1)|\;=|C_1|-|C_0\cap C_1|$$
    and 
    $$|((C_1)_{\mathscr{A}})_M\setminus ((C_0\cap C_1)_{\mathscr{A}})_M|\;=|((C_1)_{\mathscr{A}})_M|-|((C_0\cap C_1)_{\mathscr{A}})_M|.$$
    Observe that $C_0\cap C_1$, as a proper subset of both $C_0$ and $C_1$, has a matching set in both $((C_0)_{\mathscr{A}})_M$ and $((C_1)_{\mathscr{A}})_M$. But by Remark \ref{re:(G)_M-(S)_M NO S}, any matching set of $C_0\cap C_1$ must be a subset of $((C_0\cap C_1)_{\mathscr{A}})_M$. This means $|C_0\cap C_1|\leq |((C_0\cap C_1)_{\mathscr{A}})_M|$. Together with Proposition \ref{prop: cir. N=M+1}, we have 
    $$|((C_1)_{\mathscr{A}})_M\setminus ((C_0\cap C_1)_{\mathscr{A}})_M|\;< |C_1\setminus (C_0\cap C_1)|.$$
    Intuitively, it means if we begin with $C_0$, then take the union of $C_1$ with $C_0$, we are adding strictly more vertices to the $N$ side of $((C_0)_{\mathscr{A}})$ to form $((C_0\cup C_1)_{\mathscr{A}})_N$ than to the $M$ side of $((C_0)_{\mathscr{A}})$ to form $((C_0\cup C_1)_{\mathscr{A}})_M$. Thus, 
    $$|C_0\cup C_1|\;>|((C_0\cup C_1)_{\mathscr{A}})_M|+1.$$
    Our induction hypothesis is, for $1\leq k<n$, we have
    $$|C_0\cup ...\cup C_k|\;\geq |((C_0\cup...\cup C_k)_{\mathscr{A}})_M|+1.$$
    The inductive step is essentially the same as the base step. By Proposition \ref{lemma:m can't ad. to exa. one n}, a matching set for $C_{k+1}\cap (C_1\cup...\cup C_k)$ must be a subset of $(((C_{k+1})_{\mathscr{A}})_M\cap ((C_0\cup...\cup C_k)_{\mathscr{A}})_M)$. Together with Proposition \ref{prop: cir. N=M+1}, that is, $|C_{k+1}|\;=|((C_{k+1})_{\mathscr{A}})_M|+1$, we have
    $$|C_0\cup...\cup C_{k+1}|\;>|((C_0\cup ...\cup C_{k+1})_{\mathscr{A}})_M|+1.$$
    This is a contradiction to the assumption $|C|\;=|(C_{\mathscr{A}})_M|+1.$ The only way to avoid this contradiction is when $n=0$.
  \end{proof}
In general, if $C\subseteq E(M)$ in a matroid $M$ such that $M|C$ is connected, with $|C|=\operatorname{rank}(C)+1$, then $C\in\mathscr{C}(M)$(This can be deduced easily from Proposition \ref{prop:ref.connected 3}). This may look like an analogue of our preceding lemma. But notice that the condition $|(C_{\mathscr{A}})_N|=|(C_{\mathscr{A}})_M|+1$ in our content does not imply $|C|=\operatorname{rank}(C)+1.$

\begin{corollary}
    If $C$ is a subset of the ground set of a transversal matroid $M$ such that $|(C_{\mathscr{A}})_N|\;=|(C_{\mathscr{A}})_M|+1$ for some presentation $\mathscr{A}$ of $M$, then $C\in \mathscr{C}(M)$ if and only if $M|C$ contains no coloop.  
\end{corollary}
\begin{proof}
    If $M|C$ contains no coloop, arbitrarily pick an $x\in C$. Then there is a maximal independent set $I_x$ in $M|C$ such that $x\notin I_x$. This implies $I_x\cup \{x\}$ is dependent. Thus, there is a circuit $C_x\subseteq C$ with $x\in C_x$. So we may write $C=\bigcup_{x\in C}C_x.$  The desired conclusion follows from the preceding Lemma. The converse is clear.
\end{proof}

\begin{definition}
    A subgraph $S$ of a bipartite graph $G$ is called a basic-type subgraph, if for $k\ge 2$, the graph induced by $\{n_1,\ldots,n_k\}=S_N\subseteq G_N$ and $\{m_1,\ldots,m_{k-1}\}=S_M\subseteq G_M$ with $S_V=S_N\cup S_M$ contains $n_1m_1\ldots m_{k-1}n_k$, such that two vertices are adjacent when (but not necessarily only when) they are adjacent in the word $n_1m_1\ldots n_{k-1}m_{k-1}n_k$. 
    Furthermore, $S$ is called a pure basic-type subgraph if the graph induced by $S_N\cup S_M$ contains $n_1m_1\ldots m_{k-1}n_k$ such that two vertices are adjacent if and only if they are adjacent in the word $n_1m_1\ldots n_{k-1}m_{k-1}n_k$.
\end{definition}

\begin{remark}
    In the preceding definition, whether $n_1m_1\ldots m_{k-1}n_k$ is a basic type subgraph or a pure basic-type subgraph, it is possible for some element of $G_M\setminus S_M$ to be adjacent to some element of $S_N$. And it is also possible for some element of $G_N\setminus S_N$ to be adjacent to some element of $S_M$. This is because the conditions in the preceding definition only restrict adjacency between elements in $S_V$. 

    If $\mathscr{A}$ is a presentation for a transversal matroid $M$ and $S\subseteq E(M)$, and if $S_{\mathscr{A}}$ is a basic-type subgraph, or pure basic-type subgraph of $G_{\mathscr{A}}$, by Remark \ref{re:(G)_M-(S)_M NO S}, no element in $(G_{\mathscr{A}})_M\setminus (S_{\mathscr{A}})_M$ can be adjacent to any element of $S$ in $G_{\mathscr{A}}$, since $(S_{\mathscr{A}})_M$ is defined to contain all elements of $(G_{\mathscr{A}})_M$ that are adjacent to at least one element of $S$. But it is still possible for some element of $(G_{\mathscr{A}})_N\setminus S$ to be adjacent to an element in $(S_{\mathscr{A}})_M$. 
\end{remark}

\begin{figure}[h]
    \centering
    \includegraphics[width=0.3\linewidth]{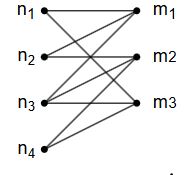}
    \caption{A basic type graph that is not a pure basic-type graph.}
    \label{fig:many basic type sub.}
\end{figure}

\begin{example}
    In Figure \ref{fig:many basic type sub.}, we have a presentation $\mathscr{A}$ of a matroid on $E=\{n_1,n_2,n_3,n_4\}$. It is easy to see that $\{n_1,n_2,n_3,n_4\}$ is a circuit of the matroid. $G_{\mathscr{A}}$ contains a basic-type subgraph $n_1m_1n_2m_2n_3m_3n_4$ (it can be read backward as $n_4m_3n_3m_2n_2m_1n_1$). The subgraph induced by $\{n_1,n_2,n_3,m_1,m_2\}$ is a basic-type subgraph, but it is not a pure basic-type subgraph. $\{n_1,n_2,n_3\}_{\mathscr{A}}$, which is defined to be the subgraph induced by $\{n_1,n_2,n_3,m_1,m_2,m_3\}$, is not a basic-type subgraph.
        Obviously, $G_{\mathscr{A}}$ is not a pure basic-type subgraph. 
\end{example}

\begin{lemma}\label{lemma:m can't ad. to exa. one n}
    If $\mathscr{A}$ is a presentation for a transversal matroid $M$, and $C\in \mathscr{C}(M)$, then in $C_{\mathscr{A}}$ every element in $(C_{\mathscr{A}})_M$ is adjacent to at least two elements of $C$.
\end{lemma}
\begin{proof}
    If $C$ is a singleton, then $C$ is an isolated vertex in $G_{\mathscr{A}}$ and $(C_{\mathscr{A}})_M$ is empty, we are finished. So suppose $C$ is not a singleton. Suppose that there is an $x\in (C_{\mathscr{A}})_M$ that is adjacent to only one element $y\in C$. Then $C\setminus \{y\}$ is an independent set and thus can be matched by a subset $J\subseteq (C_{\mathscr{A}})_M$, $x\notin J$ . But by Proposition \ref{prop: cir. N=M+1}, $|(C_{\mathscr{A}})_M\setminus \{x\}|=|C|-2$, this is absurd. 
\end{proof}

\begin{lemma}
    If $\mathscr{A}$ is a presentation for a transversal matroid $M$, and $C\in \mathscr{C}(M)$ with $|C|\geq 2$, then $C_{\mathscr{A}}$ contains a maximal (with respect to subgraph containment) pure basic type subgraph. 
\end{lemma}
\begin{proof}
    Arbitrarily choose a subset of $C$ of cardinality $|C|-1$ and denote it by $\{d_1,...,d_{|C|-1}\}$. Since it is independent and by Proposition \ref{prop: cir. N=M+1} and Remark \ref{re:(G)_M-(S)_M NO S}, it can be matched by  $(C_{\mathscr{A}})_M=\{m_1,...,m_{|C|-1}\}$. Assume that $d_i$ is matched with $m_i$, $1\leq i\leq |C|-1$. $d_{|C|}$ is adjacent to some  $m_i\in (C_{\mathscr{A}})_M$. The proof is finished by letting $A=\{d_{|C|},d_i\}$ and $B={m_i}$. We have shown the existence of a pure basic type subgraph in $C_{\mathscr{A}}$. Since $C_{\mathscr{A}}$ is a finite graph, it must contain a maximal pure basic-type subgraph.  
\end{proof}

\begin{definition}(Graph planting)
    Let $G$ and $G'$ be two undirected graphs, $x\in G_V,x'\in G'_V$, we form a graph $G*_{x,x'} G'$ such that $(G*_{x,x'}G')_V=(G_V\setminus \{x\})\cup (G_{V'}\setminus \{x'\})\cup \{y\}$, where $y\notin (G_V\setminus \{x\})$ and $y\notin (G_{V'}\setminus \{x'\})$. And $(G*_{x,x'}G')_E$ is defined as $A\cup B \cup C \cup D$, where $\\A=\{\{v_1,v_2\}\in G_E \mid \{v_1,v_2\}\cap \{x\}=\emptyset\}$, \\$B=\{\{v_1',v_2'\}\in G'_E\mid \{v_1',v_2'\}\cap \{x'\}=\emptyset\}$, 
    \\$C=\{\{y,v\}\mid \{x,v\}\in G_E \}$, and $D=\{\{y,v'\}\mid \{x',v'\}\in G'_E\}$. 
    
    $G*_{x,x'}G'$ is called the plant of $G$ from $x$ to $x'$ of $G'$, or equivalently, the plant of $G'$ from $x'$ to $x$ of $G$. And $y$ is called the joint of the planting.
\end{definition}

Intuitively, $G*_{x,x'}G'$ is the graph formed by identifying $x\in G_V$ with $x'\in G_{V'}$, while all the edges of $G$ and $G'$ are preserved. 

\begin{theorem}\label{thm:circuit, planting of pure basic type}
    If $\mathscr{A}$ is an acyclic presentation for a transversal matroid $M$ and $C\in \mathscr{C}(M)$, then $C_{\mathscr{A}}$ is a successive planting of pure basic-type subgraphs such that all joints of plantings are vertices in $G_N$.
\end{theorem}
\begin{proof}
    Arbitrarily choose a subset of $C$ of cardinality $|C|-1$ and denote it by $\{d_1,...,d_{|C|-1}\}$, and suppose it is matched by $\{m_1,...,m_{|C|-1}\}$. Without loss of generality we may assume that $m_i$ is matched to $d_i$, $1\leq i\leq |C|-1$. Since $d_{|C|}$ is not isolated, it must be adjacent to some vertex; without loss of generality, say that vertex is $m_1$. Then we ask whether $d_1$ is adjacent to any vertex in $(C_{\mathscr{A}})_M\setminus \{m_1\}$. If yes, again through renaming, without loss of generality, we may assume $d_1$ is adjacent to $m_2$. And we continue to ask whether $d_2$ is adjacent to any vertex in $((C_{\mathscr{A}})_M\setminus \{m_1,m_2\})$. Since $(C_{\mathscr{A}})_M$ is a finite set and $((C_{\mathscr{A}})_M\setminus \{m_1\})\supsetneq ((C_{\mathscr{A}})_M\setminus \{m_1,m_2\})\supsetneq\ldots$, there must be a least integer $k$, $1\leq k\leq |C|-1$, such that $d_k$ is not adjacent to any vertex in $((C_{\mathscr{A}})_M\setminus \{m_1,\ldots,m_k\})$. If $k=|C|-1$, $C_{\mathscr{A}}$ contains pure basic-type subgraph $d_{|C|-1}m_{|C|-1}d_{|C|-2}m_{|C|-2}\ldots d_1m_1d_{|C|}$, since $G_{\mathscr{A}}$ is acyclic. And it is actually equal to $C_{\mathscr{A}}$. 

    If $k<|C|-1$, then again since $\mathscr{A}$ is acyclic, the subgraph induced by \\$\{d_{|C|},d_1,\ldots ,d_k,m_1,\ldots,m_k\}$ is of pure basic-type, namely, $d_km_k\ldots d_1m_1d_{|C|}$. Now, investigating $C_{\mathscr{A}}$, we have two cases: 
    \\Case 1. Some element in $C\setminus \{d_k,\ldots,d_1,d_{|C|}\}$ is adjacent to some element of $\{m_k,\ldots,m_1\}$. Or \\Case 2. Some element in $(C_{\mathscr{A}})_M\setminus \{m_k,\ldots,m_1\}$ is adjacent to some element in $\{d_k,\ldots,d_1,d_{|C|}\}$.
    \\In case 1, suppose $d_{k+1}$ is adjacent to some $m_i$, $1\leq i\leq k$. But by Lemma \ref{lemma:m can't ad. to exa. one n}, $m_{k+1}$ cannot be adjacent only to $d_{k+1}$, and $m_{k+1}$ cannot be adjacent to any vertex of $\{d_k,\ldots,d_1,d_{|C|}\}$ because $G_{\mathscr{A}}$ is acyclic. Assume $m_{k+1}$ is adjacent to, say, $d_{k+2}$, again $m_{k+2}$ must be adjacent to some vertex in $\{d_{k+3},\ldots,d_{|C|-1}\}$. As $|C|-1$ is a finite integer, we finally have that $m_{|C|-1}$ must be adjacent to some vertex of $C\setminus \{d_{|C|-1},d_{|C|-2},\ldots,d_1,d_{|C|}\}=\emptyset$. This is not possible. 
    \\So we must proceed with Case 2. Some $d_i$, $1\leq i<k$ is adjacent to, again, without loss of generality, $m_{k+1}$. Then we ask whether $d_{k+1}$ is adjacent to some vertex in $\{m_{k+2},...,m_{|C|-1}\}$. Again, there must be a least positive integer $k_1$ such that $d_{k+k_1}$ is not adjacent to any vertex in $\{m_{k+k_1+1},\ldots,m_{|C|-1}\}$, which is empty if $k+k_1=|C|-1$. Then again, the graph induced by $\{d_i,d_{k+1},\ldots,d_{k+k_1},m_{k+1},\ldots,m_{k+k_1}\}$ is a basic-type subgraph that can be read as $d_{k+k_1}m_{k+k_1}d_{k+k_1-1}m_{k+k_1-1}\ldots d_{k+1}m_{k+1}d_i$. So far we obtain a planting of a pure basic-type graph $d_{k+k_1}m_{k+k_1}\ldots d_{k+1}m_{k+1}d_i$ from $d_i$ to $d_i$ of a pure basic-type graph $d_km_k\ldots d_1m_1d_{|C|}$, with $d_i$ the joint of planting. 
    
    Then again $m_{k+k_1+1}$ must be adjacent to some element of \\$\{d_{k+k_1},d_{k+k_1-1},\ldots,d_k,\ldots,d_1,d_{|C|}\}$, and we continue with the planting. But each time we start a new planting, we begin with an element of $(C_{\mathscr{A}})_M$ that is not yet in the planting we previously had. And at each planting, the cardinality of the subset of $(C_{\mathscr{A}})_M$ that is already in the planting strictly increases. Since $|(C_{\mathscr{A}})_M|$ is finite, the planting has to stop.
\end{proof}

\begin{figure}[h]
    \centering
    \includegraphics[width=0.2\linewidth]{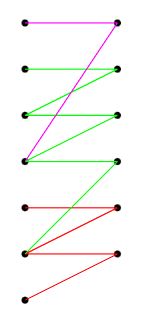}
    \caption{Plantings of pure basic-type graphs.}
    \label{fig:plantings of pure basic type}
\end{figure}

Figure \ref{fig:plantings of pure basic type} illustrates an example of a circuit described in the preceding Theorem.

\section{Circuits intersecting in a basic way}
In this section, we are going to see what happens if circuits in a transversal matroid intersect in a `nice' way.

\begin{definition}\label{def:intersect basic way}
    Two circuits $C_1$ and $C_2$ of a transversal matroid $M$ with nonempty intersection are said to intersect in a basic way in a presentation $\mathscr{A}$ if they have basic-type subgraph of the form $n_1m_1\ldots n_im_i\ldots n_l$ and $n_1m_1\ldots n_im_in'_{i+1}m'_{i+1}\ldots n'_{j}$ respectively, such that $\{n'_{i+1},\ldots,n'_{j}\}\cap \{n_1,\ldots,n_i,\ldots,n_l\}=\emptyset$, and $\{m'_{i+1},\ldots,m'_{j-1}\}\cap \{m_1,\ldots,m_{l-1}\}=\emptyset$, and $C_1=\{n_1,\ldots,n_i,\ldots,n_l\},C_2=\{n_1,\ldots,n_i,n_{i+1}',\ldots,n_j'\}$. Note that this does not say that $(C_1)_{\mathscr{A}}$ or $(C_2)_{\mathscr{A}}$ are of pure basic-type in $G_{\mathscr{A}}$.
\end{definition}

\begin{example}
Figure $\ref{fig:inter. basic way}$ illustrates the situation stated in the preceding definition. $C_1=\{n_1,n_2,n_3,n_4\}$ and $C_2=\{n_1,n_2,n_5,n_6\}$ intersect in a basic way in this presentation. Note that this bipartite graph is not acyclic and $C_1,C_2$ are also of basic-type in this presentation. 

In Figure $\ref{fig:inter. non basic way}$, $C_1=\{n_1,n_2,n_3,n_4\}, C_2=\{n_1,n_2,n_5,n_6\}$ do not intersect in a basic way in this presentation.

    \begin{figure}[h]
    \begin{subfigure}[T]{0.5\textwidth}
\includegraphics[width=0.62\linewidth, height=6.8cm]{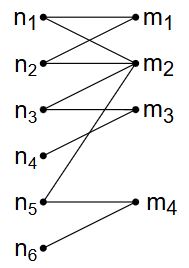} 
\caption{}
\label{fig:inter. basic way}
\end{subfigure}
\begin{subfigure}[T]{0.5\textwidth}
\includegraphics[width=0.65\linewidth, height=7.5cm]{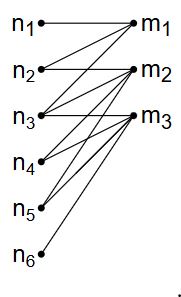}
\caption{}
\label{fig:inter. non basic way}
\end{subfigure}
\caption{Illustration of circuits intersecting in a basic way and a non-basic way.}
\label{fig:both inter. basic and not}
\end{figure}
\end{example}

\begin{remark}
Note that if $C_1,C_2\in \mathscr{C}(M)$ intersect in a basic way in a presentation $\mathscr{A}$ of a transversal matroid $M$, by Definition \ref{def:intersect basic way}, $(C_1)_{\mathscr{A}},(C_2)_{\mathscr{A}}$ contain basic type subgraphs $n_1m_1\ldots n_im_i\ldots n_k$ and $n_1m_1\ldots n_im_in'_{i+1}m'_{i+1}\ldots n'_j$, respectively, where $\{n'_{i+1},\ldots,n'_j\}\cap C_1=\emptyset$ and $\{m'_{i+1},\ldots,m'_{j-1}\}\cap ((C_1)_{\mathscr{A}})_M=\emptyset$. But if there is a $C_3\in \mathscr{C}(M)$ such that $C_1,C_3$ also intersect in a basic way in $\mathscr{A}$, that does not mean $(C_3)_{\mathscr{A}}$ must contain a basic-type subgraph of the form $n_1m_1\ldots n_rm_rn^3_{r+1}m^3_{r+1}\ldots n^3_s$ since $(C_1)_{\mathscr{A}}$ may contain other basic-type subgraphs. But there must be $\sigma \in S_k,\tau\in S_{k-1}$ so that $(C_1)_{\mathscr{A}},(C_3)_{\mathscr{A}}$ contain basic type subgraphs $n_{\sigma(1)}m_{\tau(1)}\ldots n_{\sigma(t)}m_{\tau(t)}\ldots n_{\sigma(k)}$ and $n_{\sigma(1)}m_{\tau(1)}\ldots n_{\sigma(t)}m_{\tau(t)}n^3_{t+1}m^3_{t+1}\ldots n^3_w$, where $\{n^3_{t+1},\ldots,n^3_w\}\cap C_1=\emptyset$ and $\{m^3_{t+1},\ldots,m^3_{w-1}\}\cap ((C_1)_{\mathscr{A}})_M=\emptyset$. This makes our analysis much more complicated. However, under certain assumptions, the situation can be made friendlier to deal with. We are going to return to this issue later. 
\end{remark}

\begin{figure}[h]
    \centering
    \includegraphics[width=0.35\linewidth]{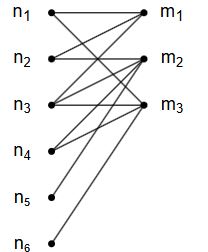}
    \caption{Various ways to read a basic-type subgraph}
    \label{fig:inter. basic way variance}
\end{figure}

\begin{example}
    Figure \ref{fig:inter. basic way variance} illustrates an example described by the preceding remark. The bipartite graph in the figure defines a presentation $\mathscr{A}$ for a  transversal matroid $M$ on the ground set $\{n_1,n_2,n_3,n_4,n_5,n_6\}$. $C_1=\{n_1,n_2,n_3,n_4\},C_2=\{n_1,n_2,n_5\}$ are both circuits of $M$ and they clearly intersect in a basic way in $\mathscr{A}$. The subgraphs $(C_1)_{\mathscr{A}}$ and $(C_2)_{\mathscr{A}}$ contain basic-type subgraphs $n_1m_1n_2m_2n_3m_3n_4$ and $n_1m_1n_2m_2n_5$ respectively. Note that $C_3=\{n_2,n_3,n_4,n_6\}\in \mathscr{C}(M)$ and it also intersects $C_1$ in a basic way in $\mathscr{A}$. This is because $C_1$ also contains basic-type subgraph $n_2m_1n_3m_2n_4m_3n_1$, and $C_3$ contains basic-type subgraph $n_2m_1n_3m_2n_4m_3n_6$ in $\mathscr{A}$.
\end{example}

\begin{lemma}\label{lemma: S_A basic type S is a cir.}
    If $\mathscr{A}$ is a presentation for a transversal matroid $M$ and $S\subseteq E(M)$, and if $S_{\mathscr{A}}$ is a pure basic-type subgraph or a basic-type subgraph in $G_{\mathscr{A}}$, then $S\in \mathscr{C}(M)$.
\end{lemma}
\begin{proof}
    Suppose $S_{\mathscr{A}}$ is of pure basic-type. We may assume that $S=\{s_1,\ldots,s_n\}$, and $S_{\mathscr{A}}$ can be read as $s_1m_1\ldots m_{n-1}s_n$ so that two vertices in $S_{\mathscr{A}}$ are adjacent if and only if they are adjacent in the word $s_1m_1\ldots m_{n-1}s_n$.$S$ cannot be matched since $|(S_{\mathscr{A}})_M|=n-1=|S|-1$, so $S$ is dependent. Our goal is to show that any subset of $S$ of cardinality $|S|-1$ can be matched. 
    \\ \noindent Case 1. $S'=S\setminus \{s_n\}$  
    \\ $S'$ can clearly be matched, with $s_i$ matched to $m_i$, for all $1\leq i\leq n-1$. 
    \\ Case 2. $S'=S\setminus \{s_1\}$ \\$S'$ can be matched as follows: $s_j$ is matched to $m_{j-1}$, for all $2\leq j\leq n$.
    \\Case 3. $S'=S\setminus\{s_k\}$ for some $k$, $1<k<n$
    \\ $S'$ can be matched as follows: $s_i$ is matched to $m_i$ for all $1\leq i<k$, and $s_j$ is matched to $m_{j-1}$ for all $k<j\leq n$. 
    \\Therefore, $S\in \mathscr{C}(M)$. Furthermore, if $S_{\mathscr{A}}$ is a basic-type subgraph, we may assume $S_{\mathscr{A}}$ is of the above form but there may be extra edges between non-adjacent vertices in the word $s_1m_1\ldots m_{n-1}s_n$, apart from the edges between adjacent vertices in $s_1m_1\ldots m_{n-1}s_n$. Then any subset of $S$ of cardinality $|S|-1$ can still be matched since adding edges does not affect a matching we already had. And $S$ still cannot be matched because $|(S_{\mathscr{A}})_M|=n-1=|S|-1$. Therefore, $S$ is still a circuit. 
\end{proof}

\begin{proposition}\label{prop:connect to one commoin part}
    Suppose two circuits $C_1,C_2$ in a transversal matroid $M$ intersect in a basic way in a given presentation $\mathscr{A}$ for $M$. Then $(C_1\cup C_2)\setminus (C_1\cap C_2)\in \mathscr{C}(M)$ if and only if $(C_1\cup C_2)\setminus (C_1\cap C_2)$ is adjacent to only one vertex in $((C_1\cap C_2)_{\mathscr{A}})_M$ in $G_{\mathscr{A}}$.
   \end{proposition}
\begin{proof}
    By our assumption, we may assume $C_1=\{n_1,\ldots,n_i,\ldots,n_k\}$ and $C_2=\{n_1,\ldots,n_i,n_{i+1}',\ldots,n_j'\}$ have basic-type subgraphs of the form $n_1m_1\ldots n_im_i\ldots n_k$ and $n_1m_1\ldots n_im_in'_{i+1}m'_{i+1}\ldots n'_j$, respectively, in $G_{\mathscr{A}}.$ And $(C_1\cup C_2)\setminus(C_1\cap C_2)=\{n_{i+1},\ldots,n_k,n'_{i+1},\ldots,n'_j\}$, $((C_1\cap C_2)_{\mathscr{A}})_M\supseteq \{m_1,\ldots,m_i\}$. 
    \\Then $((C_1\cup C_2)\setminus (C_1\cap C_2))_{\mathscr{A}}$ in $G_{\mathscr{A}}$ contains a basic-type subgraph 
    $$n_km_{k-1}n_{k-1}\ldots n_{i+1}m_in'_{i+1}m'_{i+1}\ldots m'_{j-1}n'_j.$$ Note that $(C_1\cup C_2)\setminus (C_1\cap C_2)$ is already adjacent to one vertex in $((C_1\cap C_2)_{\mathscr{A}})_M$, namely, $m_i$. 

    For the only if part, if any vertex in $(C_1\cup C_2)\setminus (C_1\cap C_2)$ is adjacent to any vertex in $((C_1\cap C_2)_{\mathscr{A}})_M\setminus \{m_i\}$, then $|(((C_1\cup C_2)\setminus (C_1\cap C_2))_{\mathscr{A}})_M|\geq |(C_1\cup C_2)\setminus (C_1\cap C_2)|$, which contradicts Proposition \ref{prop: cir. N=M+1}. 
    
    For the converse, $((C_1\cup C_2)\setminus(C_1\cap C_2))_{\mathscr{A}}$ contains basic type subgraph $n_km_{k-1}n_{k-1}\ldots n_{i+1}m_in'_{i+1}m'_{i+1}\ldots m'_{j-1}n'_j$ and $(C_1\cup C_2)\setminus (C_1\cap C_2)=\{n_{i+1},\ldots,n_k,n'_{i+1},\ldots,n'_j\}$. The assumption $(C_1\cup C_2)\setminus (C_1\cap C_2)$ contains no element from $((C_1\cap C_2)_{\mathscr{A}})_M\setminus \{m_i\}$ means that $((C_1\cup C_2)\setminus (C_1\cap C_2))_{\mathscr{A}}$ is a basic-type subgraph. By Lemma \ref{lemma: S_A basic type S is a cir.}, $((C_1\cup C_2)\setminus (C_1\cap C_2))\in \mathscr{C}(M)$.
\end{proof}

Now we address the issue raised by the first remark in this chapter.
\begin{proposition}\label{prop:rearranging inter. basic type}
    In a transversal matroid $M$, if $C_1,C_2\in \mathscr{C}(M)$, $(C_1\cup C_2)\setminus (C_1\cap C_2)\in \mathscr{C}(M)$, with $C_1=\{n_1,\ldots,n_k\}$ and $(C_1)_{\mathscr{A}}$ contains a basic-type subgraph $n_1m_1\ldots n_k$ in $G_{\mathscr{A}}$ for a presentation $\mathscr{A}$ of $M$, and also if $C_1,C_2$ intersect in a basic way in $\mathscr{A}$, so that $(C_1)_{\mathscr{A}},(C_2)_{\mathscr{A}}$ contain basic type subgraphs $n_{\sigma(1)}m_{\tau(1)}\ldots n_{\sigma(t)}m_{\tau(t)}\ldots n_{\sigma(k)}$ and $n_{\sigma(1)}m_{\tau(1)}\ldots n_{\sigma(t)}m_{\tau(t)}n'_{t+1}m'_{t+1}\ldots n'_{w}$ in $G_{\mathscr{A}}$ such that $\{n'_{t+1},\ldots,n'_w\}\cap C_1=\emptyset$ and $\{m'_{t+1},\ldots,m'_{w-1}\}\cap ((C_1)_{\mathscr{A}})_M=\emptyset$,$\sigma\in S_k,\tau \in S_{k-1}$. Then $(C_2)_{\mathscr{A}}$ contains either $n_1m_1\ldots n_tm_tn'_{t+1}m'_{t+1}\ldots n'_w$ or $n_km_{k-1}n_{k-1}\ldots m_{k-t}n'_{t+1}m'_{t+1}\ldots n'_w$ as a basic-type subgraph in $G_{\mathscr{A}}$.
\end{proposition}
\begin{proof}
    Let $C_1,C_2\in \mathscr{C}(M)$ satisfy our assumptions. Suppose $m_{\tau(t)}=m_h$ for some $1\leq h\leq k-1$, we claim that either $n_h$ or $n_{h+1}$ belongs to $C_2$. Suppose neither of them belongs to $C_2$, then note that for each $n_{h'}\in C_2,1\leq h'\leq k$, by Remark \ref{re:(G)_M-(S)_M NO S}, $m_{h'}\in (C_2)_{\mathscr{A}}$, and $m_{h'-1}\in (C_2)_{\mathscr{A}}$ if $h'\neq 1$. Then $m_{h'}$ or $m_{h'-1}$ is not $m_h$. That means $|C_1\cap C_2|=t$ but $|((C_1)_{\mathscr{A}})_M\cap ((C_2)_{\mathscr{A}})_M|> t$, which is a contradiction (Note that we assumed that $C_1$ and $C_2$ intersect in a basic way and $((C_1\cup C_2)\setminus (C_1\cap C_2))\in \mathscr{C}(M)$. As a result of Proposition \ref{prop:connect to one commoin part}, $|C_1\cap C_2|=|((C_1)_{\mathscr{A}})_M\cap ((C_2)_{\mathscr{A}})_M|=t$). Now suppose both $n_h,n_{h+1}$ belong to $C_2$. Then $m_{h-1}$ must be in $(C_2)_{\mathscr{A}}$ if $h\neq 1$ and $m_{h+1}$ must be in $(C_2)_{\mathscr{A}}$ if $h\neq k$. Then we apply Proposition \ref{prop:connect to one commoin part}, $n_{h-1}$ must be in $C_2$, and $n_{h+2}\in C_2$ if $h\neq k-1$. Continuing this process, the basic-type subgraph $n_1m_1\ldots n_k$ of $C_1$ becomes a subgraph of $(C_2)_{\mathscr{A}}$, which means $C_1$ is a proper subset of $C_2$, that is absurd. Therefore, our claim is justified. 

    Now, if $n_h\in C_2$, then $m_{h-1}\in ((C_2)_{\mathscr{A}})_M$. Again, by the preceding Proposition \ref{prop:connect to one commoin part}, $n_{h-1}\in C_2$ if $h\neq 1$, continuing this process until we reach $n_1\in C_2$. Therefore, $(C_2)_{\mathscr{A}}$ contains a basic type subgraph $n_1m_1\ldots n_hm_hn'_{t+1}m'_{t+1}\ldots n'_w$. Note that after we reach $n_1$, the process cannot continue to reach $m_{h+i},h+i\leq k$, if this happens, arguments similar to those above show that the basic type subgraph $n_1m_1...n_k$ of $C_1$ becomes a subgraph of $(C_2)_{\mathscr{A}}$, a contradiction. Now, recall that $|C_1\cap C_2|=|((C_1)_{\mathscr{A}})_M\cap ((C_2)_{\mathscr{A}})_M|=t$, so $h=t$. The case involving $n_{h+1}$ is analogous.
\end{proof}

\begin{theorem}\label{thm:admits acyclic pre.}
    If $\mathscr{A}$ is a presentation for a transversal matroid $M$ such that there are $C_1,C_2\in \mathscr{C}(M)$ intersecting in a basic way in $\mathscr{A}$ and $(C_1\cup C_2)\setminus (C_1\cap C_2)\in \mathscr{C}(M)$, moreover, for all $C\in \mathscr{C}(M)$ with $C\cap C_i\neq \emptyset$, we have that $C$ intersects $C_i$ in a basic way, for $i\in \{1,2\}$, then there is a presentation $\mathscr{A'}$ for $M$ such that $(C_1\cup C_2)_{\mathscr{A'}}$ is acyclic.
\end{theorem}
\begin{proof}
         Let $C_1,C_2\in \mathscr{C}(M)$ satisfy our assumption. Suppose $C_1=\{n_1,\ldots,n_j,\ldots,n_k\}$ and $C_2=\{n_1,\ldots,n_j,n'_{j+1},\ldots,n'_h\}$, and $(C_1)_{\mathscr{A}}$ and $(C_2)_{\mathscr{A}}$ contain basic type subgraphs
    \begin{equation}\label{eq:C_1 basic type word}
    n_1m_1\ldots n_jm_jn_{j+1}m_{j+1}\ldots m_{k-1}n_k
    \end{equation}
    and  
    \begin{equation}\label{eq:C_2 basic type word}
    n_1m_1\ldots n_jm_jn'_{j+1}m'_{j+1}\ldots m'_{h-1}n'_h
    \end{equation}
    respectively, with $\{m_{j+1},\ldots,m_{k-1}\}\cap \{m'_{j+1},\ldots,m'_{h-1}\}=\emptyset$. Note that since \eqref{eq:C_1 basic type word} and \eqref{eq:C_2 basic type word} are basic-type subgraphs contained in $(C_1)_{\mathscr{A}}$ and $(C_2)_{\mathscr{A}}$, it is possible for non-adjacent vertices in \eqref{eq:C_1 basic type word} or \eqref{eq:C_2 basic type word} to be adjacent in $(C_1)_{\mathscr{A}}$ or $(C_2)_{\mathscr{A}}$. We call those edges joining non-adjacent vertices in \eqref{eq:C_1 basic type word} or \eqref{eq:C_2 basic type word} the set of extra edges.   
    
    Figure \ref{fig:C_1,C_2 basic type surplus} is an illustration of the situation we describe, where the edges marked in red are extra edges. The edges and vertices marked in bold indicate the intersection of $(C_1)_{\mathscr{A}}$ and $(C_2)_{\mathscr{A}}$. Remove all extra edges and denote the resulting presentation as $\mathscr{A'}$. Then $(C_1)_{\mathscr{A'}}$ and $(C_2)_{\mathscr{A}'}$ become pure basic-type subgraphs and by Lemma \ref{lemma: S_A basic type S is a cir.}, $C_1$ and $C_2$ remain as circuits. Figure \ref{fig:C_1,C_2 basic type} illustrates this situation. $(C_1\cup C_2)_{\mathscr{A'}}$ is acyclic because both $(C_1)_{\mathscr{A'}}$ and $(C_2)_{\mathscr{A'}}$ are acyclic, and by Proposition \ref{prop:connect to one commoin part}, $m_j$ is the only element in $((C_1\cap C_2)_{\mathscr{A}})_M$ that can be adjacent to an element in $C_2\setminus (C_1\cap C_2)$. Because $(C_2)_{\mathscr{A'}}$ is pure basic-type, $n'_{j+1}$ is the only element of $C_2$ that can be adjacent to $m_j$. So $(n'_{j+1},m_j)$ is the only edge between $C_2\setminus (C_1\cap C_2)$ and $((C_1\cap C_2)_{\mathscr{A'}})_M$. Also recall that from Remark \ref{re:(G)_M-(S)_M NO S}, no element in $C_2\setminus (C_1\cap C_2)$ can be adjacent to any element in $((C_1)_{\mathscr{A'}})_M\setminus ((C_1\cap C_2)_{\mathscr{A'}})_M$). Therefore, $(C_1\cup C_2)_{\mathscr{A'}}$ is acyclic. 
        
\begin{figure}[h]
\begin{subfigure}{0.5\textwidth}
\includegraphics[width=0.6\linewidth, height=11.2cm]{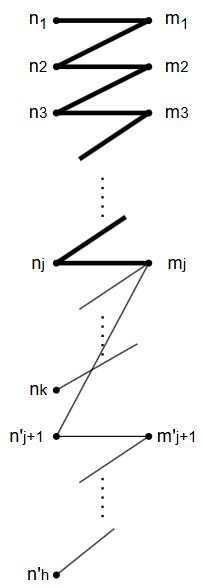} 
\caption{}
\label{fig:C_1,C_2 basic type}
\end{subfigure}
\begin{subfigure}{0.5\textwidth}
\includegraphics[width=0.6\linewidth, height=11.2cm]{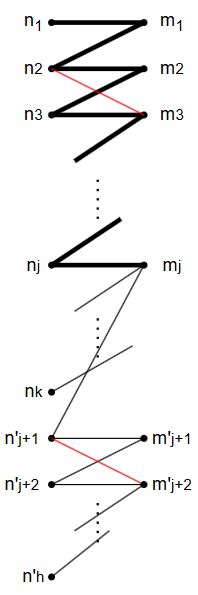}
\caption{}
\label{fig:C_1,C_2 basic type surplus}
\end{subfigure}

\caption{Illustration of $(C_1\cup C_2)_{\mathscr{A}}$}
\label{fig:C_1,C_2}
\end{figure}

Now we need to check that removing extra edges does not affect other circuits of $M$. Suppose we have $C_3\in \mathscr{C}(M)$ and $C_3\cap (C_1\cup C_2)=\emptyset$, then $(C_3)_{\mathscr{A'}}$ does not contain any of the extra edges, since every extra edge is incident to one element in $C_1$ or $C_2$. If $C_3\cap C_1\neq \emptyset$, since $C_3$ and $C_1$ intersect in a basic way in $\mathscr{A}$, then by Proposition \ref{prop:rearranging inter. basic type},  $(C_3)_{\mathscr{A}}$ contains either 
$$n_1m_1\ldots n_im_i\ldots n'_{i+1}m'_{i+1}\ldots m'_{r-1}n'_r$$
as a basic type subgraph with $C_3=\{n_1,\ldots,n_i,n'_{i+1},\ldots,n'_r\}$ or 
$$n_km_{k-1}n_{k-1}\ldots m_{k-k'}n_{k-k'+1}'m_{k-k'+1}\ldots n'_{q}$$ 
as a basic-type subgraph with $C_3=\{n_k,n_{k-1},n'_{k-k'+1},\ldots,n'_q\}$. The intersection of the basic-type subgraph $(C_3)_{\mathscr{A}}$ with word \eqref{eq:C_1 basic type word} consists of consecutive vertices in word \eqref{eq:C_1 basic type word}, thus again by Lemma \ref{lemma: S_A basic type S is a cir.}, $C_3$ remains a circuit after removing extra edges. 

But we also need to make sure that removing extra edges does not create new circuits. We are going to show this by contradiction. Suppose removing an extra edge $e=\{n_f,m_{g_1}\}$ in $(C_1)_{\mathscr{A}}$ creates a new circuit. Let the presentation we get from removing $e$ be $\mathscr{A}_0$. It means that there exists a $D\subseteq E(M)$ such that $D$ is independent in $M=T(G_{\mathscr{A}})$, but $D$ is a circuit of $T(G_{\mathscr{A}_0})$. Note that $n_f\in D$, otherwise removing $e$ does not affect the matching of $D$ and $D$ would stay as a circuit in $T(G_{\mathscr{A}_0})$. Also note that $m_{g_1}\notin (D_{\mathscr{A}_0})_M$, since $D$ is independent in $T(G_{\mathscr{A}})$ but $D$ becomes a circuit in $T(G_{\mathscr{A}_0})$, so $(D_{\mathscr{A}_0})_M\subsetneq (D_{\mathscr{A}})_M$. The only possible element lost from $(D_{\mathscr{A}})_M$ to $(D_{\mathscr{A}_0})_M$ is $m_g$. Moreover, by Proposition \ref{prop: cir. N=M+1}, $|D|=|(D_{\mathscr{A}_0})_M|+1$, therefore, $|D|=|(D_{\mathscr{A}})_M|$. Observe that whenever $m_{g_i}\notin (D_{\mathscr{A}_0})_M$, $n_{g_i}\notin D$ for all $1\leq g_i\leq k$. Now we want to show that there is a $D'\in \mathscr{C}(T(G_{\mathscr{A}}))$ with $D\subsetneq D'$.

Take $D\cup \{n_{g_1}\}$, if $n_{g_1}$ is not adjacent to any element in $((C_1)_{\mathscr{A}})_M\setminus ((D_{\mathscr{A}})_M\cup \{m_{g_1}\})$ in $G_{\mathscr{A}}$, we let $D'=D\cup \{n_{g_1}\}$. If $n_{g_1}$ is adjacent to some element, say $m_{g_2} \in ((C_1)_{\mathscr{A}})_M\setminus ((D_{\mathscr{A}})_M\cup \{m_{g_1}\})$, then take $D\cup \{n_{g_1},n_{g_2}\}$. Again, we ask whether there is an element in $\{n_{g_1},n_{g_2}\}$ that is adjacent to some $m_{g_3}\in ((C_1)_{\mathscr{A}})_M\setminus ((D_{\mathscr{A}})_M\cup \{m_{g_1},m_{g_2}\})?$ If the answer is no, let $D'=D\cup \{n_{g_1},n_{g_2}\}$. If the answer is yes, take $D\cup \{n_{g_1},n_{g_2},n_{g_3}\}$ and continue the process. Eventually, we can let $D'=D\cup \{n_{g_1},\ldots,n_{g_z}\}$ such that no element in $\{n_{g_1},\ldots,n_{g_z}\}$ is adjacent to any element in $((C_1)_{\mathscr{A}})_M\setminus ((D_{\mathscr{A}})_M\cup \{m_{g_1},\ldots,m_{g_z}\})$. We want to show that $D'\in \mathscr{C}(T(G_{\mathscr{A}}))$, and we argue that by induction. Our induction hypothesis is that whenever $D'=D\cup \{n_{g_1},\ldots,n_{g_{t-1}}\}$ then $D'\in \mathscr{C}(T(G_{\mathscr{A}}))$, with $1\leq t-1<z$. Now suppose $D'=D\cup \{n_{g_1},\ldots,n_{g_t}\}$, by the way we defined $D'$, we have $|D'|=|(D'_{\mathscr{A}})_M|+1$, thus $D'$ is dependent in $T(G_{\mathscr{A}})$. Then we need to show that any subset $S\subseteq D'$ of cardinality $|D'|-1$ can be matched in $G_{\mathscr{A}}$. 

\noindent Case 1) If $n_{g_t}\notin S$, then $S=D\cup \{n_{g_1},\ldots,n_{g_{t-1}}\}$.
\\ Imagine that we prohibit the use of all edges connecting elements in $S$ to $m_{g_t}$ in any matching of $S$, by our induction hypothesis, any subset of $S$ of cardinality $|D\cup \{n_{g_1},\ldots,n_{g_{t-1}}\}|-1$ can be matched in $G_{\mathscr{A}}$. Also note that by the way we construct $D'$, the matching of that subset does not involve $m_{g_t}$. But with those edges allowed (in particular those having $m_{g_t}$ as an endpoint) in the matching of $S$, $S$ can be matched in $G_{\mathscr{A}}$.

\noindent Case 2) $n_{g_t}\in S$, this means $n_{g_i}\notin S$ for some $n_{g_i}\in \{n_{g_1},\ldots,n_{g_{t-1}}\}$.
\\$n_{g_t}$ can be matched by $m_{g_t}$. Then we are not using any edge connecting $m_{g_t}$ to elements in $S\setminus \{n_{g_t}\}$ in the matching of $S\setminus \{n_{g_t}\}$. Then applying our induction hypothesis, $S\setminus \{n_{g_t}\}$ still can be matched. 

Now we are left to show the base case. That is, if $D'=D\cup \{n_{g_1}\}$. Again, we need to show that any subset $S\subseteq D'$ of cardinality $|D'|-1$ can be matched in $G_{\mathscr{A}}$.

\noindent Case 1) $n_{g_1}\notin S$, this means $S=D$.
\\ Certainly, $S$ can be matched in $G_{\mathscr{A}}$ because $D$ is assumed to be independent in $T(G_{\mathscr{A}})$.

\noindent Case 2) $n_{g_1}\in S$, this means $d\notin S$ for one $d\in D$. 
Since $D$ is a circuit in $T(G_{\mathscr{A}_0})$, any subset of $D$ of cardinality $|D|-1$ can be matched by some subset of $(D_{\mathscr{A}_0})_M$ in $G_{\mathscr{A}_0}$. Recall that $m_{g_1}\notin (D_{\mathscr{A}_0})_M$, thus it can be matched by some subset of $(D_{\mathscr{A}})_M$. Together with the fact $n_{g_1}$ can be matched by $m_{g_1}$, we have that $S$ can be matched in $G_{\mathscr{A}}$.

Therefore, we have $D'\in \mathscr{C}(T(G_{\mathscr{A}}))$ with $D\subsetneq D'$. Remember that erasing extra edges does not affect any circuit of $T(G_{\mathscr{A}})$, so $D'$ stays as a circuit in $T(G_{\mathscr{A}_0})$. But $D$, as a proper subset of $D'$, becomes a circuit in $T(G_{\mathscr{A}_0})$, which is a contradiction.
\end{proof}

\section{All presentations of a transversal matroid}\label{sec: Brualdi et al.}
This section is a survey of an algorithm proposed by Brualdi et al. \cite{Brualdi*}, who used the algorithm to determine whether or not a matroid is a transversal matroid. If the matroid is transversal, the algorithm yields an explicit determination of the unique maximal presentation. Furthermore, they proposed a way of obtaining all presentations of a transversal matroid. 

Let $M$ be a matroid on a finite set $E$ and $S\subseteq E$. $S$ is said to be coloop-free if $M|S$ has no coloop. The Algorithm \cite{Brualdi*}: 

\noindent Given a matroid $M$ of rank $r$ on a finite set E, 
\\(1) Let $\mathscr{F}_M^{(1)}$ consist of the coloop-free hyperplanes each with multiplicity 1. 
\\(2) Construct $\mathscr{F}_M^{(2)}$ which is to include $\mathscr{F}_M^{(1)}$ as a subfamily. Let $F$ be a coloop-free $r-2$ flat contained in $p$ coloop-free hyperplanes. If $p\geq 2$, $F$ is not a member of $\mathscr{F}_M^{(2)}$; otherwise $F$ is a member of $\mathscr{F}_M^{(2)}$ with multiplicity $2-p$

$$\vdots$$

\noindent (k) Construct a family $\mathscr{F}_M^{(k)}$ which is to include $\mathscr{F}_M^{(k-1)}$ as a subfamily. Let $F$ be a coloop-free $(r-k)$-flat, suppose $F$ is contained in $p$ members of $\mathscr{F}_M^{(k-1)}$. If $p\geq k$, $F$ is not a member of $\mathscr{F}_M^{(k)}$; otherwise $F$ is a member of $\mathscr{F}^{(k)}$ with multiplicity $k-p.$
$$\vdots$$

\noindent (r) Construct a family $\mathscr{F}_M^{(r)}$ which is to include $\mathscr{F}_M^{(r-1)}$ as a subfamily. Suppose $\mathscr{F}_M^{(r-1)}$ has $p$ members, if $p\geq r$, then $\emptyset$ is not a member of $\mathscr{F}_M^{(r)}$; otherwise $\emptyset$ is a member of $\mathscr{F}_M^{(r)}$ with multiplicity $r-p$. 

\noindent Let $\mathscr{F}_M=\mathscr{F}_M^{(r)}.$$\mathscr{F}_M$ is a family of coloop-free flats with $|\mathscr{F}_M|\geq r.$

If $M$ is a transversal matroid of rank $r$ with presentation $\mathscr{A}=\{A_1,\ldots, A_r\}$, then $\mathscr{A}$ is a \textit{maximal presentation} of $M$ if 
$$(A'_1,\ldots,A'_r), A_i\subseteq A'_i (1\leq i\leq r) \mbox{ is a presentation of } M$$
implies $A'_i=A_i (1\leq i\leq r).$

\begin{theorem}[{\cite{Brualdi*}}]
    If $M$ is a transversal matroid of rank $r$, then the maximal presentation is unique. If $(F_1,\ldots,F_n)$ is the family obtained in the algorithm, then $n=r$ and $\{E(M)\setminus F_1,\ldots, E(M)\setminus F_r\}$ is the maximal presentation of $M$.
\end{theorem}

Once they figure out the maximal presentation, they may obtain all other presentations:

\begin{theorem}[{\cite{Brualdi*}}]
    Let $M$ be a transversal matroid of rank $r$ with maximal presentation $\{M_1,\ldots, M_r\}$. Let $\mathscr{A}=\{A_1,\ldots, A_r\}$ be a family of subsets of $E(M)$. Then $\mathscr{A}$ is a presentation of $M$ if and only if both of the following hold:
    \begin{itemize}
        \item there exists a permutation $\sigma$ of $[r]$ such that $A_k\subseteq M_{\sigma(k)}, (1\leq k\leq r),$
        \item $\rho\left(\bigcap_{i\in I}(E(M)\setminus A_i)\right)\leq r-|I|, (I\subseteq [r]).$
    \end{itemize}
\end{theorem}

Lastly, we finish by reflecting on what we established in the previous section. Note that due to Proposition \ref{prop:connect to one commoin part}, if we construct a transversal matroid such that all of its circuits are of pure basic type in a presentation, then its number of circuits is easy to control. For example,
\begin{figure}[h]
    \centering
    \includegraphics[width=0.4\linewidth]{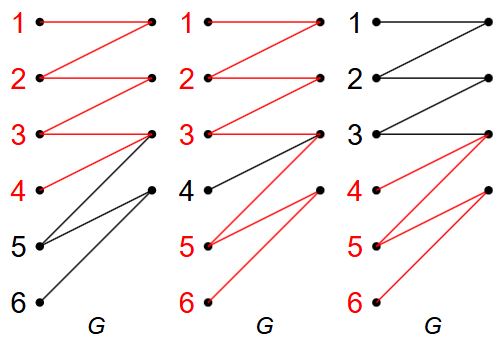}
    \caption{All circuits of $T(G)$ are of pure basic type in this presentation}
    \label{fig:pure basic type easy to control}
\end{figure}
\noindent in Figure \ref{fig:pure basic type easy to control}, we draw circuits $\{1,2,3,4\}$ and $\{1,2,3,5,6\}$ of pure basic type. Then $\{4,5,6\}$ is also a circuit of pure basic type. Moreover, we can make use of the Brualdi-Dinolt's Algorithm to construct a transversal matroid such that its dual is not transversal:
\begin{figure}[h]
    \centering
    \includegraphics[width=0.8\linewidth]{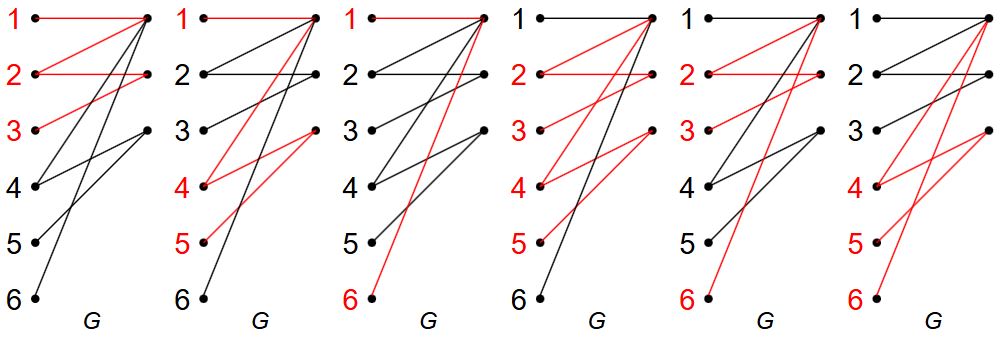}
    \caption{A transversal matroid of rank 3 such that its dual is non-transversal.}
    \label{fig:T(G) dual is rank 3 non tran}
\end{figure}

\noindent In Figure \ref{fig:T(G) dual is rank 3 non tran}, there are six circuits (colored in red) in $T(G)$, $T(G)^*$ is of rank 3 and has six hyperplanes:$$\{4,5,6\},\{2,3,6\},\{2,3,4,5\},\{1,6\},\{1,4,5\},\{1,2,3\}.$$

\noindent $\mathscr{F}_{T(G)^*}$ has at least five members, so $T(G)^*$ is non-transversal.

\chapter{Further topics of research}
In this chapter, we pose some questions for directions of further research.
\begin{question}
    Can we find other ways of ``randomly generating'' a matroid that provide more freedom than the method in Definition \ref{M(V)} ?
\end{question}
By having more ``freedom'', we mean a way of generating that is less restrictive than that in Definition \ref{M(V)}, in which we force every possible exchange of subsets (with the same cardinality) of $E_n$ we choose to be included. Certainly, randomly generated bipartite graphs (thus transversal matroids) provide more freedom than that, but that is still restricted to the class of transversal matroids. 

\begin{question}
    Can we find a necessary and sufficient condition for a bipartite graph $G$ such that $T(G)$ contains $U_{2,4}$ as a minor?
\end{question}
In general, if $\mathscr{A}$ is a minimal presentation of a matroid $M$, it is neither necessary nor sufficient for $G$ to contain $G_{\mathscr{A}}$ as a subgraph for $T(G)$ to have $M$ as a minor. Can we find such a necessary and sufficient condition for the case of $U_{2,4}$?  

\begin{question}
    If $\mathscr{A}$ is a presentation for a transversal matroid $M$, and $C\in \mathscr{C}(M)$, can we find a way to describe $C_{\mathscr{A}}?$
\end{question}
In Theorem \ref{thm:circuit, planting of pure basic type}, we show that in an acyclic presentation $\mathscr{A}$ of a transversal matroid $M$, the circuits of $M$ can be identified as successive plantings of pure basic type subgraphs of $G_{\mathscr{A}}$. If $M$ has no acyclic presentation, can we identify the circuits of $M$ in $G_{\mathscr{A}}$, where $\mathscr{A}$ is a presentation of $M$?

\begin{question}
    If we weaken the assumption of Theorem \ref{thm:admits acyclic pre.}, for instance, not every circuit $C$ of $M$ having nonempty intersection with $C_i$ intersects $C_i$ in a basic way, then we may not be able to obtain a presentation $\mathscr{A}'$ for $M$ such that $(C_1\cup C_2)_{\mathscr{A}'}$ is acyclic. But there may still be edges in the original presentation $\mathscr{A}$ for $M$, removing which still results in a presentation for $M$. Can we identify those edges?
\end{question}

\begin{question}
    A set $C\subseteq E(M)$ is a circuit of a matroid $M$ if and only if $E(M)\setminus C$ is a hyperplane of $M^*$ (see \cite[Theorem 2, p.39]{Welsh}). If $\mathscr{A}$ is a presentation for a transversal matroid $M$, and we study how to spot the circuits of $M$ via $G_{\mathscr{A}}$, then we get information on the hyperplanes of $M^*$ via $G_{\mathscr{A}}$. What can we know about the hyperplane arrangement of $M^*$ via $G_{\mathscr{A}}$? 
\end{question}

\clearpage
\phantomsection

\end{document}